\documentclass[aps,preprint,groupedaddress,showkeys]{revtex4-1}
\usepackage{graphicx}
\usepackage{dcolumn}
\usepackage{bm}
\usepackage{amssymb}
\usepackage{amsthm}
\usepackage{amsmath}
\usepackage{graphicx}
\usepackage{amsthm,amsmath,amssymb} \usepackage{mathrsfs}
\usepackage{epstopdf}

\newcommand{\rmd}{\textrm{d}}

\newtheorem{theorem}{\hskip\parindent\bf Theorem}
\newtheorem{lemma}{\hskip\parindent\bf Lemma}
\newtheorem{proposition}{\bf Proposition}
\newtheorem{remark}{\hskip\parindent\bf Remark}
\newtheorem{definition}{\hskip\parindent\bf Definition}
\begin{document}

\preprint{}

\title[Global dynamics in a predator-prey model with cooperative hunting and Allee effect and  bifurcation induced by  diffusion and  delays]{Global dynamics in a predator-prey model with cooperative hunting and Allee effect and  bifurcation induced by  diffusion and  delays}

\author{Yanfei Du}
\affiliation{School of Mathematics, Harbin Institute of Technology, Harbin, 150001, China.}

\affiliation{Shaanxi University of Science and Technology, Xi'an 710021,  China.}

\author{Ben Niu*}

\affiliation{Department of Mathematics, Harbin Institute of Technology, Weihai 264209,  China.\\*Corresponding author, niu@hit.edu.cn}

\author{Junjie Wei}

\affiliation{School of Mathematics, Harbin Institute of Technology, Harbin, 150001, China.}
\affiliation{School of Science, Jimei University,  Xiamen, Fujian,  361021, P. R. China}

\date{\today}

\begin{abstract}
We consider the local bifurcation and global dynamics of a predator-prey model with cooperative hunting and Allee effect. For the  model with weak cooperation,  we prove the existence of  limit cycle,  heteroclinic cycle at a threshold of conversion rate $p=p^{\#}$. When $p>p^{\#}$, both species go extinct, and when $p<p^{\#}$, there is a separatrix. The species with initial population above the separatrix finally become extinct; otherwise, they coexist or oscillate sustainably.    In the case  with strong cooperation,  we exhibit the complex dynamics of system in three different cases, including limit cycle, loop of heteroclinic orbits among three equilibria, and homoclinic cycle.  Moreover, we find diffusion may induce Turing instability and Turing-Hopf bifurcation, leaving the system with  spatially inhomogeneous distribution of the species, coexistence of two different spatial-temporal oscillations.  Finally, we investigate Hopf and double Hopf bifurcations of the diffusive system induced by two delays.
\end{abstract}

\keywords{Hunting cooperation, Allee effect, connecting orbit, invariant manifold,  bifurcation, coexistence}
                           \maketitle

\section{Introduction}
\label{}

A  predator-prey model can be described by the following equations \cite{Cosner}
\begin{equation}
\label{ode1}
     \left\lbrace            \begin{array}{lll}
\dot{u}  &=& f(u)u   -G(u,v)v, \\
\dot{v}  &=&     dG(u,v)v-mv  , \\
\end{array}
\right.\end{equation}
where $u$ and $v$ stand for the densities of prey and predator, respectively. $f(u)$ represents
the per capita prey growth rate in the absence of
predators, and the function $G(u,v)$ is the functional response charactering predation. $d$ describes the rate of biomass conversion from predation, and $m$ is the death rate of predator.   Many kinds of functional response have been proposed \cite{Jeschke}.  Among them Holling type I, II, III \cite{C. S. Holling,N. Kazarinov,P. Turchin} are discussed widely, which are prey-dependent functional response.

Cooperative behavior within a  species  is  very common in nature, of which cooperative hunting is the
most widely distributed form  in animals. In mammals, the most famous examples include wolves \cite{Schmidt}, African wild dogs \cite{Courchamp}, lions \cite{Scheel},  and chimpanzees \cite{C. Boesch}, who  cooperate for hunting preys.   Cooperative behaviors are also widespread in  other living organisms,   such as   aquatic organisms \cite{Bshary}, spiders \cite{Uetz}, birds \cite{D. P. Hector},  and ants \cite{Moffett},  who find, attack and move their preys together. Due to hunting cooperation, the attack rate of predators  increases with predator density, which   obviously   depends on both prey and predator densities.  Recently,   a few literatures have paid attention to derive functional response to describe  the  cooperative hunting  \cite{Cosner,Berec,Alves}.
Cosner et al.  \cite{Cosner} proposed a
functional response to explore the effects of predator aggregation when
predators encounter a cluster of prey.
Berec \cite{Berec}  generalized the Holling type-II functional
response to  a family of functional responses by considering attack rate and handling time
of predators varies with predator density  to interpret the foraging facilitation among predators.
 Alves and Hilker  \cite{Alves} considered  a special
 case of the more general functional response  proposed by Berec \cite{Berec}.
They added a cooperation term to the attack rate for representing the benefits that hunting
cooperation brings to the predator population, which makes the functional response $G(u,v)=(b+cv)u$. Here $c > 0$  describes the strength of predator cooperation in hunting, where $cy$ is referred to as the cooperation term.
They found that cooperative hunting can
improve persistence of the predator,
  and it is a form of
 foraging facilitation which can induce strong Allee effects.
After their work, many predator-prey model with cooperative hunting in the predator have been investigated \cite{Jang,S. Pal,D. Sen,S. Yan,Wu D,Danxia Song}.

The most well-known
example of $f(u)$ is the logistic form
  $f(u)= r  \left(  1-\frac{u }{K}\right)$,  where  $r$ is the intrinsic growth rate, $K$ denotes the environmental carrying capacity.  Then, in the absence of predator,   the prey population is governed by   $\dot{u} = ru \left(  1-\frac{u }{K}\right)$. We can conclude from the equation that the species may increase in size   when the density  is low.  However, for many species,  low population density may induce  many problems, such as mate difficulties, inbreeding and predator avoidance of defence.   It turns out  that the growth rate of the low density population is not always positive, and it may be negative when the density of population is less than  the  minimum number $a$ for the survival of the population, which is called the Allee threshold.
The phenomenon is known as Allee effect, which was first observed by Allee \cite{Allee}, and  has been observed on various organisms, such as vertebrates, invertebrates and plants.  Such a population can be described by   $\dot{u} = ru \left(  1-\frac{u }{K}\right)(u-a)$.  Allee effect may cause the extinction of low-density population. Because  of  increasing fragmentation of habitats, invasions of exotic species,
Allee effect has aroused more and more attention (see for example \cite{Courchamp F,Wang2010,Conway} and the references cited therein).

Jang  et al. \cite{Jang} considered both the cooperation in the predator and Allee  effects in the prey, and proposed the following   predator-prey model
\begin{equation}
\label{odepredator}
     \left\lbrace            \begin{array}{lll}
\dfrac{{\rm d}u}{{\rm d}t} &=& r_1u\left[  1-\frac{u}{K_1}\right]\left[u-a_1\right] -\left[ b_1+c_1v\right] uv, \\
\dfrac{{\rm d}v}{{\rm d}t} &=&     p_1[b_1+c_1v]uv-m_1v , \\
\end{array}
\right.\end{equation}
where $r_1$ is the per capita intrinsic growth rate of the prey, $K_1 $ is the enviromental carrying capacity for the prey, $a_1$ ($a_1<K_1$) is the Allee threshold of the prey population, $b_1$ is the attack rate per predator and prey, $c_1$ measures the degree of cooperation of predator, $p_1$ is the prey conversion to predator, and  $m_1$ is the per capita death rate of predator. They discussed the extinction and and coexistence of the species when the parameters are in different ranges, and found  that cooperation hunting might change the number of equilibria and change that stability of the interior equilibria. They also devised a best strategy for culling the predator by maximizing the prey population and minimizing the predators along with the costs associated with the control.

With a nondimensionalized change of variables:
$$\widehat{u}=\frac{1}{K}u,~\widehat{v}=bv,$$
and dropping the hats for simplicity of notations, system (\ref{odepredator})   takes the form
\begin{equation*}
\left\{
\begin{array}{l}
\dfrac{{\rm d}u}{{\rm d}t}= K_1r_1u[1-u] (u-\frac{a_1} {K_1}) -(1+\frac{c_1} {b_1^2}v)uv,\\
\dfrac{{\rm d}v}{{\rm d}t}= m_1v\left[  \frac{b_1p_1K_1} {m_1} u\left(  1+\frac{c_1} {b_1^2}v\right) -1 \right]  .
\end{array}
\right.
\end{equation*}
Let
$$
r=K_1r_1,~~  a=\frac{a_1}{K_1},~~  c=\frac{c_1} {b_1^2},~~\mbox{and}~~ p=\frac{b_1p_1K_1} {m_1}.
$$
Then we obtain the simplified dimensionless system
\begin{equation}
\label{ODE system}
\left\{
\begin{array}{l}
\dfrac{{\rm d}u}{{\rm d}t}=  ru(1-u) (u-a) -(1+cv)uv, \\
\dfrac{{\rm d}v}{{\rm d}t}=  mv\left[  p u\left( 1+cv\right) -1 \right] .\\

\end{array}
\right.
\end{equation}
Recall the definition of Allee effect, one always has $0<a<1$.
In fact, for fixed $K_1$, $r_1$, $a_1$, $b_1$ and $m_1$, we have $c$ and $p$   directly proportional to $c_1$ and $p_1$, respectively. Thus, we refer to $c$ and $p$ as the degree of cooperative hunting of the predator and the conversion rate  from the prey to the predator  respectively in the context. By the explanation above, we make the following assumption
always:
$$
(H)~~~~~~r,~~a,~~c,~~m,~~p~\mbox{and}~c~~\mbox{are all positive, and}~~0<a<1.
$$

Motivated by \cite{Jang}, in this paper, we will consider a complete global analysis of the model (\ref{ODE system}) by investigating the stable/unstable manifolds of saddles, then the existence of some connecting orbits is obtained. Hence, a partition of the phase space is given, i.e., the attraction basin of each  equilibrium is obtained.


 This paper is organized as follows. In section \ref{sectionode}, we investigate the global dynamics of the ODE predator-prey model (\ref{ODE system}) with  weak cooperation and strong cooperation respectively. Firstly, for the case of weak cooperation, the existence and stability of equilibria are discussed.  Taking $p$ as the bifurcation parameter,  the existence of Hopf bifurcation  and loop of  heteroclinic orbits is proved,  and the global dynamics are investigated.  For the case with strong cooperation, the existence of Hopf bifurcation, loop of heteroclinic orbits, and homoclinic cycle  are observed by theoretical analysis or numerical simulation.
In  section \ref{sectionTuring}, we consider the diffusive system (\ref{diffusion}), and investigate Turing instability and Turing-Hopf bifurcation induced by diffusion. We illustrate some complex dynamics of system, including the existence of spatial inhomogeneous steady state, coexistence of two  spatial inhomogeneous periodic solutions.
In  section \ref{sectiondelays}, we consider the diffusive system  with two delays, and Hopf bifurcation and double Hopf bifurcation induced by two delays are analyzed on the center manifold via normal form approach.

\section{Global  dynamics of the ODE system}\label{sectionode}
In this section, we consider the ODE system (\ref{ODE system}) with two different cases. It is found that the strength of cooperation heavily affects the number of interior equilibria. Thus, the theoretical results are given on two cases:   when  $c<\frac{1}{r(1-a)}$, we say the cooperation among predators is weak;  when  $c>\frac{1}{r(1-a)}$, the cooperation  is strong.  We investigate the local and global dynamics in both cases, together with some numerical illustrations, as well as biological interpretations.

\subsection{The  system  with weak cooperative hunting}\label{odedynamicsweak}

 In this section, we first consider the existence of boundary equilibria and interior equilibria, respectively. Then, taking $p$ as bifurcation parameter, we investigate the Hopf  bifurcation near the unique interior equilibrium.  Through studying the stable manifold and unstable manifold of saddles,   we prove the existence of loop of heteroclinic orbits, and   get the global dynamics of system (\ref{ODE system}).

For model (\ref{ODE system}), the first quadrant is  invariant since $\{(u,v):u=0\}$ and $\{(u,v):v=0\}$  are invariant manifolds for (\ref{ODE system}). We can get the following result.
\begin{lemma}\label{bounded}
	The solution of (\ref{ODE system}) with positive initial value is positive and bounded.
\end{lemma}
\noindent Proof.  For any $u(0)>1$, $u'=ru(1-u)(u-a)-(1+cv)uv<0$ if $u>1$. On $u=1$, $u'<-(1+cv)v<0$. Noticing  that there is  no equilibrium in the region $\{(u,v):u>1,v\geq0\}$,  any positive solution satisfies $u(t)\leq\max\{u(0),1\}$ for $t\geq0$. From (\ref{ODE system}), we obtain that $(mpu+v)'=mpru(1-u)(u-a)-mv\leq mpru(1-u)(u-a)+m^2pu-m(mpu+v)\leq \zeta -m(mpu+v)$, where  $\zeta=\max_{t\geq0}\{mpru(1-u)(u-a)+m^2pu\}$.  Then we have
$mpu(t)+v(t)\leq (mpu(0)+v(0))e^{-mt}+\frac{\zeta}{m}(1-e^{-mt})$, which means that $v(t)$ is bounded.

\subsubsection{Existence  of equilibria}

  System (\ref{ODE system}) has three boundary equilibria:  (i)  $E_0(0,0)$, which means  the extinction of both the species;  (ii)  $E_a(a,0)$, which is induced by Allee effect; (iii) $E_1(1,0)$,  which means the extinction of the predator and the survival of  the prey, achieving its carrying capacity.

  Now we discuss the existence of interior equilibria similar as the discussion in \cite{Jang}. The existence of interior equilibria depends on the position of $u-$nullcline $r(1-u)(u-a)=(1+cv)v$ and  $v-$nullcline $pu(1+cv)=1$.
 In fact, $r(1-u)(u-a)=(1+cv)v$ is an ellipse, sitting in the first and forth quadrant,  and intersecting the horizontal axis at $a$  and $1$.  $pu(1+cv)=1$ is hyperbolic, with its right branch sitting in  the first and forth quadrant and intersecting the horizontal axis at $\frac{1}{p}$, and is decreasing and  concave in $(0,\frac{1}{p})$.
Denote the $u-$nullcline curve $r(1-u)(u-a)=(1+cv)v$ in the first quadrant as $v=f(u)$,     the $v-$nullcline curve $pu(1+cv)=1$ in the first quadrant as $v=g(u)$.
Obviously, an intersection  of $v=f(u)$ and $v=g(u)$ is an interior equilibrium, thus any interior equilibrium has components $u>0$ and $v>0$ satisfying
\begin{equation}\label{equili3}
\begin{array}{rl}
r(1-u)(u-a)&=(1+cv)v,\\
pu(1+cv)&=1,
\end{array}
 \end{equation}
i.e., \begin{equation}\label{equili}
\begin{array}{rl}
r(1-u)(u-a)&=\frac{1-pu}{cp^2u^2},\\
v&=\frac{1-pu}{cpu}.
\end{array}
 \end{equation}

The existence of interior equilibria may have the following different cases.

\begin{figure}
  	\centering
  	a) \includegraphics[width=0.45\textwidth]{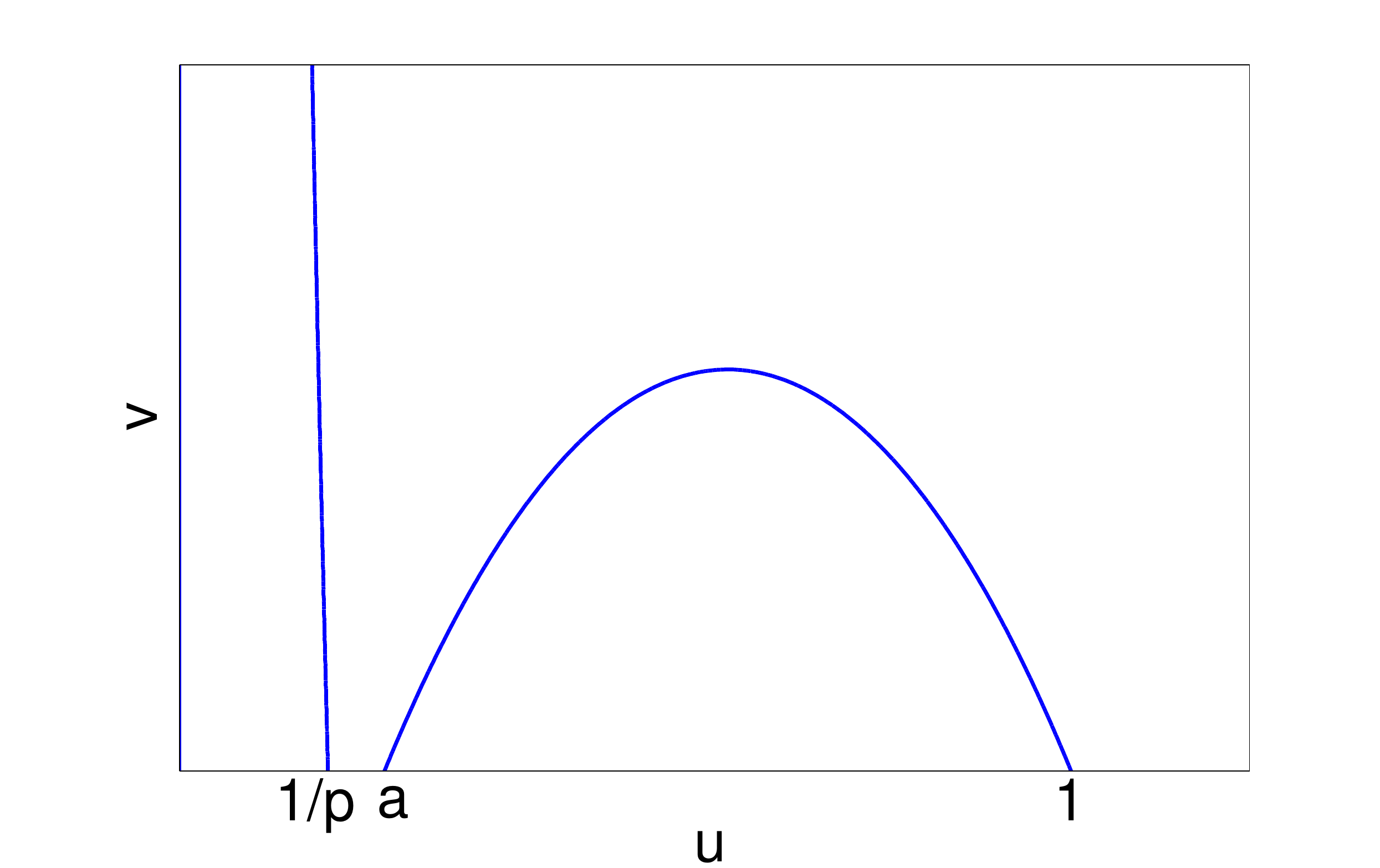}
    	b)   \centering
  	\includegraphics[width=0.45\textwidth]{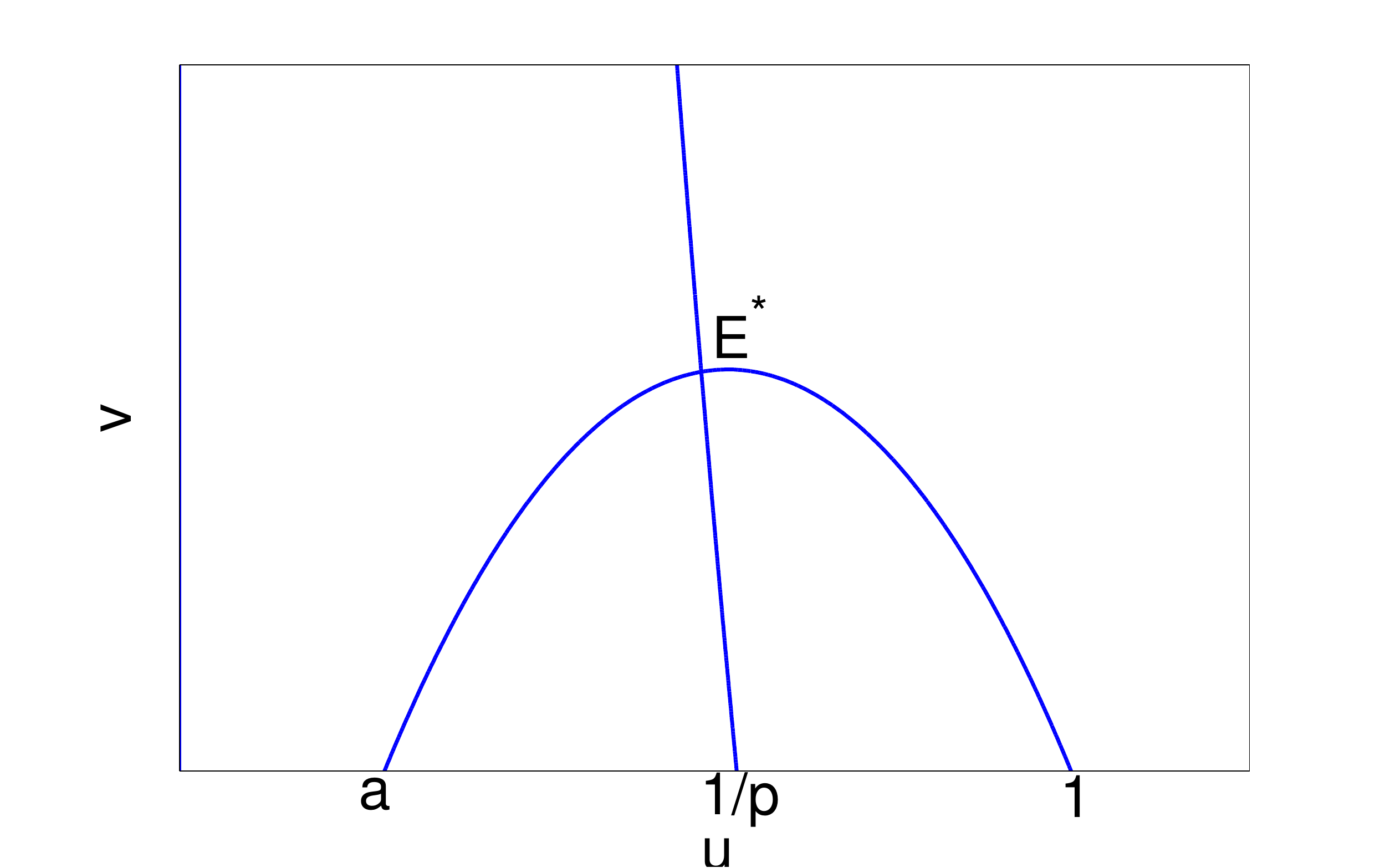}
  	\caption{
  	a) When  $p>\frac{1}{a}$,  there is  no interior equilibrium. b) When $1<p<\frac{1}{a}$,   there is a unique interior equilibrium.
  	}
  	\label{fig:akequilibrium}
  \end{figure}

\begin{proposition}\label{numberofequilibrium}Suppose that $c<\frac{1}{r(1-a)}$.\\
 $~~~~~~$	(i)  When $p\geq\frac{1}{a}$, system (\ref{ODE system})  has  no interior equilibria.\\
$~~~~~~$	(ii) When $1<p<\frac{1}{a}$,  system (\ref{ODE system}) has a unique positive constant equilibrium $E^*(u^*,v^*)$ with $a<u^*<\frac{1}{p}$.\\
$~~~~~~$	(iii) When $p\leq 1$, system (\ref{ODE system})  has  no interior equilibria.%
\end{proposition}
\noindent Proof.
 In fact, when $a<\frac{1}{p}<1$, two nullclines have two intersections, $E^*$ in the first quadrant (see Fig. \ref{fig:akequilibrium} b)) and $E_R^*$ in the forth quadrant. Thus there is a unique interior equilibrium $E^*(u^*,v^*)$ with $u^*$ and $v^*$ satisfying (\ref{equili3}).  When $\frac{1}{p}$ decreases to $a$, $E^*$ collides with $E_a$, and  $E_R^*$ is still in the forth quadrant. When $\frac{1}{p}<a$, $E^*$ moves into the forth quadrant, leaving no interior equilibria (see Fig. \ref{fig:akequilibrium} a)).

\begin{figure}
      	\centering
      	a) \includegraphics[width=0.45\textwidth]{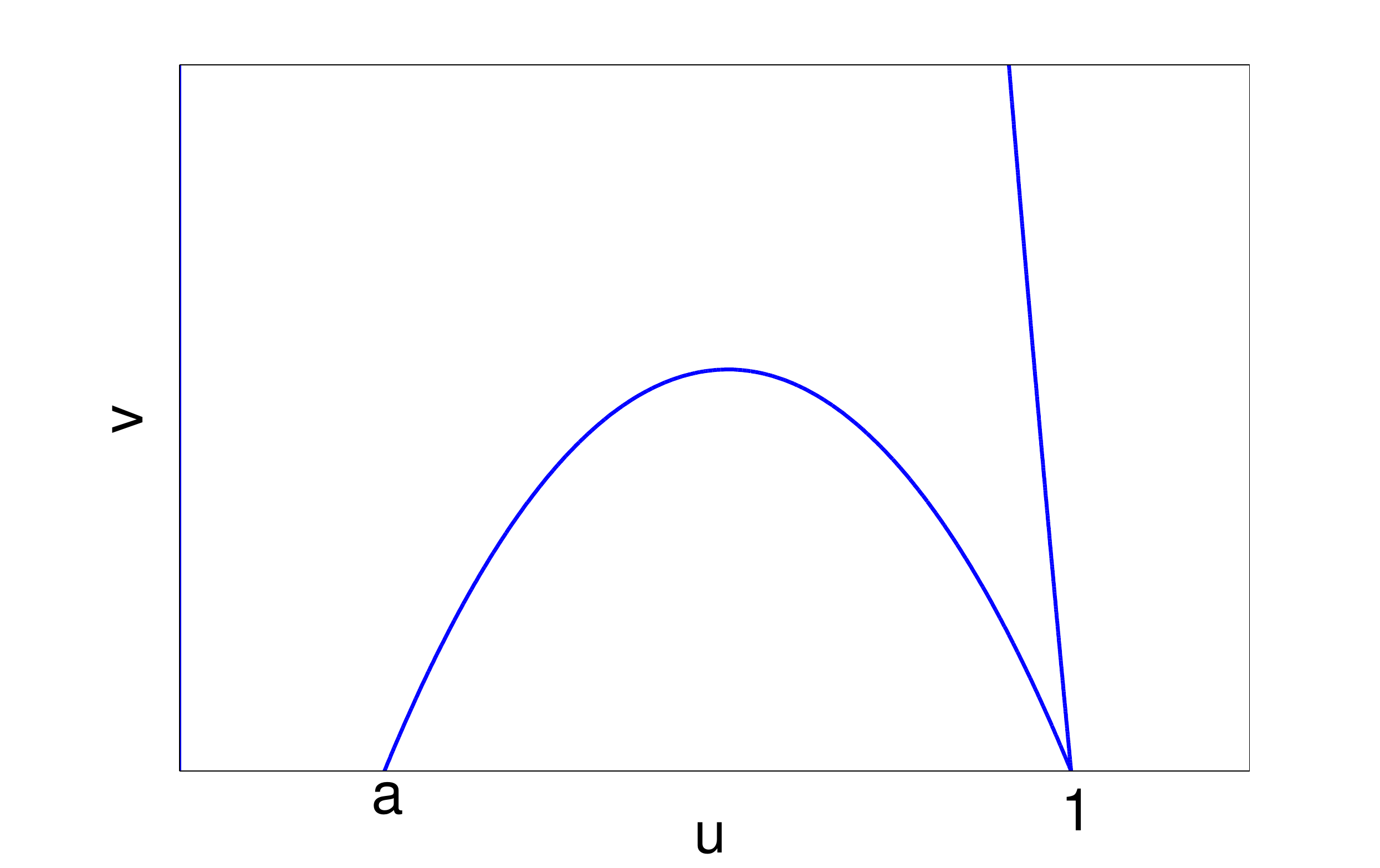}
                   	b)   \centering
             	\includegraphics[width=0.45\textwidth]{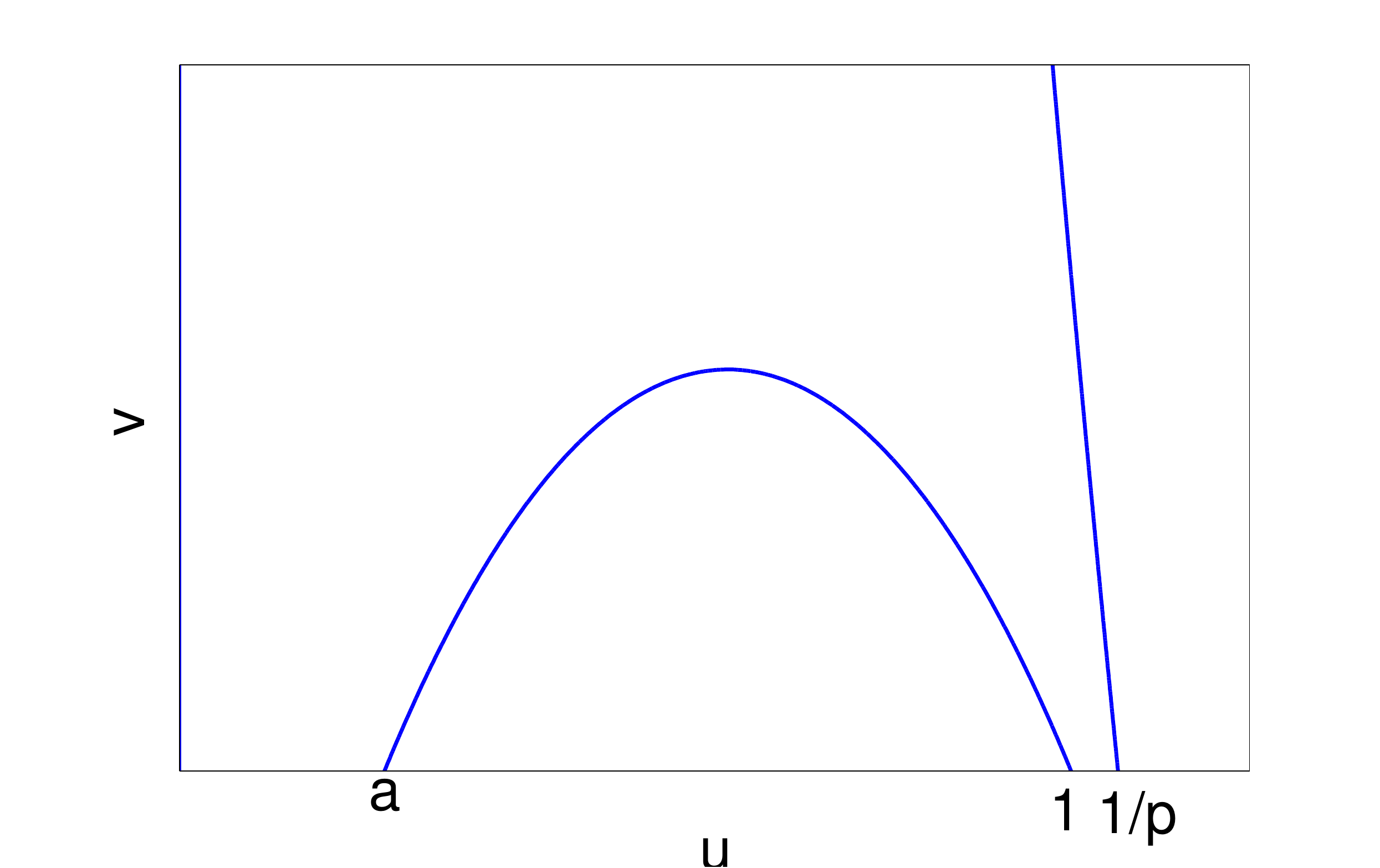}
      	\caption{	
      		If $c<\dfrac{1}{r(1-a)}$,    there is  a) no interior equilibrium when $ p=1$;  b) no interior  equilibrium when $p<1$.}
      	\label{fig:pequals1}
      \end{figure}

  When $\frac{1}{p}$ increases to $1$, the $v-$nullcline $v=g(u)$ and $u-$nullcline  $v=f(u)$  intersect at $E_1(1,0)$.  The slope of tangents of $v-$nullcline and $u-$nullcline at $E_1(1,0)$ are $-\frac{1}{c}$ and $-r(1-a)$ respectively.
 If $-\frac{1}{c}< -r(1-a)$,  $E^*$ collides with $E_1$ when $p=1$ (see Fig. \ref{fig:pequals1} a)),  $E_R^*$ is in the forth quadrant, and there is no interior equilibria. Moreover, when  $p<1$, there is no  interior equilibria (see Fig. \ref{fig:pequals1} b)). $~~~~~~~~~~~~~~~~~~~~~~~~~~~~~~~~~~~~~~~~~~~~~~~~~~~~~~~~~~~~~~\Box$

 There is always a unique positive equilibrium $E^*(u^*,v^*)$  when  $1<p<\frac{1}{a}$.   Now we wonder  the impact of $c$ on the value of  the component of the unique equilibrium $E^*(u^*,v^*)$  when  $1<p<\frac{1}{a}$, and we have the following conclusion.
 \begin{proposition} \label{ucmontone}
 	When $1<p<\frac{1}{a}$, for	the unique positive  equilibrium $E^*(u^*,v^*)$ ($a<u^*<\frac{1}{p}$),  $u^*$ is monotonically   decreasing with respect to $c$;  $v^*$ is monotonically  decreasing with respect to $c$ when $a<u^*<\frac{a+1}{2}$, and  it is  monotonically increasing when $\frac{a+1}{2}<u^*<1$.
 \end{proposition}
 \noindent  Proof.
 Consider both sides of the first equation of (\ref{equili}) as functions of $u$, denoted by  $y=\frac{1-pu}{cp^2u^2}$ and $y=r(1-u)(u-a)$.  Obviously,  The $u-$component of   interaction of $y= r(1-u)(u-a)$ and  $y=\frac{1-pu}{cp^2u^2}$  moves left when the value of $c$ increases.   It means that   $u^*$ is monotonically  decreasing with respect to $c$.

 Now we focus on the effect of $c$ on the     $v^*$ component of the interior equilibrium.   Since    $r(1-u^*)(u^*-a)=\frac{1-pu^*}{cp^2u^{*2}}
 $, and
 $v^*=\frac{1-pu^*}{cpu^*}$, then $v^*= rpu^*(1-u^*)(u^*-a)$.  With the increasing of $c$, $u^*$ decreases. If $\frac{a+1}{2}>u^*>a$,
 $v^*$ is increasing with respect to $u^*$, thus it is decreasing with $c$.
 If $\frac{a+1}{2}<u^*<1$, $v^*$ is decreasing with $u^*$, and thus it is increasing with $c$.  $~~~~~~~~~~~~~~~~~~~~~~~~~~~~~~~~~~~~~~~~~~~~~~~~~~~~~~~~~~~~~~~~~~~~~~~~~~\Box$

 \begin{remark}\label{conequilibrium}
 	In the absence of cooperative hunting within the predator, i.e., $c=0$,   system (\ref{ODE system}) always has a unique interior equilibrium $E_0^*=(\frac{1}{p}, r(1-\frac{1}{p})(\frac{1}{p}-a))$,  and $E_0^*$ exists if and only if  $1<p<\frac{1}{a}$.   The cooperative hunting  $c$  may change the density of prey and predator. It is not surprising that  the increasing of cooperative hunting will decrease the density of the prey.   When $\frac{a+1}{2}<u^*<1$,   cooperative hunting is beneficial to  the density of predator.
 	However, when $a<u^*<\frac{a+1}{2}$,  the  cooperative hunting will  decrease the  stationary  density of the predator.
 \end{remark}

  \subsubsection{Stability of all equilibria and Hopf bifurcation at $E^*$}\label{sectionHopf}

 The Jacobian matrices of the function on the right-hand of  system (\ref{ODE system}) around $E_0$, $E_a$, and $E_1$ are respectively,
 \begin{equation*}\label{Jacoode0}
 J_{E_0}=\left( \begin{array}{cc}
 -ra & 0\\0 & -m
 \end{array}\right),
 \end{equation*}
 \begin{equation*}\label{Jacoodea}
 J_{E_a}=\left( \begin{array}{cc}
 ra(1-a) &-a\\0 &  m(pa-1)
 \end{array}\right),
 \end{equation*}
 \begin{equation*}\label{Jacoode1}
 J_{E_1}=\left( \begin{array}{cc}
 -r(1-a) &-1\\0 &  m(p-1)
 \end{array}\right),
 \end{equation*}
 from which we can easily get the local stability of the boundary equilibria.

 \begin{lemma}\label{boundary}
   For system (\ref{ODE system}),
 \\$~~~~~~$(i) $E_0(0,0)$ is a stable node;
 \\$~~~~~~$(ii) $E_a(a,0)$ is an unstable node  if $p>\frac{1}{a}$, and it is a saddle if $p<\frac{1}{a}$;
 \\$~~~~~~$(iii) $E_1(1,0)$ is a stable node if $p<1$,    and it is a saddle if $p>1$.
 \end{lemma}

   For the interior equilibrium $E^*(u^*,v^*)$, the corresponding Jacobian matrix is
   \begin{equation}\label{Jacoode}
     J_{E^*}=\left( \begin{array}{cc}
     ru^*(1+a-2u^*)&-2cu^*v^*-u^*\\mpv^*(1+cv^*) &mpcu^*v^*
     \end{array}\right),
     \end{equation}
     and thus
     \begin{equation*}
     \begin{array}{l}
          {\rm tr} J_{E^*}=-ru^*(2u^*-a-1)+mpcu^*v^*=-ru^*(2u^*-a-1)+m(1-pu^*),\\
      {\rm det} J_{E^*}=mpu^*v^*\left[ rcu^*(1+a-2u^*)+(1+cv^*)(1+2cv^*)\right]  .
      \end{array}
    \end{equation*}

   In fact,   $u-$nullcline $v=f(u)$ and $v-$nullcline $v=g(u)$ intersect at $E^*$,  where the slope of $v=f(u)$ is larger than that of $v=g(u)$, i.e.,  $\dfrac{1+cv^{*}}{-cu^{*}}< \dfrac{r(1+a-2u^{*})}{1+2cv^{*}}$.  Thus,   ${\rm det} J_{E^*}>0$.
   Thus $E^*$ may be node or focus,  and its stability depends on  the trace \begin{equation}{\rm tr}J_{E^*}=-ru^*(2u^*-a-1)+m(1-pu^*).\end{equation}
    Thus, we have the following conclusions on the stability of $E^*$.

\begin{figure}
    	\centering
      \includegraphics[width=0.5\textwidth]{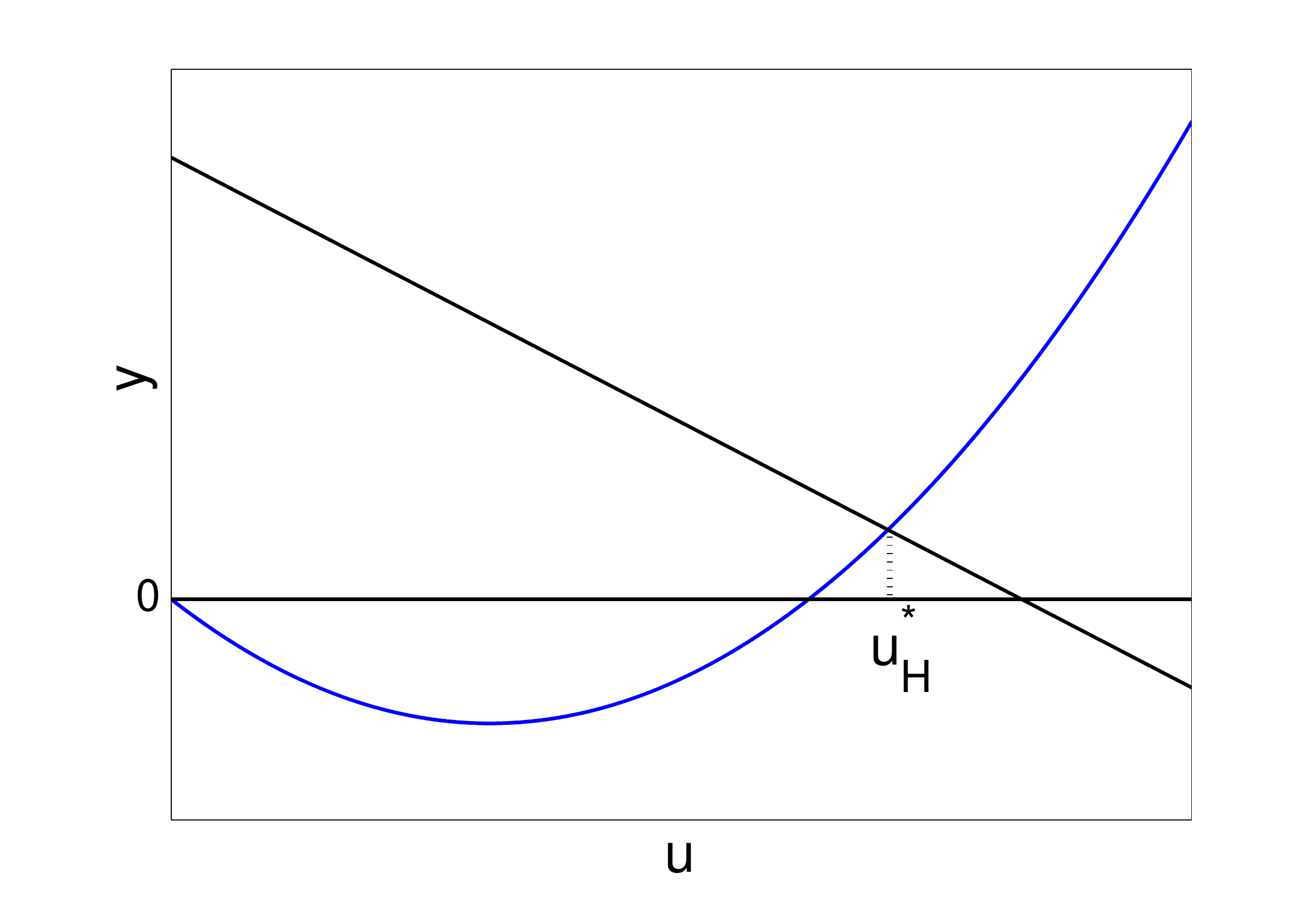}

    	\caption{
    	  The blue curve represents  $y=ru^*(2u^*-a-1)$, and the black  line represents $y=m(1-pu^*)$. When $u=u_H^*$,  ${\rm tr}J_{E^*}=0$.
    	}
    	\label{fig:conExingtr}
    \end{figure}

 \begin{theorem}\label{HopfExing}
	If $c<\frac{1}{r(1-a)}$, then there exists a unique $p_H\in (1,\frac{2}{a+1})$ such that the unique interior equilibrium $E^*$ is locally asymptotically stable when $1<p<p_H$, and unstable when  $p_H<p<\frac{1}{a}$. Moreover,
 	system (\ref{ODE system}) undergoes a Hopf bifurcation at $E^*$ when  $p=p_H$.
 \end{theorem}
 \noindent Proof.  Solving for ${\rm tr}J_{E^*}=0$,  when $\frac{a+1}{2}<\frac{1}{p}<1$, we obtain a unique $u_{H}^*\in (\frac{a+1}{2},1)$,  such that ${\rm tr}J_{E^*}>0$ on
              $(0,u_{H}^*)$, and ${\rm tr}J_{E^*}<0$ on  $(u_{H}^*,1)$ (see Fig. \ref{fig:conExingtr}).
          Moreover, it is easy to verify that $\frac{\partial u^*}{\partial p}<0$,  thus, corresponding to $u_{H}^*$,   we obtain a unique $p_H$ ($1<p_H<\frac{2}{a+1}$),    such that ${\rm tr}J_{E^*}>0$ when  $p\in( p_H,\frac{1}{a})$, and ${\rm tr}J_{E^*}<0$ when  $p\in(1,p_H)$.

         Now we verify the transversality condition. Let $\mu=\alpha(p)\pm i \omega(p)$ be the roots of  $\mu^2-{\rm tr}J_{E^*}\mu+{\rm det}J_{E^*}=0$ when $p$ is near $p_H$. We have $\alpha'(p)=\frac{1}{2}\frac{{\rm d}{\rm tr}J_{E^*}}{{\rm d}p}=\frac{1}{2} [ -ru^{*'}(p)(2u^*-a-1)-2ru^*u^{*'}(p)$$ -mu^*-mpu^{*'}(p)]$ , where $ u^{*}(p)$  is a function of $p$ determined by  (\ref{equili}).
  Taking the derivative of both sides of (\ref{equili}) with respect to $p$, we have $$-ru^{*'}(p)(2u^*-a-1)=\dfrac{(-u-pu^{*'}(p))[cp^2u^{*2}+2cpu^*(1-pu^*)]}{c^2p^4u^{*4}},$$ thus $$-u-pu^{*'}(p)\mid_{p=p_{H}}=\dfrac{-ru^{*'}(p)(2u^*-a-1)c^2p^4u^{*4}}{cp^2u^{*2}+2cpu^*(1-pu^*)}\mid_{p=p_{H}}>0,$$ since $\frac{a+1}{2}<u_H^*<\frac{1}{p}$  and  $u^{*'}(p)<0$.
  Thus,  $\alpha'(p)\mid_{p=p_{H}}>0$.
$~~~~~~~~~~~~~~~~~~~~~~~~~~~~~~~~~~~~~~~~~~~~~~\Box$

 \begin{remark}\label{ptop}
  When there is no cooperative hunting in predator, i.e., $c=0$,
  the $v-$nullcline is $u=\frac{1}{p}$, which is vertical.  $
       {\rm tr }J_{E^*}\arrowvert_{c=0}=
       -\frac{r}{p}(\frac{2}{p}-a-1),
     {\rm det}J_{E^*}\arrowvert_{c=0}=       mv^*.       $
       The system  undergoes a Hopf bifurcation near $E^*$  when $\frac{1}{p}=\frac{a+1}{2}$, which is the top of     the $u-$nullcline   $v=r(1-u)(u-a)$. In fact, it has been discussed widely that the stability of positive equilibrium can be stated graphically by the  $u-$nullcline when the $v-$nullcline is vertical: $(u^*,v^*)$ is unstable if the
$v-$nullclines intersect to the left of a local maximum of the $u-$nullcline, and
stable if they intersect to the right    \cite{Rosenzweig,Rosenzweig1969,Hsu SB,Wang2010},
  and thus, Hopf bifurcation occurs at the ``top of the hump" of   $u-$nullcline.
   However, for system (\ref{ODE system})  in this paper, the  $v-$nullcline is not vertical, and we can verify that  the intersection of the $u$-nullcline and the $v$-nullcline corresponding to  the Hopf bifurcation point $p=p_H$ is on the right of the ``top of the hump" of $u$-nullcline.      In fact, the ``top of the hump" of   $u-$nullcline $v=f(u)$ achieves at $u=\frac{a+1}{2}$, at which the corresponding  value of $p$ is denoted by $p_{{\rm top}}$.
  	     When $p=p_{\rm top}$, $E^*$ can be also proved to be unstable.    Comparing the results in systems with and without cooperation, we can conclude  that  cooperative hunting is more likely to  bring  instability into $E^*$.   
  \end{remark}

In order to determine the direction of Hopf bifurcation and the stability of the Hopf bifurcating periodic solution, we need to calculate the normal form near the Hopf bifurcation point.  Through direct calculation following the steps in \cite{Wiggins S }, we can get the following truncated normal form
 \begin{equation*}
 \begin{array}{l}
\dot{r}=\alpha'(p_H)pr+a(p_H)r^3+O(r|p-p_H|^2,r^3|p-p_H|,r^5),\\
\dot{\theta}=\omega(p_H)+\omega'(p_H)(p-p_H)+c(p_H)r^2+O(|p-p_H|^2,r^2|p-p_H|,r^4).
 \end{array}
 \end{equation*}

   Recalling that    $\alpha'(p_H)>0$,     the direction of Hopf bifurcation and stability of Hopf bifurcating periodic solution are determined by the first Lyapunov coefficient $a(p_H)$, and we have the following conclusion.

 \begin{theorem}\label{Hopfdirection}
 System  (\ref{ODE system}) undergoes a Hopf bifurcation at $E^*$ when $p=p_H$.\\ $~~~~~~$ (i) If $a(p_H)<0$,  the bifurcating periodic solution is orbitally asymptotically stable, and it is bifurcating from  $E^*$ as $p$ increases and passes $p_H$.  \\
$~~~~~~$ (ii) If $a(p_H)>0$,  the bifurcating periodic solution is unstable, and it is bifurcating from  $E^*$ as $p$ decreases and passes  $p_H$.
 \end{theorem}

\subsubsection{The global dynamics of system (\ref{ODE system}) with weak cooperative hunting}

From the previous discussion, we have known  the stability of all equilibria when $c<\frac{1}{r(1-a)}$,  which is  listed in  Table 1. Motivated by the work of \cite{Wang2010,Conway}, we investigate the  global dynamics of system (\ref{ODE system}) when $c<\frac{1}{r(1-a)}$.

 \begin{table}[tbp]\label{Jacobian1}
 \caption{Stability of  equilibria of (\ref{ODE system})  when $c<\frac{1}{r(1-a)}$.}\centering	\begin{tabular}{ccccc}\hline	Equilibrium  & $p>\frac{1}{a}$   & $p_H<p<\frac{1}{a}$ & $1<p<p_H$ &  $p<1$ \\   \hline
 $E_0 (0,0)$ &stable node &stable node& stable node & stable node \\         

 $E_a (a,0)$  & unstable node & saddle & saddle & saddle  \\        

 $E_1(1,0)$  &
saddle & saddle & saddle& stable node \\

  $E^*(u^*,v^*)$  &
 does not exist & unstable node or focus & stable  node or focus&  does not exist \\ \hline
  \end{tabular} \end{table}

When $p<\frac{1}{a}$, $E_a$ is a saddle. The eigenvector corresponding to $\lambda_1=ra(1-a)>0 $  is $(1,0)$, which means that  the unstable manifold of $E_a$ is on the $u$-axis.  For $\lambda_2=m(pa-1)<0$,  the  corresponding eigenvector  is  $(1,r(1-a)-\frac{m(pa-1)}{a})^T$.  Thus the tangent vector of the stable manifold of $E_a(a,0)$ (denoted by $\Gamma_p^s$)  at $E_a(a,0)$ is  $k_1=r(1-a)-\frac{m(pa-1)}{a}$.  Comparing  with the tangent vector $k_2=r(1-a)$ of $v=f(u)$  at $E_a(a,0)$, we have $k_1>k_2$, thus $\Gamma_p^s$  is above the nullcline $v=f(u)$ near $E_a$.   Moreover, from the vector field in (\ref{ODE system}) on $v=f(u)$,  $\Gamma_p^s$ is always above $v=f(u)$ before $\Gamma_p^s $ meets the $v-$nullcline.

Similarly, when $p>1$, $E_1$ is a saddle. Its stable manifold is on the $u-$axis, and its unstable manifold, $\Gamma_p^u$, is above the nullcline $v=f(u)$  before it meets the $v-$nullcline.

 When $1<p<\frac{1}{a}$, $E_0$ and $E_a$ are all saddle points. To figure out the global dynamics of system (\ref{ODE system}), we should first study the stable manifold $\Gamma_p^s$ of $E_a$ and unstable manifold  $\Gamma_p^u$ of $E_1$.

 \begin{proposition}\label{manifold}
 	Suppose that $c<\frac{1}{r(1-a)}$, and $1<p<\frac{1}{a}$.   	\\
 $~~~~~~$	(i) The orbit  $\Gamma_p^s$ meets the $v-$nullcline $v=g(u)$ at a point $(u_p^s,S(p))$, where $S(p)\geq v^*:=f(u^*)$, and $S(p)$ is a monotone decreasing function for $p\in (1,\frac{1}{a})$.\\
 $~~~~~~$ (ii) The orbit  $\Gamma_p^u$ meets the $v-$nullcline $v=g(u)$ at a point $(u_p^u,U(p))$, where $U(p)\geq v^*:=f(u^*)$, and $U(p)$ is a monotone increasing function for $p\in (1,\frac{1}{a})$.
  		\end{proposition}
 	\noindent Proof. Inspired by \cite{Wang2010}, we prove (i),  and the proof for (ii) is similar.  		
 	We have proved that   $\Gamma_p^s$ approaches $E_a$ from the region $\{(u,v):v>f(u)\}$. Moreover,  $\Gamma_p^s$ is always above $v=f(u)$ before $\Gamma_p^s $ meets the $v-$nullcline. 	
 	In order to show  that $\Gamma_p^s$ meets the $v-$nullcline, we still need to prove that  it  remains bounded  for $u>a$. In fact,  for $u>a$, $\Gamma_p^s$ is the graph of a function $v(u)$,  satisfying $$\dfrac{{\rm d}v(u)}{{\rmd}u}=\dfrac{mv[1-pu(1+cv)]}{u[(1+cv)v-r(1-u)(u-a)]}.$$ If there is a $u_b<u^*$ such that  $v(u)\rightarrow \infty$ as $u\rightarrow u_b^-$, then $(1+cv)v-r(1-u)(u-a)$ is bounded below for $a+\varepsilon\leq u\leq u_b$ and any $\varepsilon>0$. Thus for $a+\varepsilon\leq u\leq u_b$,
  	$$\dfrac{{\rm d}v}{{\rm d}u}\leq Av$$ for some positive  constant $A$. It means that $v(u)$ is bounded as  $u\rightarrow u_b^-$, which is a contradiction. Thus it cannot blow up before it extends to the $v-$nullcline $v=g(u)$. Therefore, $S(p)$ exists for all $p\in (1,\frac{1}{a})$ and $S(p)\geq v^ *$.

  	Notice that $S(p)=v^*$ only when  $\Gamma_p^s\rightarrow E^*$ as $t\rightarrow-\infty$, which means $E^*$ must be an unstable node. When $1<p<p_H$, $E^*$ is locally stable, thus $S(p)>v^*$ if $p<p_H$.  When $p_H<p<\frac{1}{a}$ and near $p_H$, $E^*$ is an unstable spiral.
  	Thus, the set $K=\{p\in (1,\frac{1}{a}):S(p)>v^*\}$ is a nonempty open set containing $(1, p_H+\varepsilon)$.

  	Now we show that $S(p) $ is monotonically decreasing on any component of $K$. Let $p_1<p_2$ be two points in some interval in  $K$.  Denote the $v-$nullcline for $p_1$ and $p_2$ as $v=g_1(u)$ and $v=g_2(u)$. From the expression of $v=g(u)$, we know that the curve of $v=g_1(u)$ is on the right of $v=g_2(u)$. The stable manifold of $E_a(a,0)$ for $p_1$ and $p_2$ are graphs of function $v_1(u)$ and $v_2(u)$ defined for $u>a$ respectively. The corresponding eigenvectors at $E_a$ are $x_1=(1,r(1-a)-\frac{m(p_1a-1)}{a})^T$,  $x_2=(1,r(1-a)-\frac{m(p_2a-1)}{a})^T$, respectively. Hence, $v_1(u)>v_2(u)$ for $u$ sufficiently near $a$. Suppose that $v_1(u)=v_2(u)$ for some $u$ with $a<u\leq\frac{1}{p_2}$. Let $\overline{u}$ be the smallest such value. Then we must have  $v'_2(\overline{u})\geq v'_1(\overline{u})\geq 0$. But
  	$\dfrac{mv_1(\overline{u})[1-p_1\overline{u}(1+cv_1(\overline{u}))]}{\overline{u}[(1+cv_1(\overline{u}))v_1(\overline{u})-r(1-\overline{u})(\overline{u}-a)]}\leq \dfrac{mv_2(\overline{u})[1-p_2\overline{u}(1+cv_2(\overline{u}))]}{\overline{u}[(1+cv_2(\overline{u}))v_2(\overline{u})-r(1-\overline{u})(\overline{u}-a)]}$ implies that $p_1\geq p_2$, which is a contradiction.
  	Thus, the intersection of $v_2(u)$ and $g_2(u)$, $S(p_2)$, is below the intersection  of $v_1(u)$ and $g_2(u)$.   For $v_1(u)$ between  $g_2(u)$  and $g_1(u)$,  from the vector field on the left of $g_1(u)$,  the intersection of $v_1(u)$ and $g_1(u)$, $S(p_1)$, is higher than that of  $v_1(u)$ and $g_2(u)$. Thus, $S(p_2)<S(p_1)$.
  	 Notice that the argument above shows that if $p_2\in K$ and $S(p_2)>v^ *$, then any $p\in(1,p_2)$ also belongs to $K$ and $S(p)>S(p_2)$.
  	
  	 We can claim now $K=(1,p_a)$ for some $p_a\in(p_H,\frac{1}{a})$. For  $p\in (1,p_a)$, $S(p)>v^*$, and for $p\in[p_a,\frac{1}{a})$, $S(p)=v^*$.
  	 It remains to prove $S(p)$ is decreasing when  $p\in [p_a,\frac{1}{a})$.  From the vector field in (\ref{ODE system}),  $\Gamma_p^s$ moves towards the upper right, backward, before it meets the $v-$nullcline. Then, for $p=p_{{\rm top}}$, we have $S(p)>v^*$,  and thus $p_{\rm top}<p_a$.
  	   It is obvious that $v^*$ is decreasing with respect to $p$  for $p\in(p_{{\rm top}},\frac{1}{a})$, and thus so does for $ p\in[p_a,\frac{1}{a}) $. Therefore, $S(p)$ is decreasing for both $(1,p_a)$ and $[p_a,\frac{1}{a})$.
  	    $~~~~~~~~~~~~~~~~~~~~~~~~~~~~~~~~~~~~~~~~~~~~~~~~~~~~~\Box$
  	
  	 From the monotonicity of $U(p)$ and $S(p)$, we have the following result.
\begin{proposition}\label{hetero}
If $c<\frac{1}{r(1-a)}$, then there exists a unique $p^{\#}\in(1,\frac{1}{a})$, such that $\Gamma_{p^{\#}}^s=\Gamma_{p^{\#}}^u$, forming a heteroclinic orbit from  $E_1$ to $E_a$.
	\end{proposition}
  \noindent Proof.  Notice that $$\lim_{p\rightarrow \frac{1}{a}^-}(S(p)-U(p))<0,   {\rm and} \lim_{p\rightarrow 1^+}(S(p)-U(p))>0. $$	 	
From the  monotonicity of $S$ and $U$, there exists a unique $p^{\#}$ such that $S(p^{\#})=U(p^{\#})$.  $\Box$

 In fact, there is another  heteroclinic orbit from $E_a$ to $E_1$, which is formed by the unstable manifold of $E_a$ and the stable manifold of $E_1$  on the $u-$axis. Thus, there is a loop of heteroclinic orbits from $E_1$ to $E_a$, and then back to $E_1$.

\begin{figure}[htb]\centering
            a)\includegraphics[width=0.45\textwidth]{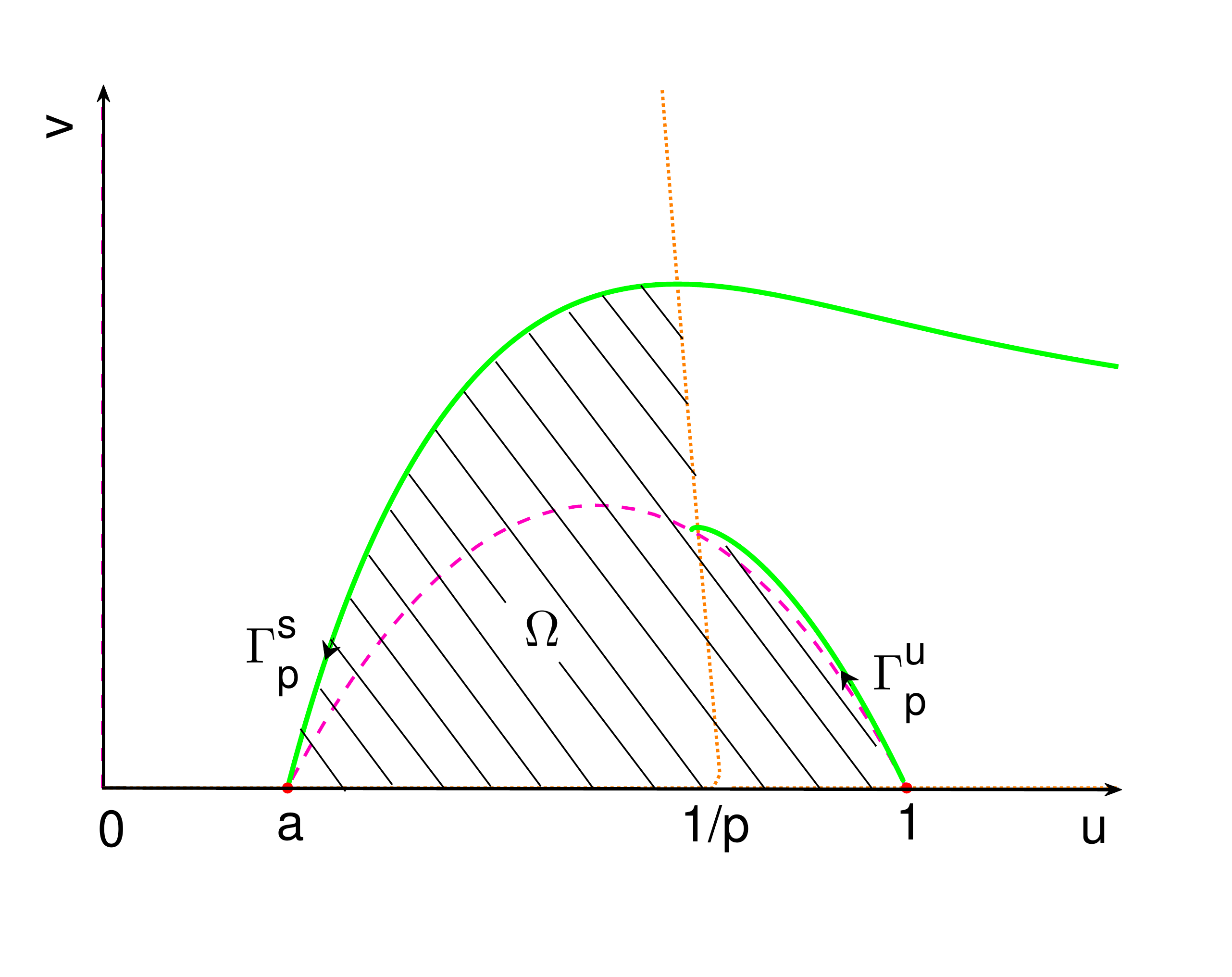}
           b) \includegraphics[width=0.45\textwidth]{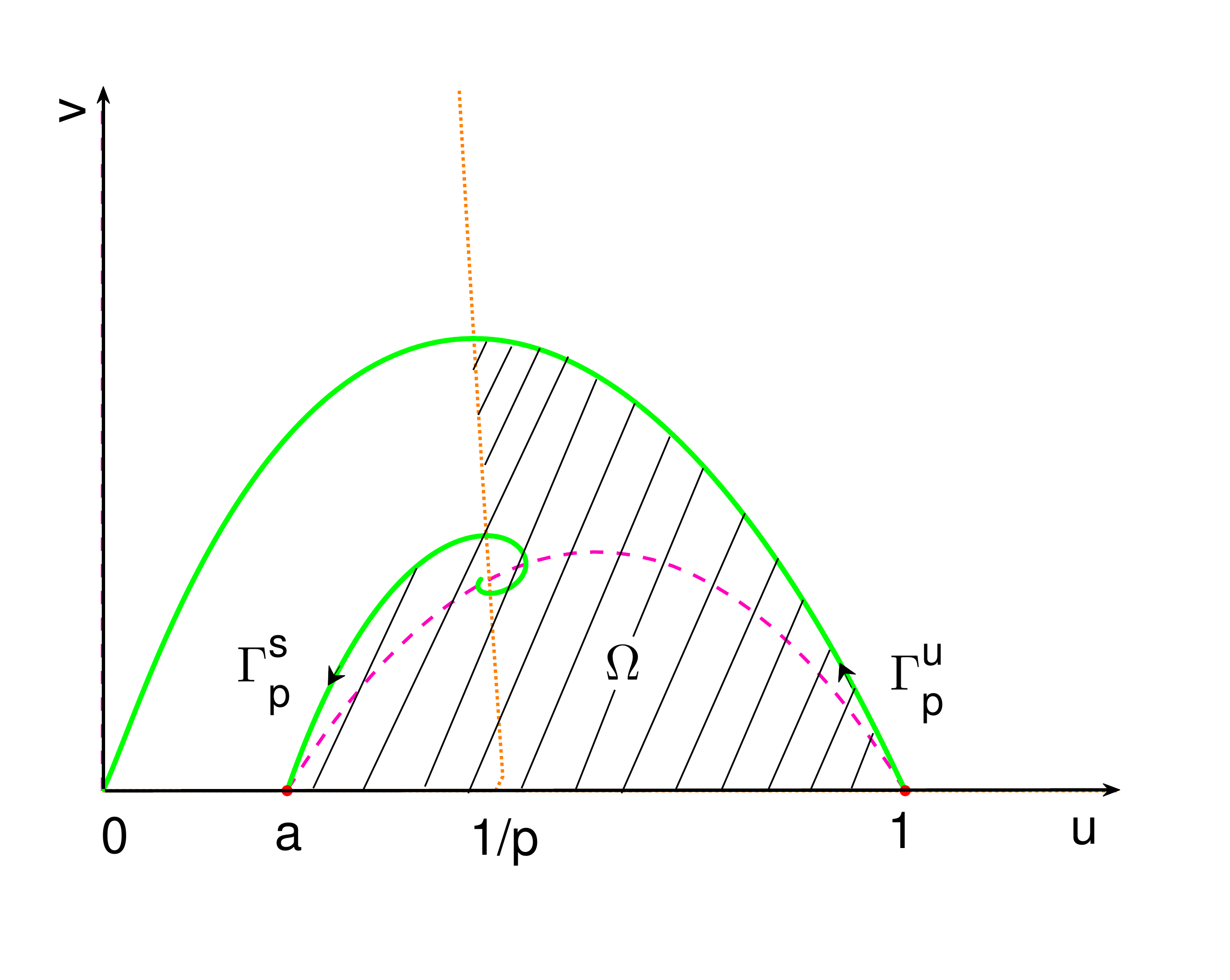}
         \caption{The dashed curve is the $u-$nullcline, and the dotted curve is the $v-$nullcline.    	  Different positions of the stable manifold $\Gamma_p^s$  of $E_a$ and the unstable manifold $\Gamma_p^u$ of $E_1$ for a)$1<p<p^{\#}$  and    b) $p^{\#}<p<\frac{1}{a}$.  	
              	 }	 	\label{us}	
               \end{figure}

  Obviously, $p^{\#}$ is a threshold value of the property of $\Gamma_p^s$ and $\Gamma_p^u$. When $p<p^{\#}$,  $\Gamma_p^s$ is above $\Gamma_p^u$, and when  $p>p^{\#}$, $\Gamma_p^s$ is below $\Gamma_p^u$. Let $\Omega_1$ denotes the bounded open subset of the positive quadrant, bounded by $\Gamma_p^s$, $v-$nullcline between $\Gamma_p^s$ and $\Gamma_p^u$, $\Gamma_p^u$, and the segment from $E_1$ to $E_a$ on the $u-$axis (see Fig.\ref{us}).

  \begin{proposition}\label{pxing}Suppose that $c<\frac{1}{r(1-a)}$.\\
 $~~~~~~$ (i) If $1<p<p^{\#}$, then $S(p)>U(p)$. Moreover, all orbits in the positive quadrant above $\Gamma_p^s$ converge to  $E_0$, and all orbits below $\Gamma_p^s$ have their $\omega-$limit sets in $\Omega_1$, which is a positive invariant set. \\
 $~~~~~~$ (ii) If $p^{\#}< p <\frac{1}{a}$, then $S(p)<U(p)$, and $\Gamma_p^u$ tends to $E_0$. Moreover, all orbits in the positive quadrant above $\Gamma_p^u$ converge to $E_0$, and all orbits below $\Gamma_p^u$ have their $\alpha-$limit sets in $\Omega_1$, which is a negative invariant set.
  \end{proposition}
  \noindent Proof. We prove the case for  $1<p<p^{\#}$, and the other case can be proved similarly. Propositions \ref{manifold} and \ref{hetero} directly lead  to $S(p)>U(p)$.  From Proposition \ref{manifold}, $\Gamma_p^s$ enters $E_a$ from the region above the  $u-$nullcline $v=f(u)$, and meets the  $v-$nullcline at $(u_p^s,S(p))$. On the right of  $(u_p^s,S(p))$, $\Gamma_p^s$ is still above $v=f(u)$, since $\Gamma_p^u$ is above the $u-$nullcline and  $S(p)>U(p)$. Moreover,  $\Gamma_p^s$ can not tend to $E_1$. In fact, from the direction of the vector field in system (\ref{ODE system}), $\Gamma_p^s$ goes to  the lower right   as $t\rightarrow-\infty$. Thus, $\Gamma_p^s$ divides the first quadrant into two regions. Noticing that the first quadrant is invariant, we can get the results from Poncar\'{e}-Bendixson theorem.
  $\Box$

  Moreover, Poincar\'{e}-Bendixson theorem yields the following conclusion.
     \begin{proposition}\label{he}Suppose that $c<\frac{1}{r(1-a)}$.\\
     $~~~~~~$ (i) If $1<p<p^{\#}$, then either $E^*$ is stable or $\Omega_1$ contains a periodic orbits which is stable from the outside (both may be true).\\
     $~~~~~~$ (ii) If   $p^{\#}< p <\frac{1}{a}$, then either $E^*$ is unstable or $\Omega_1$ contains a  periodic orbits which is unstable from the outside (both may be true).
      \end{proposition}

 From proposition \ref{he}, in order to figure out the dynamics in $\Omega_1$ in detail, we have to consider the existence and nonexistence of periodic orbits. In  section \ref{sectionHopf}, we have discussed the existence of periodic orbits bifurcating from Hopf bifurcation. Now, we consider the nonexistence of periodic orbits.
\begin{proposition}\label{nonperi}Suppose that $c<\frac{1}{r(1-a)}$.\\
$~~~~~~$(i) There is an $\varepsilon_1>0$ such that if $\frac{1}{a}-\varepsilon_1<p<\frac{1}{a}$,  there is no periodic orbits, and $\Gamma_p^s$ connects $E^*$ to $E_a$, which is a heteroclinic orbit. \\
$~~~~~~$(ii)  There is an $\varepsilon_2>0$ such that if $1<p<1+\varepsilon_2$,    there is no periodic orbits, and  $\Gamma_p^u$ connects $E_1$ to $E^*$.
	\end{proposition}	
 \noindent Proof.   	 	
We prove the case for $p<\frac{1}{a}$ and near $\frac{1}{a}$, and the latter case can be proved similarly.  If $S(p)=v^*$, $\Gamma_p^s$ connects $E^*$ to $E_a$, and there is no periodic orbits. When $S(p)>v^*$, denote the intersection of $\Gamma_p^s$ and $v-$nullcline as $P (u_p^s,S(p))$. It is easy to confirm that there is a $p_1$ such that  $S(p)= f(u_{p_1})$, where $f(u_{p_1})$ is the component of the intersection of $u-$nullcline and $v-$nullcline corresponding to $p_1$ ( denoted by $P_1 (u_{p_1}, f(u_{p_1})$). Consider the region with vertices $E_a$, $P$, $P_1$, $Q_1(\frac{1}{p_1},0)$, which is a negative invariant region. $P_1$ is well defined when $p$ is close to $\frac{1}{a}$ since $S(p)<f(p_{top})$.  Since $E^*$ is the unique equilibrium in this region, thus if there are periodic orbits, they  must encircle $E^*$ and lie wholly in this region. However,  the divergence of the vector field of (\ref{ODE system}) is positive for $p\rightarrow \frac{1}{a}^-$, since it is
$ra(1-a)+m(pa-1) $ at $E_a$, which is positive. From Bendixson's criterion, there is no periodic orbits in the region. Therefore, according to Poincar\'{e}  Bendixson theorem,  $E^*$ is the $\alpha-$limit set of $\Gamma_p^s$. $~~~~~~~~~~~~~~~~~~~~~~~~~~~~~~~~~~~~~~~~~~~~~~~~~~~~~~~~~~~~~~~~~~~~~~~~~~~~~~~\Box$

  \begin{figure}[htb]\centering
         a) \includegraphics[width=0.45\textwidth]{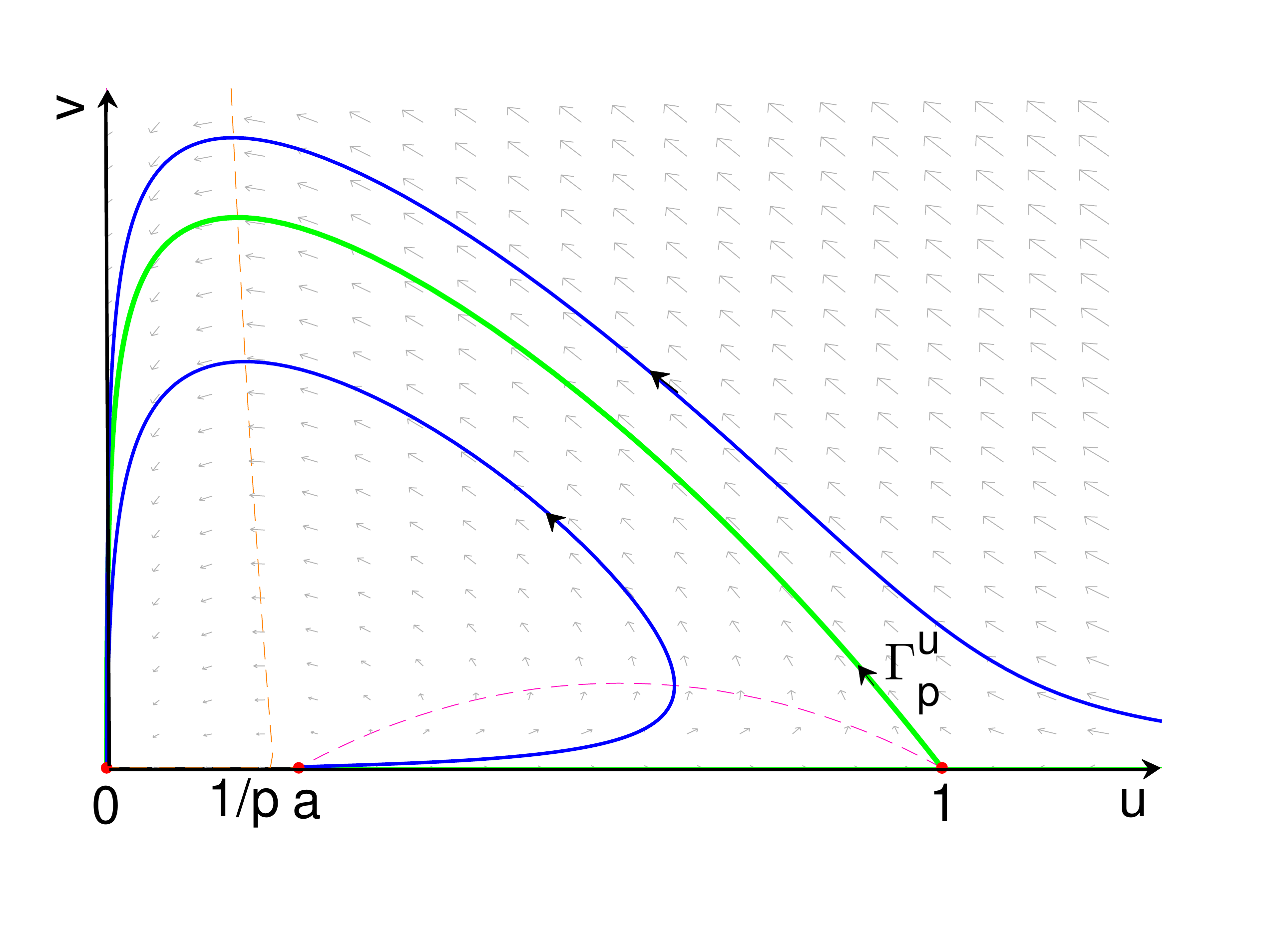}
        b)\includegraphics[width=0.45\textwidth]{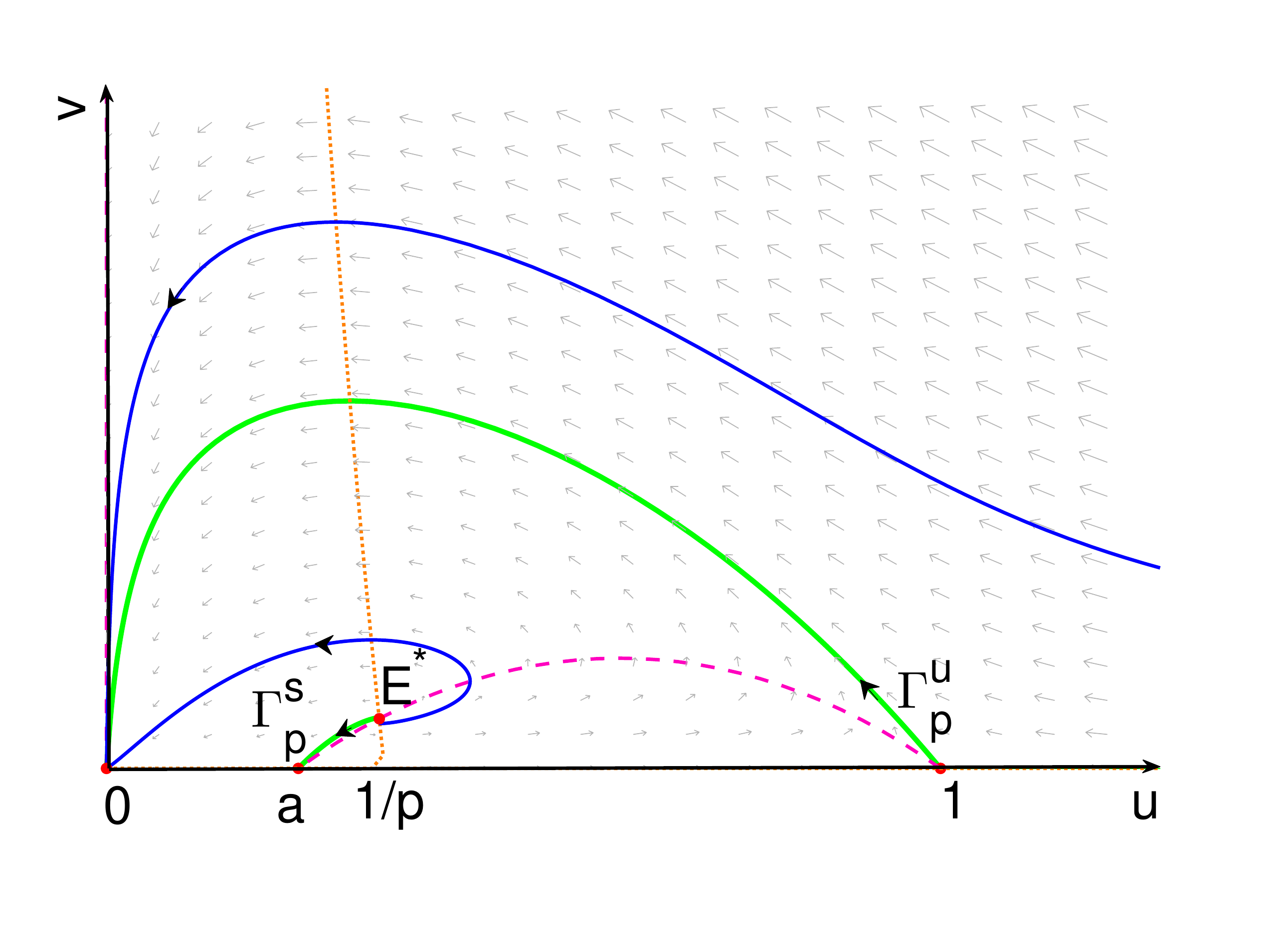}
                c)\includegraphics[width=0.45\textwidth]{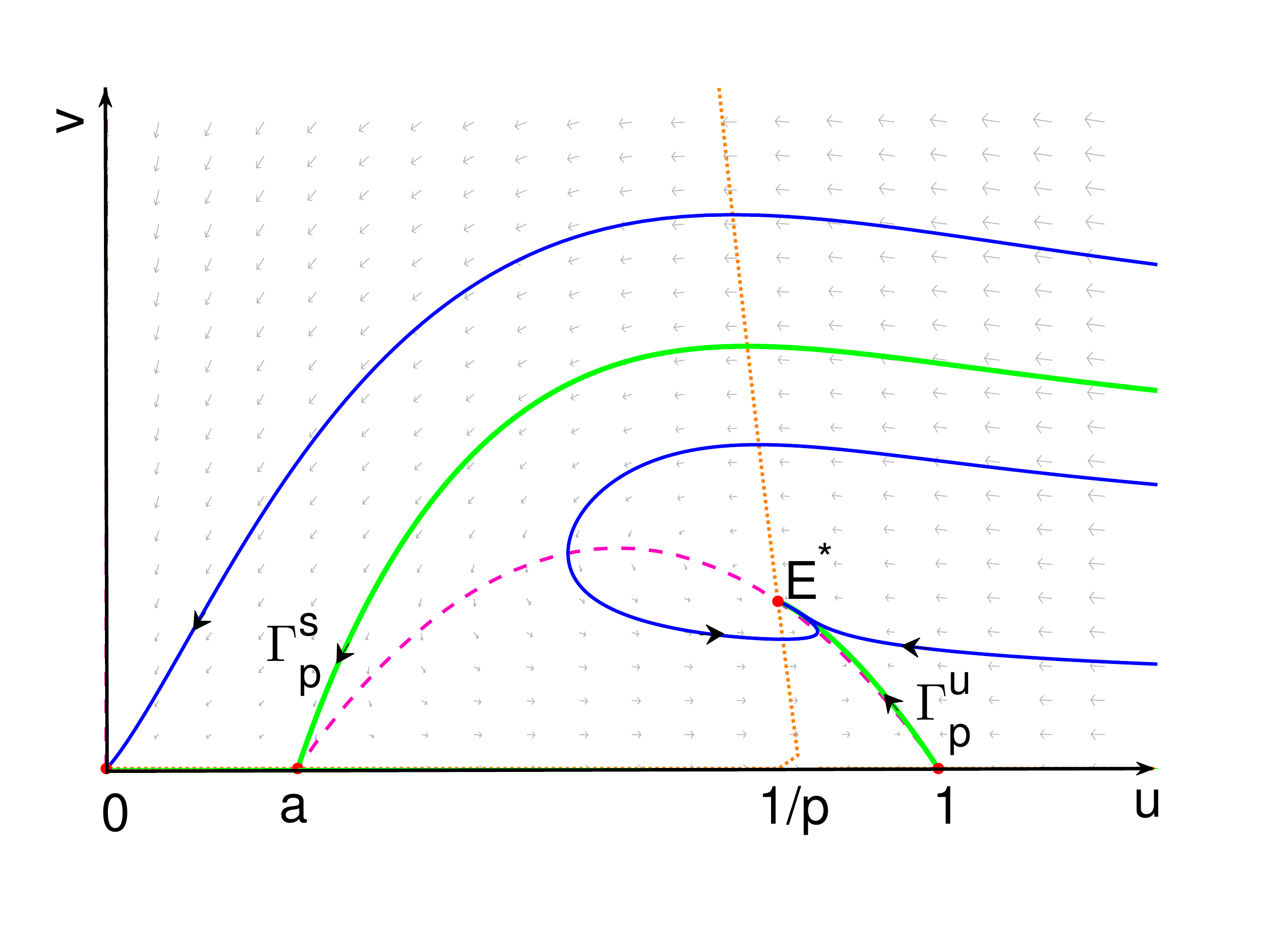}
            d) \includegraphics[width=0.45\textwidth]{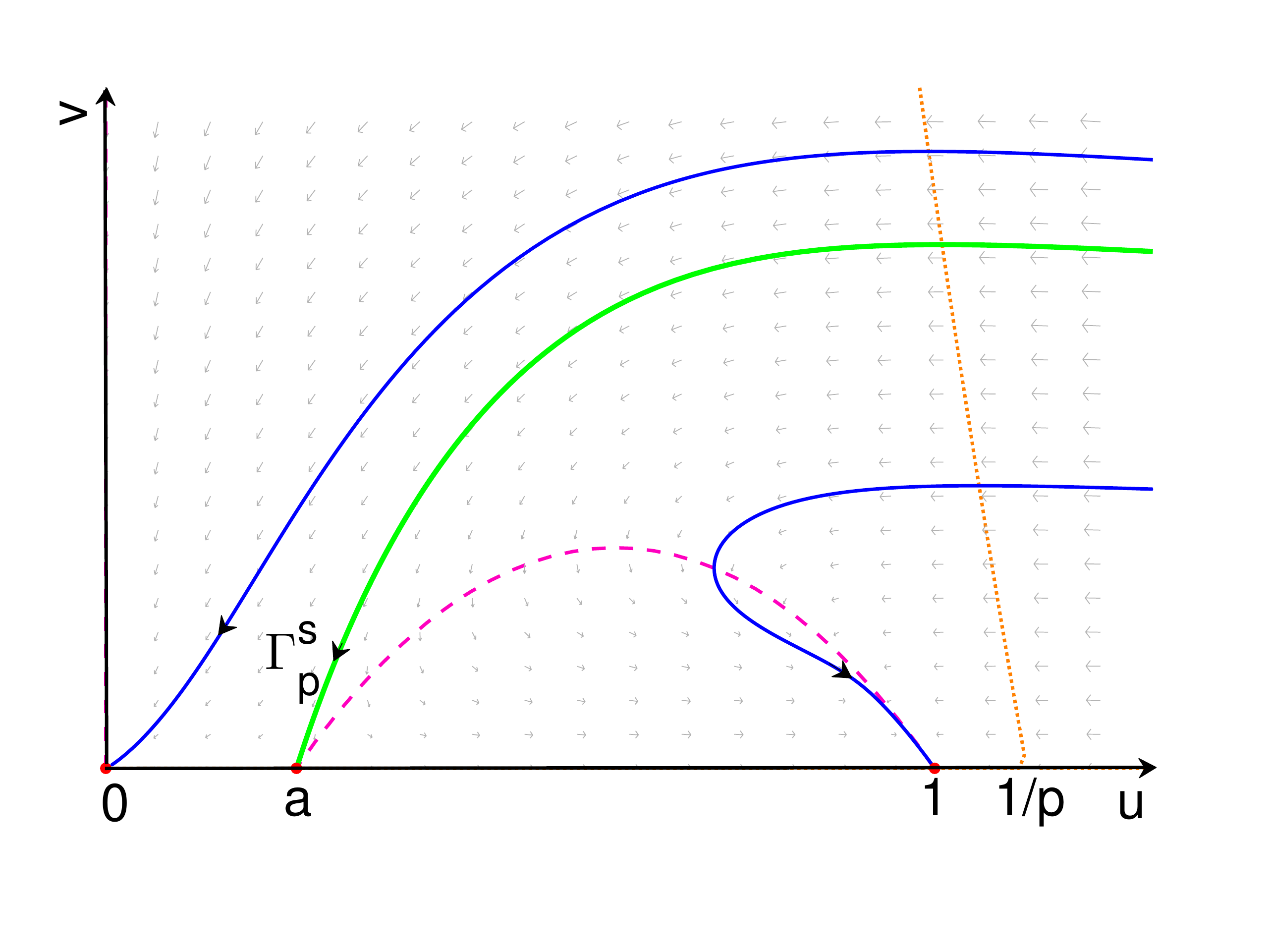}
            	\caption{ 	Phase portrait of system (\ref{ODE system}) for $c<\dfrac{1}{r(1-a)}$ when a) $p>\dfrac{1}{a}$;    b) $p<\dfrac{1}{a}$ and  close to $\dfrac{1}{a}$; c) $p>1$ and close to $1$; d) and $p<1$.     	
            	 }	 	\label{fig:ode01equipxiaoyu1}	
             \end{figure}
             	
\begin{theorem} \label{cxiaoyunormal}Suppose that $c<\frac{1}{r(1-a)}$.\\
$~~~~~~$ (i) If $p\geq\frac{1}{a}$, $E_0(0,0)$ is globally asymptotically stable (see Fig. \ref{fig:ode01equipxiaoyu1}	 a) ).\\
$~~~~$	(ii) There is an $\varepsilon_1>0$ such that   if $\frac{1}{a}-\varepsilon_1<p<\frac{1}{a}$,   $\Gamma_p^s$ connects $E^*$ to $E_a$, and the extinction equilibrium $E_0(0,0)$ is globally asymptotically stable (see Fig. \ref{fig:ode01equipxiaoyu1}	 b) ).\\
 $~~~~$	(iii) There is an $\varepsilon_2>0$ such that if $1<p<1+\varepsilon_2$, then $\Gamma_p^u$ connects $E_1$ to $E^*$. The orbits through any point above $\Gamma_p^s$ converge to $E_0$, and the orbits   through any point below $\Gamma_p^s$ converge to $E^*$ (see Fig. \ref{fig:ode01equipxiaoyu1}	 c) ).	 	 \\
 $~~~~$ (iv) If $p<1$,
 the orbits through any point above $\Gamma_p^s$ converge to $E_0$,
 and   the orbits   through any point below $\Gamma_p^s$ converge to  $E_1$  (see Fig. \ref{fig:ode01equipxiaoyu1}	 d) ).\\
\end{theorem}
\noindent Proof. (i) If $p\geq\frac{1}{a}$, there is no interior equilibrium, then there is no periodic orbit in the first quadrant. Thus, every orbit converges to a boundary equilibrium. We have known that $E_0$ is a stable node,   $E_1$ is a saddle, and  $E_a$ is a unstable node if $p>\frac{1}{a}$  and  a nonhyperbolic repellor if $p=\frac{1}{a}$. Noting that the first quadrant is  positive invariant, then $E_0$ is globally asymptotically stable.

(ii) and (iii) follows from proposition \ref{pxing} and \ref{nonperi} and the fact that $E^*$ is a sink when $p>1$ and close to $1$ and a source when $p<\frac{1}{a}$ and close to $\frac{1}{a}$.

(iv)  If $p$ decreases to $p=1$, $E^*$ collides with $E_1$.  A transcritical bifurcation involving $E^*$ and $E_1$ occurs. $E_1$ changes its stability from a saddle point to a stable node.  If $p<1$,  there is no interior equilibrium, and there is no periodic orbit in the first quadrant.  $E_0$ and $E_1$ are both stable node.  $E_a$ is a saddle, and its stable manifold $\Gamma_p^s$ divides the first quadrant into two regions. The region above $\Gamma_p^s$ is the attractive basin of $E_0$, and  the region below $\Gamma_p^s$ is the attractive basin of $E_1$.
 $~~~~~~~~~~~~~~~~~~~~~~~~~~~~~~~~~~~~~~~~~~~~~~~~~~~~~~~~~~~~~~~~~~~~~~~~~~~~~~~~~~~~~~~~~~~~~~~~~~~~~~~~~~~~~~~~~~~~~~~~~~~~~~~~~~~~~~~~~~~~~~~~~~~~~~~~~~~~~\Box$

Theorem \ref{cxiaoyunormal} provides a global description of the dynamical behaviour of system  (\ref{ODE system}) for $p\in \mathbb{R}^+ \backslash(1+\varepsilon_2,\frac{1}{a}-\varepsilon_1)\stackrel{\vartriangle}{=}\mathbb{R}^+\backslash I$. When $p\in I$, we can go further with the method in  \cite{Conway}.

  \begin{proposition}\label{phpjing}Suppose that $c<\frac{1}{r(1-a)}$.\\
         (i) If $p_H<p^{\#}$, 
         then \\
         $~~~~~~$($i_1$) for $p_H<p<p^{\#} $, 
           there is at least one periodic solution in  $\Omega_1$    which is stable from the outside and one (perhaps the same)  stable from inside;\\
          $~~~~~~$($i_2$) if $a(p_H)>0$, then there is an $\varepsilon>0$ such that if $p_H-\varepsilon<p<p_H$,
          there are at least two distinct periodic solutions, the inner of which is unstable while the outer is stable from the outside.\\
         (ii) If 
         $p^{\#} < p_H $,
         then\\
         $~~~~~~$($ii_1$)  for $p^{\#}<p<p_H$,
         there is at least one periodic solution in  $\Omega_1$    which is unstable from the outside and one (not necessarily distinct) unstable from inside;\\ $~~~~~~$($ii_2$)   if $a(p_H)<0$, then there is an $\varepsilon>0$ such that if $p_H<p<p_H+\varepsilon$,
         there are at least two distinct periodic solutions, the inner
          is stable while the outer is unstable from the outside.
         \end{proposition}
 \noindent Proof. We prove (i), and (ii) can be proved analogously.    If $p_H<p<p^{\#} $,  $E^*$ is a repellor. Since $\Omega_1$ is  a positive invariant region, ($i_1$) follows from the Poincar\'{e}-Bendixson theorem.
  If $a(p_H)>0$, an unstable periodic solution bifurcates from $E^*$ as $p$ decreases past $p_H$.  
  Again $\Omega_1$ is  a positive invariant region, from  Poincar\'{e}-Bendixon theorem, there is a second periodic orbit exterior to the unstable bifurcating  periodic orbit.
 $~~~~~~~~~~~~~~~~~~~~~~~~~\Box$

\begin{figure}[htb]\centering
             a) \includegraphics[width=0.46\textwidth]{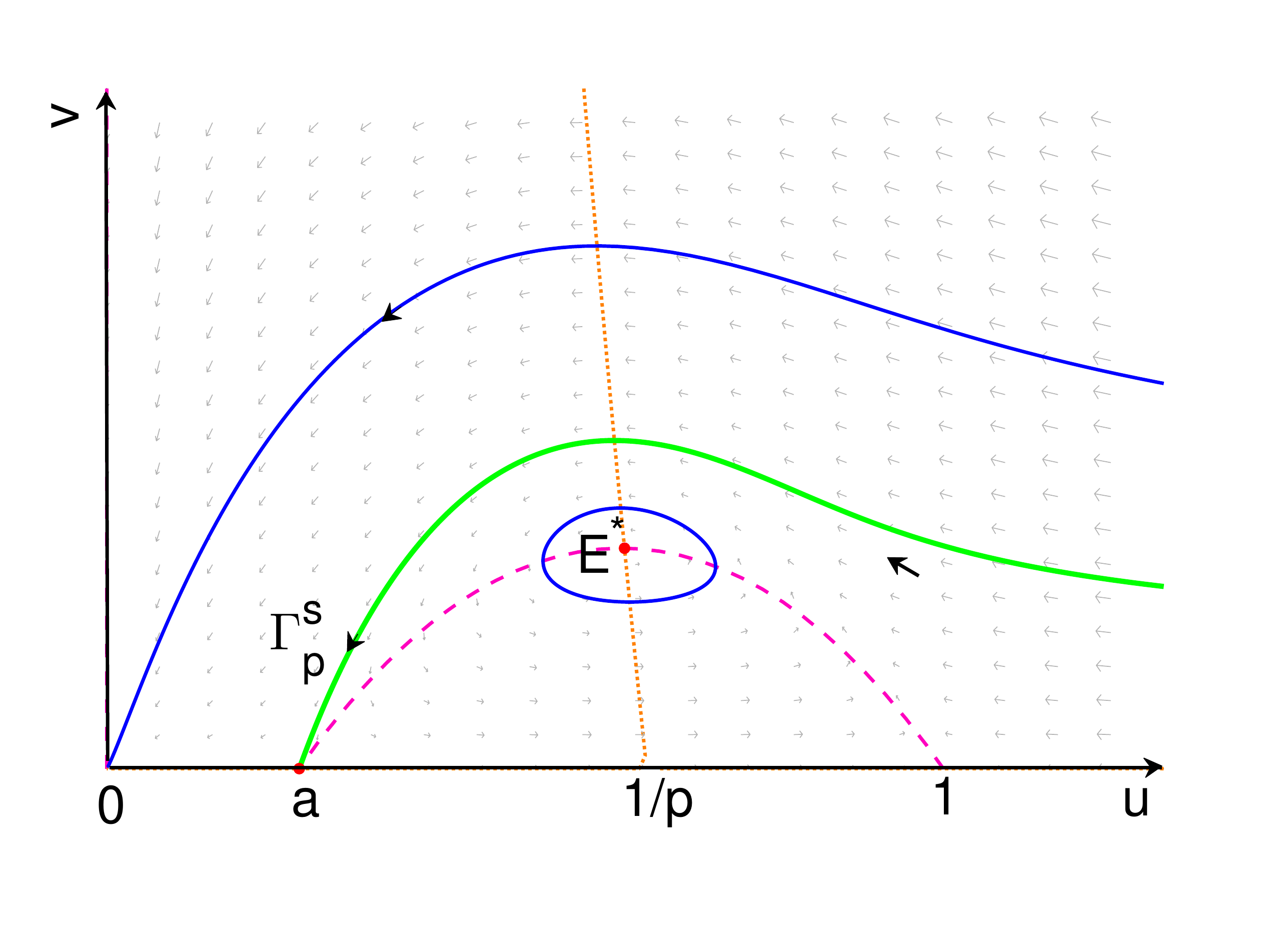}
              b) \includegraphics[width=0.46\textwidth]{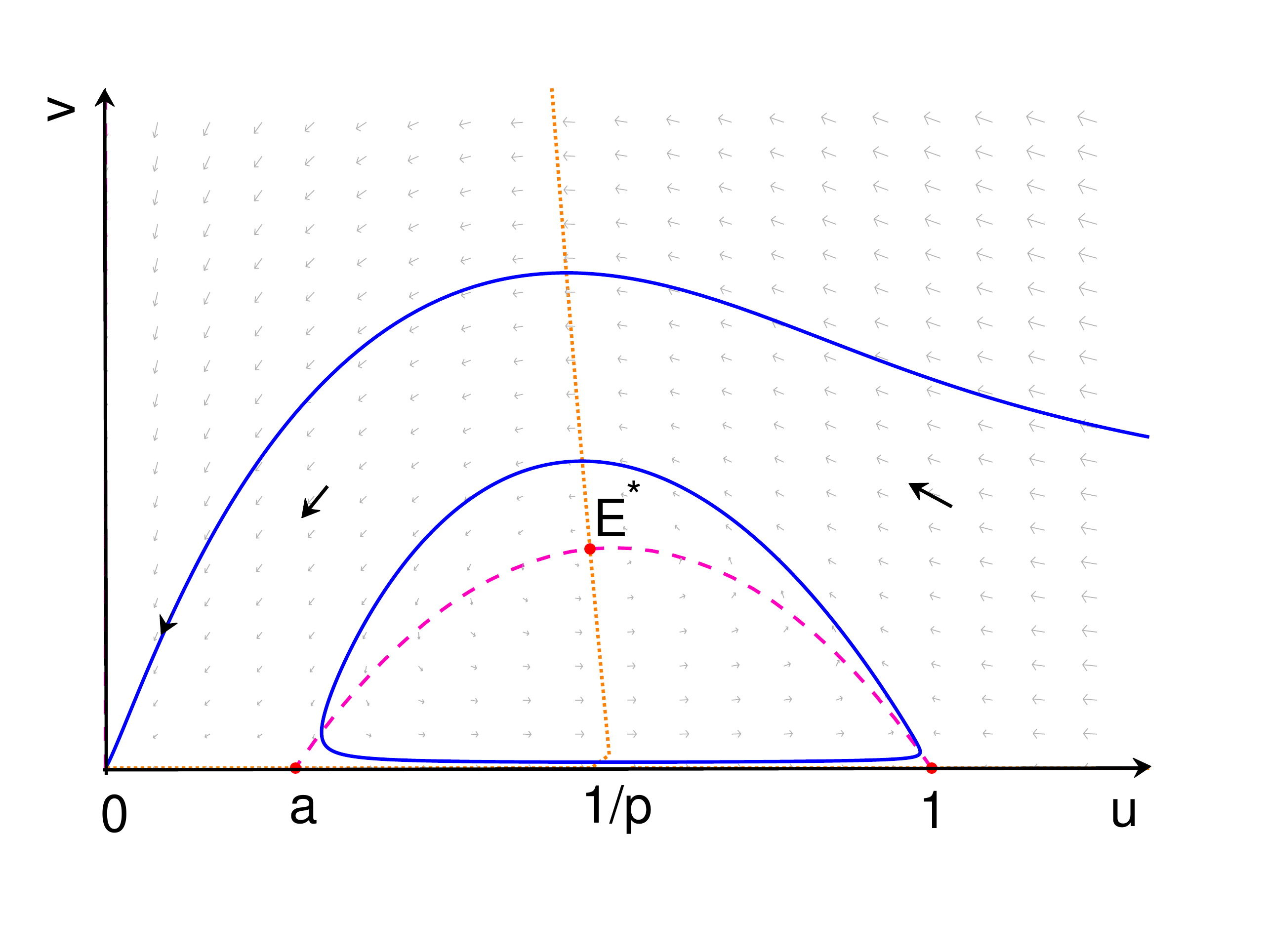}
               c) \includegraphics[width=0.46\textwidth]{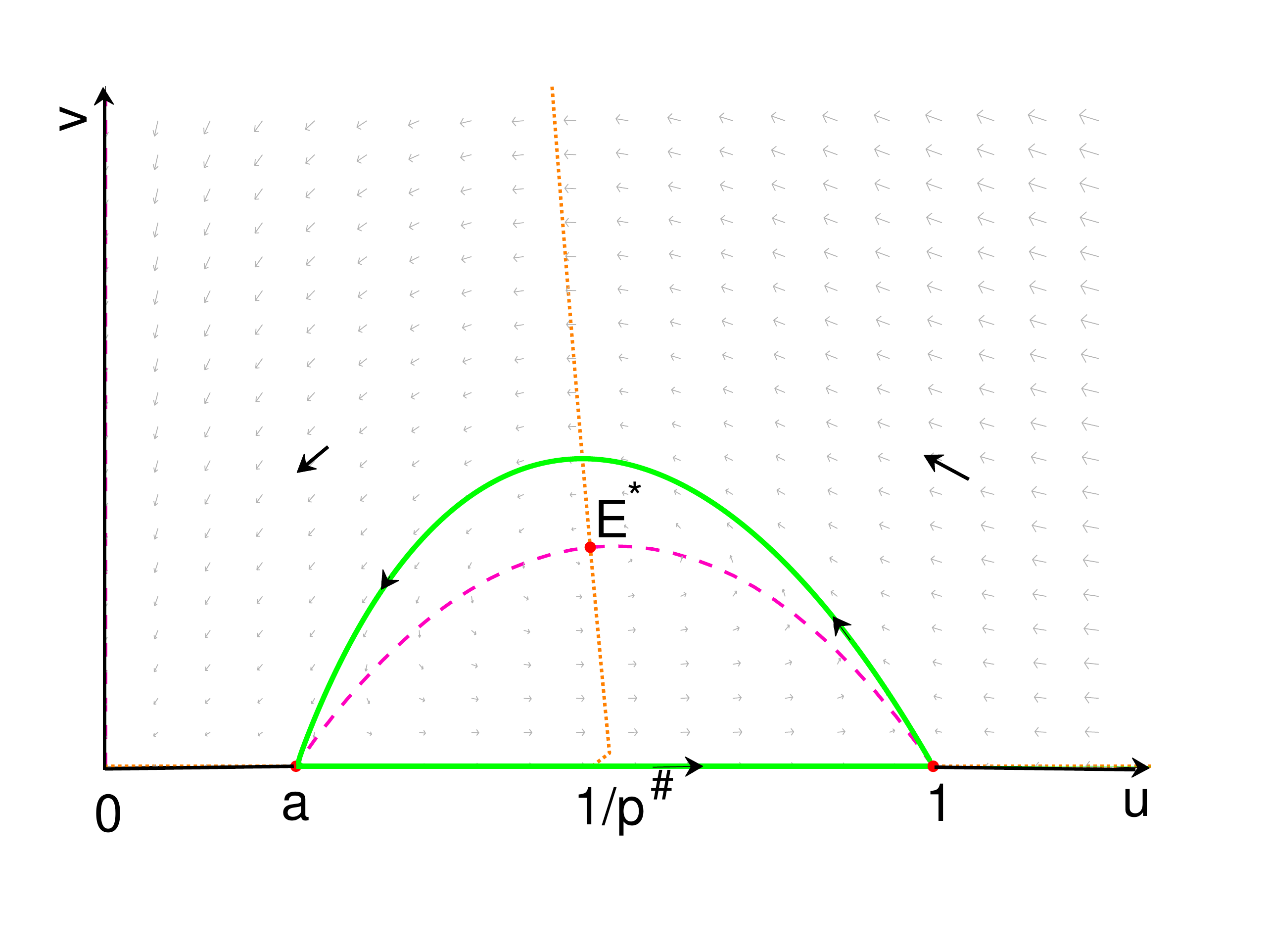}
                d) \includegraphics[width=0.46\textwidth]{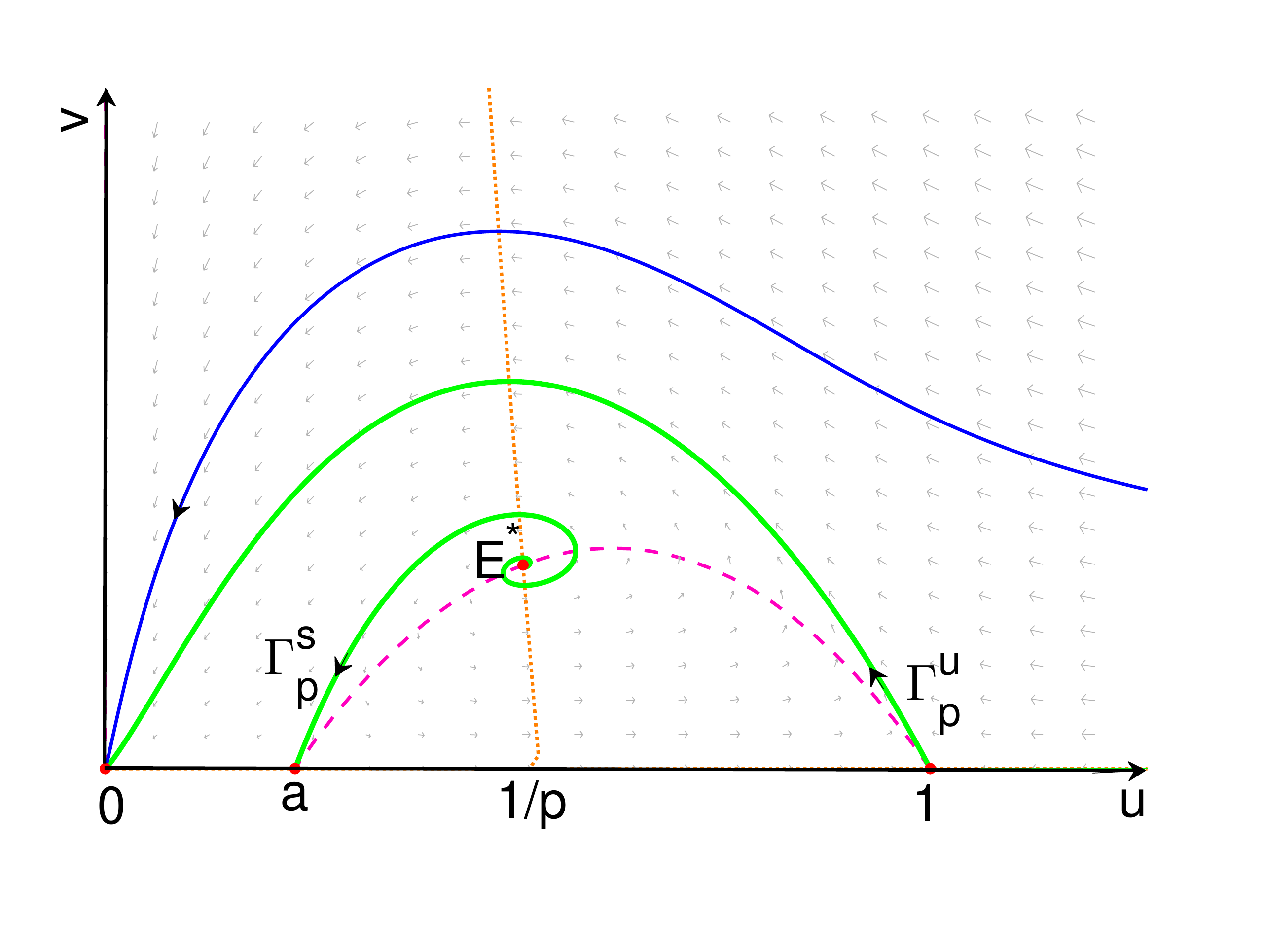}
                 \caption{   	Phase portrait of (\ref{ODE system}) for $c<\dfrac{1}{r(1-a)}$ when a)    $p_H<p<p^{\#}$ and close to $p_H$; b)      $p_H<p<p^{\#}$ and close to $p^{\#}$;  c)   $p=p^{\#}$;   d)    $p^{\#}<p<\frac{1}{a}$. }	 	\label{fig:ode01equi}
          \end{figure}

   If every periodic orbit of system (\ref{ODE system}) is orbitally stable, thus there can be at most one such orbit, and we have the following conclusion.
   \begin{theorem}\label{cxiaoyudynamics} Suppose that $c<\frac{1}{r(1-a)}$, the first Lyapunov  coefficient $a(p_H)<0$, and  every periodic orbit of system (\ref{ODE system}) is  orbitally stable.  Then $p_H< p^{\#}  $.\\
 $~~~~~~$(i) If $1< p <p_H $, $\Gamma_p^u$ connects $E_1$ to $E^*$. The orbits through any point above $\Gamma_p^s$ converge to $E_0$, and the orbits   through any point below $\Gamma_p^s$ converge to $E^*$.	 	  \\
 	 $~~~~~~$(ii) If   $p_H < p <p^{\#} $,  $E^*$ is a repellor, and there is a unique limit cycle under $\Gamma_p^s$. The orbits   through any point below $\Gamma_p^s$ converge to the  limit cycle.  (see Fig. \ref{fig:ode01equi} a) b)).\\
 	 $~~~~~~$(iii) If $p=p^{\#} $, $\Gamma_p^s=\Gamma_p^u$, there are two  heteroclinic orbits forming a loop of  heteroclinic orbits from $E_1$ to $E_a$ and back to $E_1$ (see Fig. \ref{fig:ode01equi} c)). The orbits through any point exterior to the loop converge to $E_0$, and the orbits through any point interior to the cycle  converge to the loop. \\
  $~~~~~~$(iv) If $p^{\#}<p<\frac{1}{a}$, $\Gamma_p^s$ connects $E^*$ to $E_a$, and the extinction equilibrium $E_0(0,0)$ is globally asymptotically stable (see Fig. \ref{fig:ode01equi} d)).
 \end{theorem}	
 \noindent Proof.
  If $a(p_H)<0$, from theorem \ref{Hopfdirection},  stable periodic orbits bifurcating from Hopf
    bifurcation  appear  when $p>p_{H}$. According to proposition \ref{phpjing} (ii), if $p_{\#}<p_H$, there will be at least two distinct periodic orbits,  which   contradicts   our assumption about the uniqueness of periodic orbit. Thus, we have $p_H\leq p_{\#}$. Since $E^*$ is asymptotically stable for any $1<p<p_H$, if there is a periodic orbit, it must be unstable from inside, which contradicts our assumption about the orbital  stability of periodic orbit.
 Therefore, there is no periodic orbit when $1<p<p_H$, and the conclusion (i) follows from proposition \ref{pxing} (i).

 For the case of $p^{\#}<p<\frac{1}{a}$,   any orbit must encircle $E^*$ and lie wholly in $\Omega_1$, which is negatively invariant by proposition \ref{pxing} (ii). Thus, if there is a periodic orbit, it must be unstable from outside,  which contradicts our assumption about the orbital  stability of periodic orbit.    Therefore, the conclusion (iv) follows from proposition \ref{pxing} (ii). Moreover, together with the fact that  Hopf      bifurcating periodic orbit appears when $p>p_{H}$,  the nonexistence of periodic orbits for $p^{\#}<p<\frac{1}{a}$ implies that $p_H$ is strictly less that $p_{\#}$.

  If $p_H < p <p^{\#} $,  $E^*$ is a repellor, thus it follows from proposition \ref{he} (i) that there is a periodic orbit for all $p_H < p <p^{\#} $.  Noticing that we have proved there is no cycle for $p>p^{\#} $, according to the work of \cite{Alexander JC,Chow SN,Wu J} , the period of the unique cycle must tend to infinity as $p\rightarrow p^{\#} $.
Proposition \ref{pxing} (i) completes the proof of (ii).

 When $p=p^{\#}$, then $p>p_H$, and $E^*$ is a repellor, from Poincar\'{e}-Bendixson theorem, the  loop of heteroclicic orbits is the $\omega$-limit set of orbits through any point in its interior. $~~~~~~~~~~~~~~~~~~~~\Box$

\begin{remark}
From the former conclusion on the global dynamics of model (\ref{ODE system}) with weak cooperation, we find that no matter how we choose the value of $p$, both of the two species are extinct if the ratio of  predator to prey is high.

Now we discuss the population behavior when the ratio of  predator to prey is low.   When $p<1$, because of the low conversion rate the predator is  extinct, while the prey reach the capacity of the environment.
When $1<p<p_H$, the predator and prey coexist,   and tends to a stable value. When $p_{H}<p<p^{\#}$,  the predator and prey coexist, but oscillate sustainably.  The amplitude of oscillation increases as $p$ increases, and the minimum of the predator population is close to zero when $p$ is close to $p^{\#}$, which will increase of risk on the extinction of predator.
 When the value of $p$ is too high,  high growth rate of prey leads to over hunting, and finally, both of the two species are extinct.
\end{remark}

 \subsection{The system  with strong cooperative hunting} \label{section3cases}
  By strong cooperation  we mean $c>\frac{1}{r(1-a)}$. In this section, we  first analyze the existence and  stability of interior equilibrium, then we prove the existence of Hopf bifurcation and give the condition of the existence and nonexistence of loop of heteroclinic orbits.  We also exhibit the complex   dynamics of model (\ref{ODE system}), such as limit cycle, loop of heteroclinic orbits, and homoclinic cycle.

\subsubsection{Existence and stability of equilibria}
\begin{figure}
    	\centering
    		a)   \centering
    	      	\includegraphics[width=0.45 \textwidth]{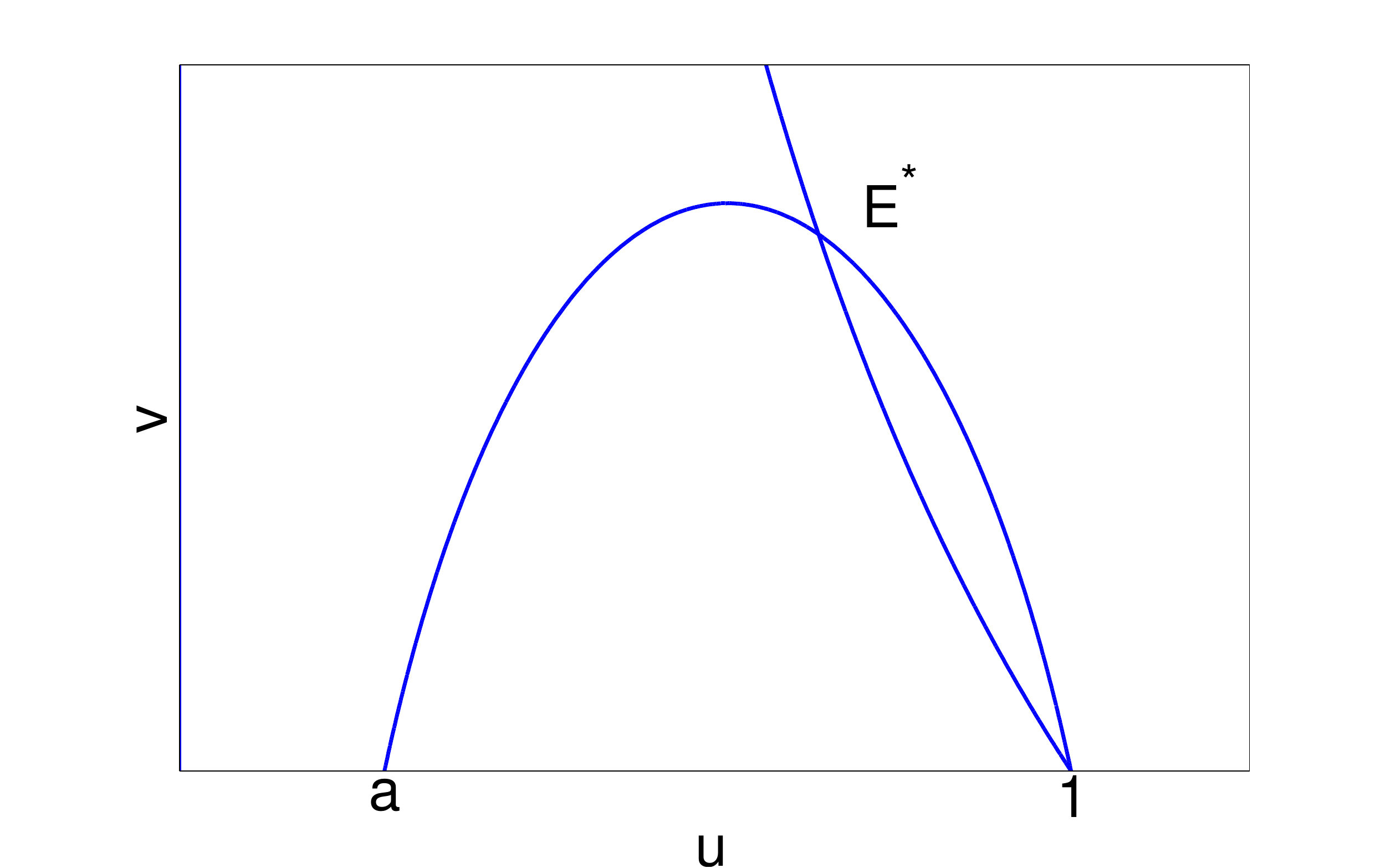}
    b)\centering
    		\includegraphics[width=0.45 \textwidth]{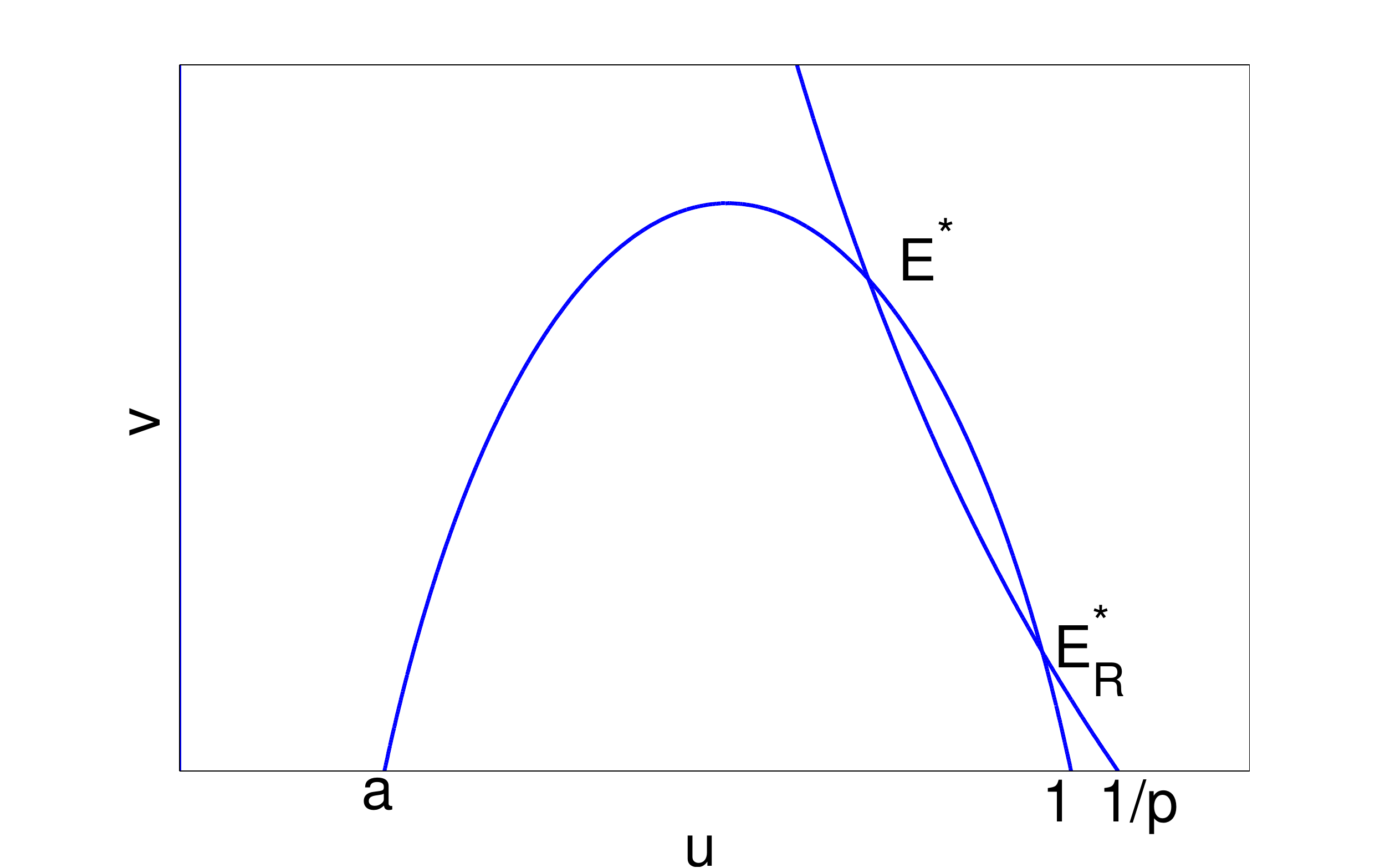} \\  c)   \centering
    	\includegraphics[width=0.45 \textwidth]{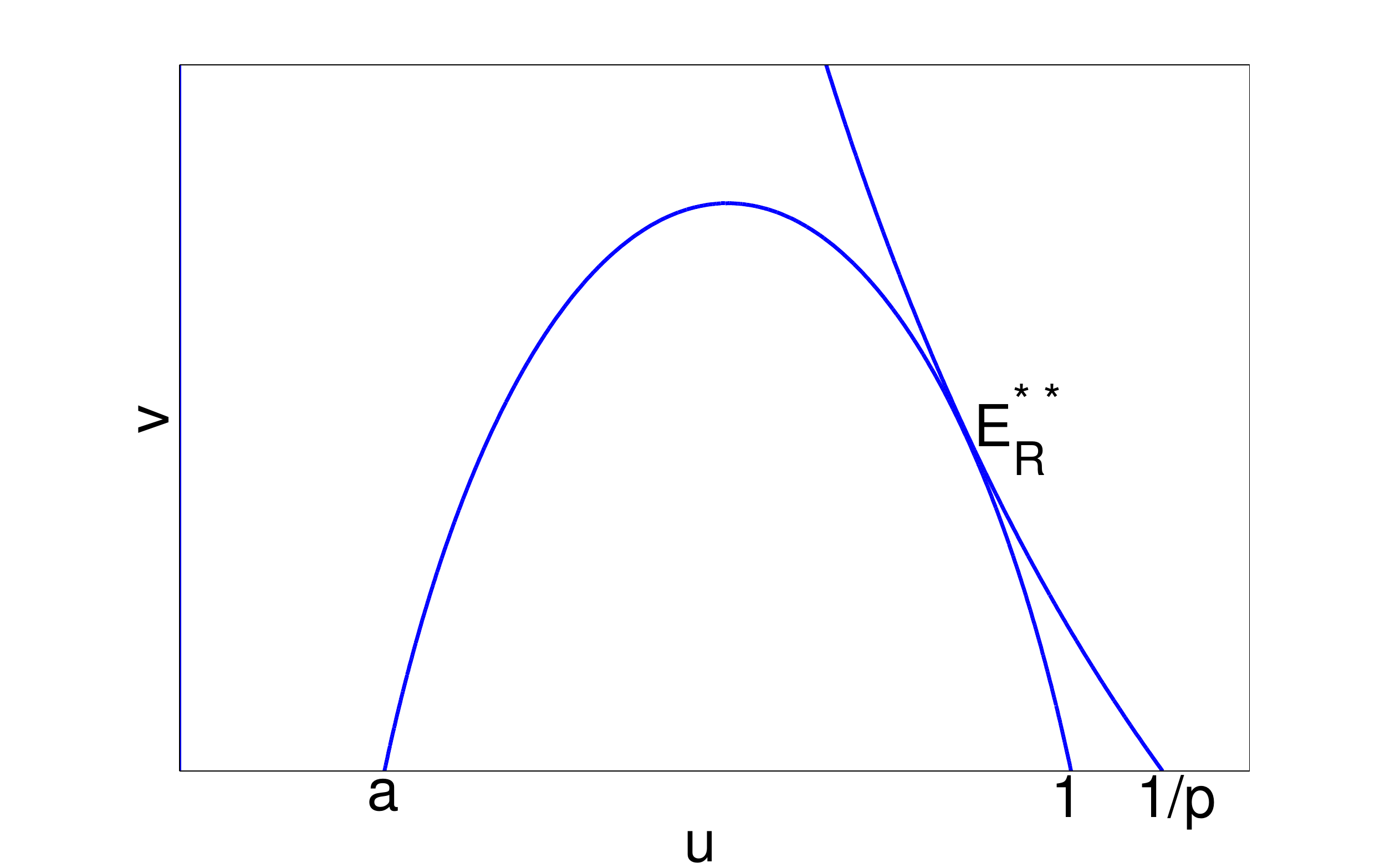}
    	d)    \includegraphics[width=0.45 \textwidth]{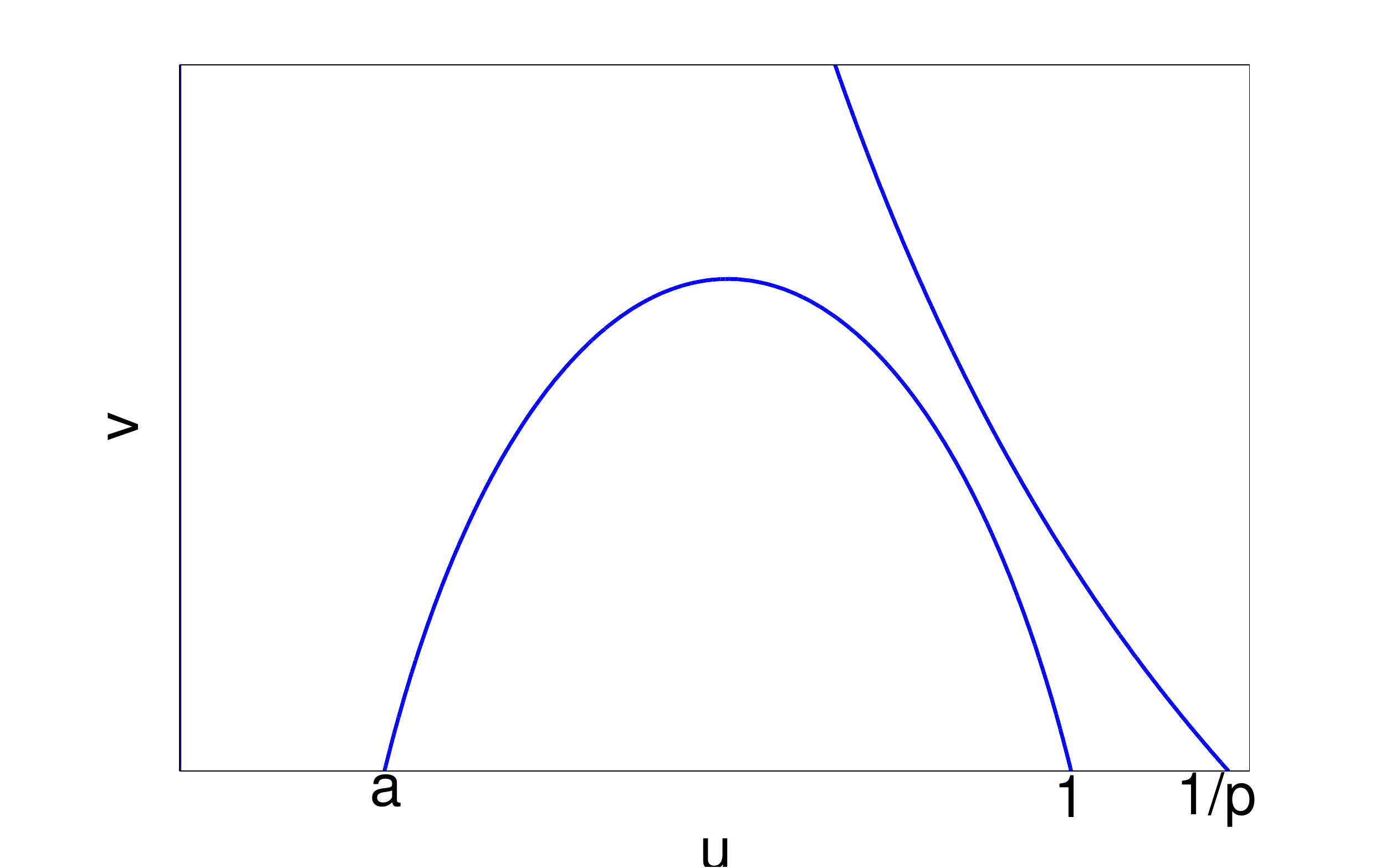}
    	\caption{If
    	 $c>\dfrac{1}{r(1-a)}$,	there is a) one interior equilibrium when $p=1$; b)  two interior equilibria when $p_{SN}<p<1$; c) one equilibrium when $ p=p_{SN}$; and d)  no interior equilibria when $p<p_{SN}$.}
    	\label{fig:kxiaoyu1equilibrium}
    \end{figure}

 	When $c>\frac{1}{r(1-a)}$, the existence and stability of the  boundary equilibria are the same as the case of $c<\frac{1}{r(1-a)}$, which   have been discussed in lemma     \ref{boundary}. For interior equilibria, it is more complicated, and we have the following conclusion.
\begin{proposition}\label{numberofequilibriatrong}
	Suppose that $c>\frac{1}{r(1-a)}$. \\
$~~~~~~$	(i)  When $p\geq\frac{1}{a}$, system (\ref{ODE system})  has  no interior equilibria.\\
$~~~~~~$	(ii) When $1\leq p<\frac{1}{a}$,  system (\ref{ODE system}) has a unique positive constant equilibrium $E^*(u^*,v^*)$.\\
$~~~~~~$ 	(iii) There exists a $p_{SN}<1$, such that system (\ref{ODE system})  has   \\
$~~~~~~$	($a$) two interior equilibria $E^*$ and $E^*_R$ when $p_{SN}<p<1$;\\
$~~~~~~$	(b)  one interior equilibrium  $E_R^{**}$ when $p=p_{SN}$; and \\
 $~~~~~~$ (c) no interior equilibria when $p<p_{SN}$. %
\end{proposition} 	
\noindent Proof. The proofs for $p\geq\frac{1}{a}$ and $1< p<\frac{1}{a}$ are the same as that of proposition \ref{numberofequilibrium}.
  When $\frac{1}{p}$ increases to $1$, the $v-$nullcline $v=g(u)$ and $u-$nullcline  $v=f(u)$  intersect at $E_1(1,0)$.  The tangents of $v-$nullcline and $u-$nullcline at $E_1(1,0)$ are $-\frac{1}{c}$ and $-r(1-a)$ respectively.  If $-\frac{1}{c}> -r(1-a)$,  $E_R^*$ (sitting in the forth quadrant when $1<p<\frac{1}{a}$) collides with $E_1$ when $p=1$, and there is a unique interior equilibrium  $E^*$ in the first quadrant (see Fig. \ref{fig:kxiaoyu1equilibrium} a)).
  When $p<1$ and close to $1$, $E_R^*$ moves into the first quadrant, and there are two  interior equilibria $E^*$ and $E_R^*$ (see Fig. \ref{fig:kxiaoyu1equilibrium} b)).   If  $p<1$ goes on to decrease, there exists a $p_{SN}$, such that the $u-$nullcline $v=f(u)$ and $v-$nullcline $v=g(u)$ tangent at $E_R^{**}$, where $E^*$ and $E_R^*$ collide, and there is a unique interior equilibrium $E_R^{**}$ (see Fig. \ref{fig:kxiaoyu1equilibrium} c)). When $p<p_{SN}$, there is no interior equilibria (see Fig. \ref{fig:kxiaoyu1equilibrium} d)).
 $~~~~~~~~~~~~~~~~~~~~~~~~~~~~~~~~~~~~~~~~~~~~~~~~~~~~~~~~~~~~~~~~~~~~~~~~~~~~~~~~~~~~~~~~~~~~~~~~~~~~~~~~~~~~~~~~~~~~~~~~~~~~~~~~~~~~~~~~~~~~~~~~~~~~~~~~~~~~~~~~\Box$

  \begin{remark}\label{conequili}
  	In the absence of cooperative hunting within the predator, i.e., $c=0$,   system (\ref{ODE system}) always has a unique interior equilibrium $E_0^*=(\frac{1}{p}, r(1-\frac{1}{p})(\frac{1}{p}-a))$.  The predator and prey  coexists if and only if  $1<p<\frac{1}{a}$, and the predator is distinct if $p<1$. However, from proposition \ref{numberofequilibriatrong}, the predator and prey  still coexist when $p_{SN}<p<1$ with strong cooperation, which means that strong cooperation is beneficial for the survival of  predator.
  \end{remark}

   Now we consider the stability of interior equilibria. For  $E^*$ or $E_R^*$, the corresponding Jacobian matrix is
   \begin{equation}\label{Jacoode}
     J=\left( \begin{array}{cc}
     ru(1+a-2u)&-2cuv-u\\mpv(1+cv) &mpcuv
     \end{array}\right),
     \end{equation}
     and thus
     \begin{equation*}
     \begin{array}{l}
          {\rm tr}J=-ru(2u-a-1)+mpcuv=-ru(2u-a-1)+m(1-pu),\\
      {\rm det}J=mpuv\left[ rcu(1+a-2u)+(1+cv)(1+2cv)\right] =\frac{mv}{u}\left[rcpu^3(1+a-2u)+\frac{2}{p}-u \right] .
      \end{array}
    \end{equation*}
         \begin{figure}
    	\centering
     \includegraphics[width=0.5\textwidth]{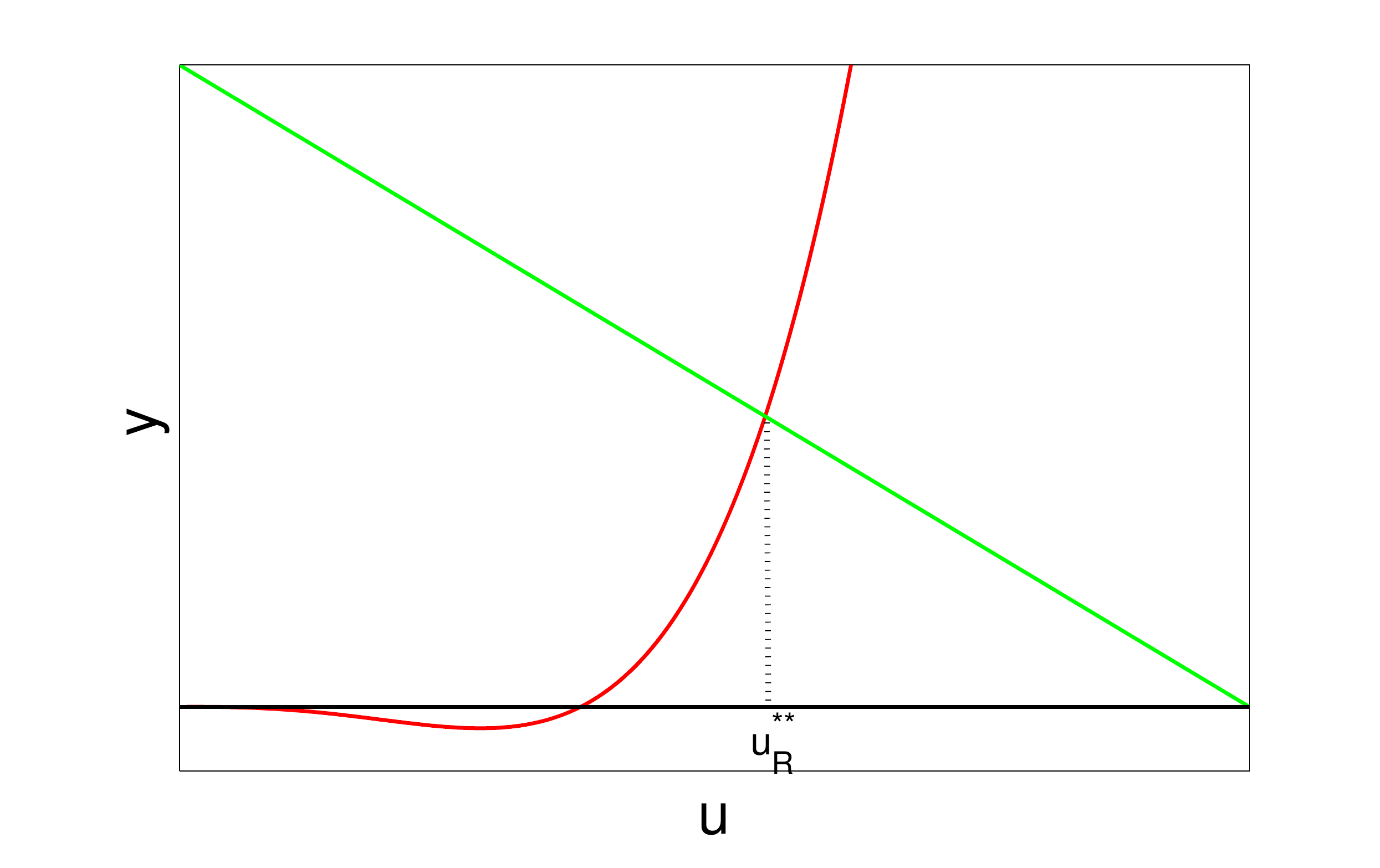}
     	\caption{
    		The red curve  represents  the  figure of $y=rcpu^{3}(2u-a-1)$, and the green  straight line represents for $y=2/p-u$.     	}
    	\label{fig:conExingstability}
    \end{figure}
   If $ {\rm det}J=0$,  then $ \frac{r(1+a-2u)}{1+2cv}=\frac{1+cv}{-cu}$,   which means that $v=f(u)$  is tangent to $v=g(u)$.  In fact, the point of tangency is $E_R^{**}$,    thus ${\rm det}J_{E_R^{**}}=0$.   It is clear that  ${\rm det}J>0$ for $u<u_R^{**}$, and
    ${\rm det}J<0$ for $u>u_R^{**}$ (see Fig. \ref{fig:conExingstability}).
      The component $u^*$ of $E^*$ is always satisfying $u^*<u_R^{**}$, and  the component $u_R^*$ of $E_R^*$ is  satisfying $u_R^*>u_R^{**}$ (see Fig.  \ref{fig:kxiaoyu1equilibrium}).  It follows that for $E^*$, ${\rm det}J_{E^*}>0$, while for $E_R^*$, ${\rm det}J_{E_R^*}<0$.
     Thus,   $E_R^*$ exists when $p_{SN}<p<1$, and it is a saddle.  $E^*(u^*,v^*)$ exists when $p_{SN}<p<\frac{1}{a}$,  and ${\rm det}J_{E^*}>0$.  Thus, $E^*$ may be node or focus,  and its stability depends on  the trace \begin{equation}{\rm tr}J_{E^*}=-ru^*(2u^*-a-1)+m(1-pu^*).\end{equation}

 Similar as the proof in Theorem \ref{HopfExing},  solving for ${\rm tr}J_{E^*}=0$,  when $\frac{a+1}{2}<\frac{1}{p}<\frac{1}{p_{SN}}$, there is a unique $u_{H}^*\in (\frac{a+1}{2},\frac{1}{p_{SN}})$  such that ${\rm tr}J_{E^*}>0$ on
              $(0,u_{H}^*)$, and ${\rm tr}J_{E^*}<0$ on  $(u_{H}^*,\frac{1}{p_{SN}})$ (see Fig. \ref{fig:conExingtr}).
         Recalling that $\frac{\partial u^*}{\partial p}<0$, corresponding to $u_{H}^*$,   we obtain a unique $p_H$ ($p_{SN}<p_H<\frac{2}{a+1}$),    such that ${\rm tr}J_{E^*}>0$ when  $p\in( p_H,\frac{1}{a})$, and ${\rm tr}J_{E^*}<0$ when  $p\in(p_{SN},p_H)$.

          Let $\mu=\alpha(p)\pm i \omega(p)$ be the roots of  $\mu^2-{\rm tr}J_{E^*}\mu+{\rm det}J_{E^*}=0$ when $p$ is near $p_H$. We can prove $\alpha'(p)\mid_{p=p_{H}}>0$, in a similar way as in Theorem \ref{HopfExing}. Therefore, we can get the following  theorem.  	
 	
 	\begin{theorem}\label{cen0deER}
 	When  $c>\frac{1}{r(1-a)}$, there exists a unique $p_H\in (p_{SN},\frac{2}{a+1})$ such that
   $E^*$  is locally asymptotically stable if $p_{SN}<p<p_H$, and unstable if  $p_H<p<\frac{1}{a}$. Moreover,	system (\ref{ODE system}) undergoes a  Hopf
 	 	bifurcation at $E^*$ when  $p=p_H$.
 	  \end{theorem}


 \subsubsection{The global dynamics of system  (\ref{ODE system}) with strong cooperative hunting}
 	When $c>\frac{1}{r(1-a)}$, system (\ref{ODE system}) may have one or two interior equilibria, which brings a few complications to the dynamics of  model  (\ref{ODE system}).  There are  three different cases of global dynamics when $c$ takes different sizes.  In this section, we  exhibit the dynamics  combining theoretical analysis and numerical simulations.
 	
 	Firstly, we consider the properties of the stable manifold of $E_a$ and unstable manifold of $E_1$, denoted by $\Gamma_p^s(E_a)$ and $\Gamma_p^u(E_1)$, respectively.
 	 \begin{proposition}\label{manifoldstrong}
 	 	Suppose that $c>\frac{1}{r(1-a)}$, and $p_{SN}<p<\frac{1}{a}$.   	\\
 	 $~~~~~~$	(i) The orbit  $\Gamma_p^s(E_a)$ meets the $v-$nullcline $v=g(u)$ at a point $(u_p^s,S_a(p))$, where $S_a(p)\geq v^*:=f(u^*)$, and $S_a(p)$ is a monotone decreasing function for $p\in (p_{SN},\frac{1}{a})$.\\
 	 $~~~~~~$ (ii) The orbit  $\Gamma_p^u(E_1)$ meets the $v-$nullcline $v=g(u)$ at a point $(u_p^u,U_1(p))$, where $U_1(p)\geq v^*:=f(u^*)$, and $U_1(p)$ is a monotone increasing function for $p\in (1,\frac{1}{a})$.
 	 		\end{proposition}
 	 	\noindent Proof.	The proof is similar as in proposition \ref{manifold}.
 	 	
 	\begin{theorem}  \label{heterostrong} 	Suppose that $c>\frac{1}{r(1-a)}$.\\
  	$~~~~~~$ (i) If $U_1(1^+)<S_a(1^+)$, there exists a unique $p^{\#}\in(1,\frac{1}{a})$, such that $\Gamma_{p^{\#}}^s(E_a)= \Gamma_{p^{\#}}^u(E_1)$, forming a heteroclinic orbit from $E_1$ to $E_a$; \\
  	$~~~~~~$ (ii) If $U_1(1^+)>S_a(1^+)$,
  	$\Gamma_{p}^s(E_a) \neq \Gamma_{p}^u(E_1)$ for any
  	$p \in(0,+\infty)$.
 	\end{theorem} 	
 	\noindent Proof. (i) If $U_1(1^+)<S_a(1^+)$, then
 	\begin{equation*}
\lim_{p\rightarrow \frac{1}{a}^-}(S_a(p)-U_1(p))<0,   {\rm and} \lim_{p\rightarrow 1^+}(S_a(p)-U_1(p))>0.
 	\end{equation*}
 	 From the  monotonicity of  $S_a(p)$ and $U_1(p)$, there exists a unique $p^{\#}$ such that $S_a(p^{\#})=U_1(p^{\#})$.
 	
 	(ii)  If $U_1(1^+)>S_a(1^+)$, 	from the  monotonicity of $S_a(p)$ and $U_1(p)$,  $U_1(p)>S_a(p)$ for any $ p\in (1,\frac{1}{a})$. Noticing that $E_1$ becomes a stable node when $p<1$ and $E_a$ is an unstable node if $p>\frac{1}{a}$, there is no $p^{\#}$, such that $\Gamma_{p^{\#}}^s(E_a)= \Gamma_{p^{\#}}^u(E_1)$. $~~~~~~~~~~~~~~~~~~~~~~~~~~\Box$	
 	
In fact,  the unstable manifold of $E_a$ on the $u-$axis connects $E_a$ to $E_1$, which is another heteroclinic orbit.  Then, if $U_1(1^+)<S_a(1^+)$,  there is a loop of heteroclinic orbits from $E_1$ to $E_a$, and then back to $E_1$.  If $U_1(1^+)>S_a(1^+)$, such a loop of heteroclinic orbits does not exist.
 	
 	Let $\Omega_{21}$ denotes the bounded open subset of the positive quadrant, boundary with $\Gamma_p^s(E_a)$, $v-$nullcline between $\Gamma_p^s(E_a)$ and $\Gamma_p^u(E_1)$, $\Gamma_p^u(E_1)$, and the segment from $E_1$ to $E_a$ on the $u-$axis.

 	  \begin{proposition}\label{pxingstrong} 	Suppose that $c>\frac{1}{r(1-a)}$, and $1<p<\frac{1}{a}$.\\
 	  $~~~~~~$ (i) Assume $U_1(1^+)<S_a(1^+)$.\\
 	 	 $~~~~~~$ (a) If $1<p<p^{\#}$, $S_a(p)>U_1(p)$. All orbits  in the positive quadrant above $\Gamma_p^s(E_a)$ converge  to  $E_0$. All  orbits below $\Gamma_p^s(E_a)$ have their $\omega-$limit sets in $\Omega_{21}$, which is a positive invariant set. \\
 	 	 $~~~~~~$ (b) If  $p^{\#}< p <\frac{1}{a}$, $S_a(p)<U_1(p)$, and $\Gamma_p^u(E_1)$ enters $E_0$. All orbits  in the positive quadrant above $\Gamma_p^u(E_1)$ converge  to  $E_0$. All orbits below $\Gamma_p^u(E_1)$ have their $\alpha-$limit sets in $\Omega_{21}$, which is a negative invariant set.\\
 	 	  	 $~~~~~~$(ii) Assume $U_1(1^+)>S_a(1^+)$. $S_a(p)<U_1(p)$ for all $p\in(1,\frac{1}{a})$, and $\Gamma_p^u(E_1)$ enters $E_0$. All orbits  in the positive quadrant above $\Gamma_p^u(E_1)$ converge  to  $E_0$. All orbits below $\Gamma_p^u(E_1)$ have their $\alpha-$limit sets in $\Omega_{21}$, which is a negative invariant set.
 	 	  \end{proposition}
 	  \noindent Proof.  The proof is similar as proposition \ref{pxing}.
 	  $~~~~~~~~~~~~~~~~~~~~~~~~~~~~~~~~~~~~~~~~~~~~~~~~~~~~~~~~~~~~~~~~~~~~~~~~~~~~~~~~~~\Box$

  	If $c>\frac{1}{r(1-a)}$, $E_R^*$ exists when $p_{SN}<p<1$, and it is a saddle. Now   we consider the properties of the stable manifold and unstable manifold of $E_R^*$, denoted by $\Gamma_p^s(E_R^*)$ and $\Gamma_p^u(E_R^*)$, respectively.

 	 \begin{proposition}\label{manifoldercing}
 	 	 	Suppose that $c>\frac{1}{r(1-a)}$, and $p_{SN}<p<1$.    The downward unstable manifold $\Gamma_p^u(E_R^*)$ connects $E_1$. The right stable manifold $\Gamma_p^s(E_R^*)$  enters $E_R^*$ from the lower right.	\\
 $~~~~~~$	(i) The orbit of left  $\Gamma_p^s(E_R^*)$ meets the $v-$nullcline $v=g(u)$ at a point $(u_p^s,S_R(p))$, where $S_R(p)\geq v^*:=f(u^*)$.\\
 	 $~~~~~~$ (ii) The orbit of upward  $\Gamma_p^u(E_R^*)$ meets the $v-$nullcline $v=g(u)$ at a point $(u_p^u,U_R(p))$, where $U_R(p)\geq v^*:=f(u^*)$.
 	 	  	  	 	 		\end{proposition}
  	\noindent Proof.
  	  For the negative eigenvalue $\lambda_1 $ of the Jacobian of $E_R^*$, the corresponding eigenvector is $(1,\frac{mpv(1+cv)}{-mpcuv+\lambda_1})$.  	Noticing that the tangent vector  of $v=g(u)$  at $E_R^*$ is $\frac{mpv(1+cv)}{-mpcuv}$,
  	 	 	 the left part of  $\Gamma_p^s(E_R^*)$ near $E_R^*$  is below the $v-$nullcline, and   the right part of  $\Gamma_p^s(E_R^*)$ near $E_R^*$  is above the $v-$nullcline. Obviously,  the right part of  $\Gamma_p^s(E_R^*)$ enters $E_R^*$ from the lower right. From the vector field for (\ref{ODE system}), before  the left part of  $\Gamma_p^s(E_R^*)$ meets the $v-$nullcline, the curve under the $u-$nullcline directs lower right; the curve above the $u-$nullcline directs lower left.  Since it can not cross the stable manifold $\Gamma_p^s(E_a)$, thus, it  is bounded before it meets the  $v-$nullcline. Thus   the left part of  $\Gamma_p^s(E_R^*)$ must meet the $v-$nullcline  at a point, denoted by  $(u_p^s,S_R(p))$.  Obviously, $S_R(p)\geq v^*$.
  	 	 	
  	 	 	 (ii) can be proved similarly. $~~~~~~~~~~~~~~~~~~~~~~~~~~~~~~~~~~~~~~~~~~~~~~~~~~~~~\Box$

 Let $\Omega_{22}$ denotes the bounded open subset of the positive quadrant,  boundary with the left $\Gamma_p^s(E_R^*)$, $v-$nullcline between $\Gamma_p^s(E_R^*)$ and $\Gamma_p^u(E_R^*)$, upward $\Gamma_p^u(E_R^*)$.  Similar as proposition \ref{pxingstrong}, we have the following conclusion. 	\begin{proposition}\label{omegaalphastrong} 	$~~~~~~$ (i) If 	 $S_R(p)>U_R(p)$. All orbits inside the region boundary with $\Gamma_p^s(E_R^*)$ have their $\omega-$limit sets in $\Omega_{22}$, which is a positive invariant set.\\ 	$~~~~~~$ (ii) If  	 $S_R(p)<U_R(p)$. All orbits below the left $\Gamma_p^u(E_R^*)$ have their $\alpha-$limit sets in $\Omega_{22}$, which is a negative invariant set. 		\end{proposition}

\begin{figure}[htb]\centering
 	 	a)\includegraphics[width=0.3 \textwidth]{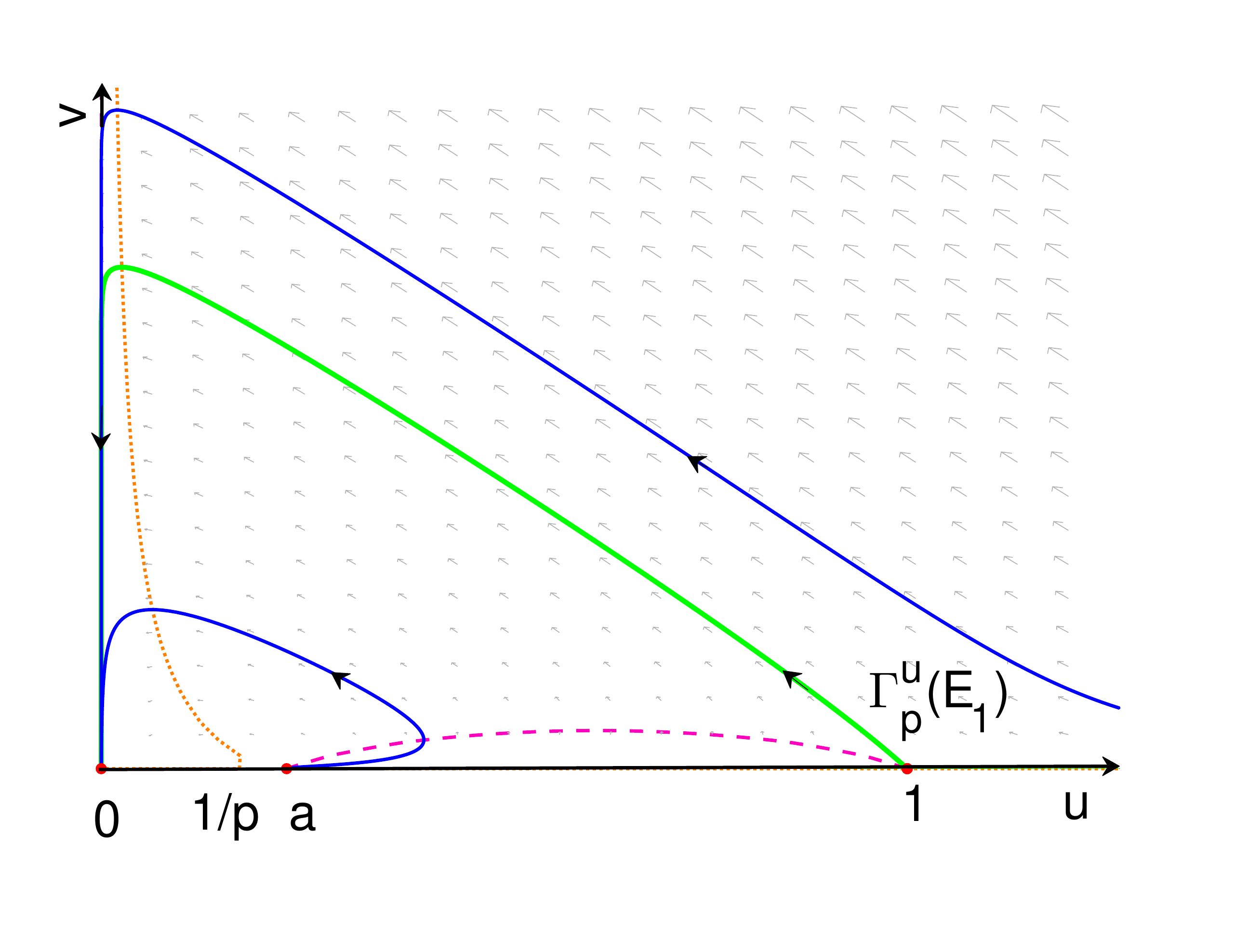} 
 	 	   	b)\includegraphics[width=0.3 \textwidth]{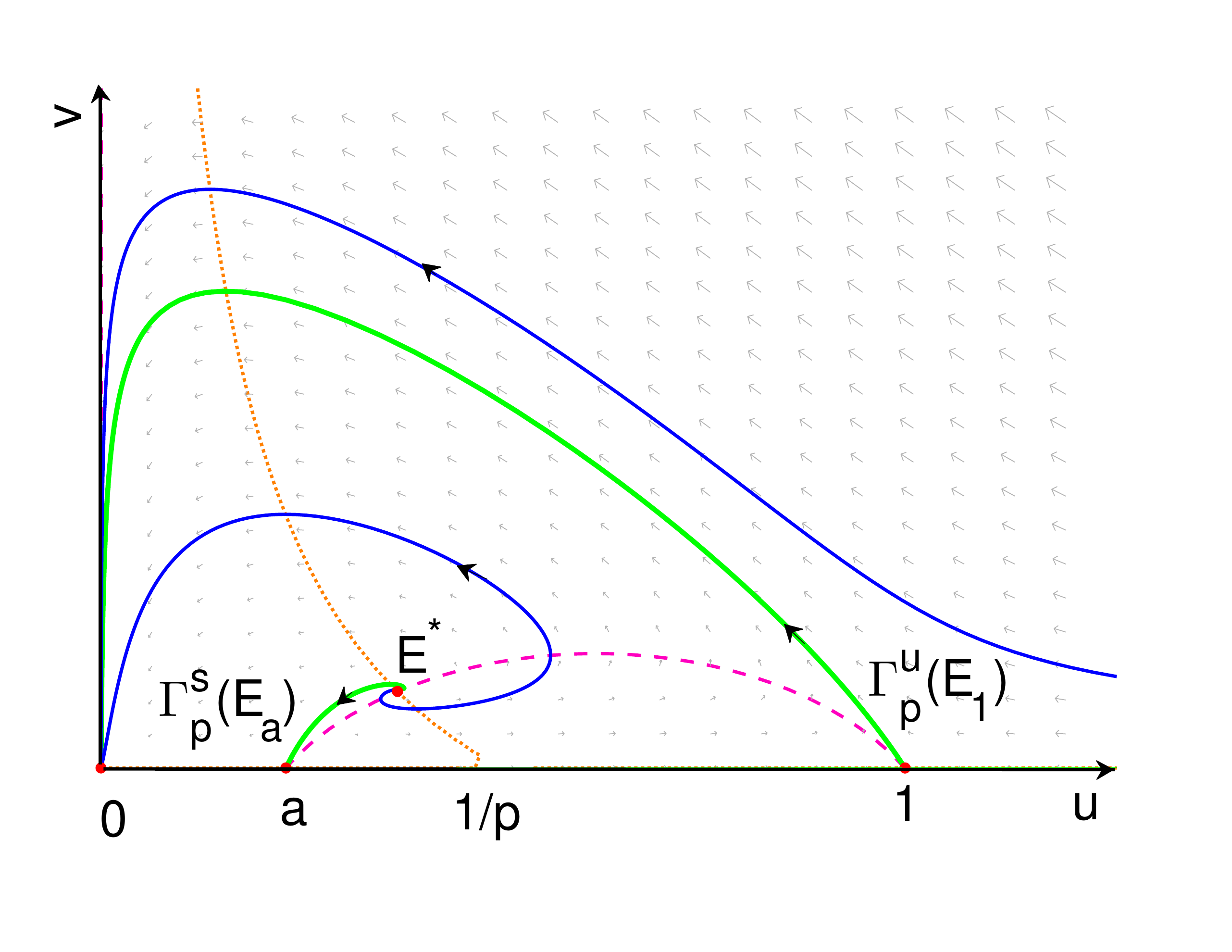} 
 	 	   	 	c)\includegraphics[width=0.3 \textwidth]{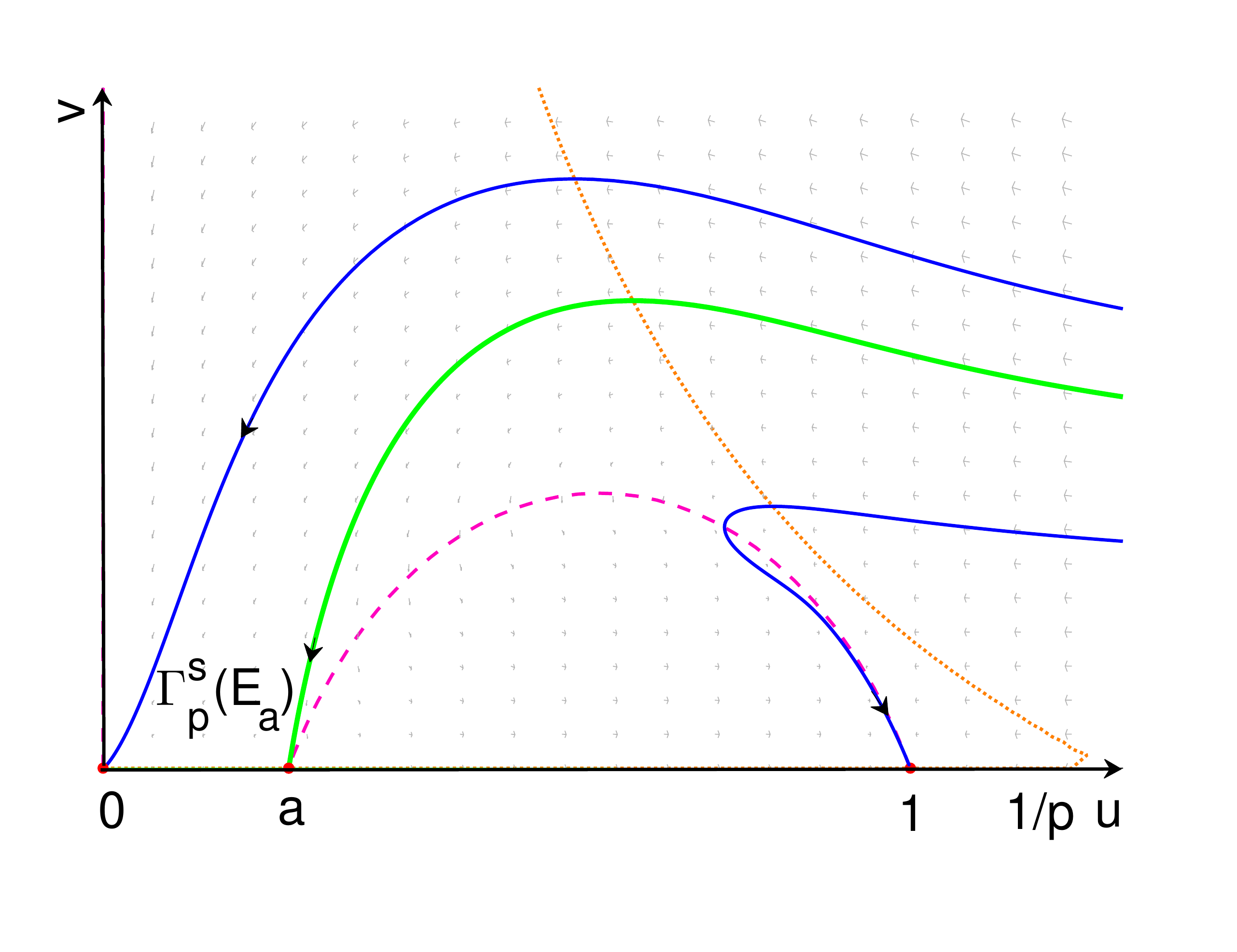}
 	\caption{  If $c>\frac{1}{r(1-a)}$, 	phase portrait of (\ref{ODE system})  when a) $p\geq\dfrac{1}{a}$;    b) $p<\dfrac{1}{a}$ and  close to $\dfrac{1}{a}$;  and c) $p<p_{SN}$.   	 }	 	 \label{fig:czhong1}
 	 	   	 		  \end{figure}

 	
 	From the previous discussion, we can determine the global dynamics of model  (\ref{ODE system}) when $p$ is chosen  in the following ranges. Similar  as in Theorem \ref{cxiaoyunormal}, we can prove the following conclusion.	
 	\begin{theorem} Suppose that $c>\frac{1}{r(1-a)}$.\\
 	$~~~~~~$ (i) If $p\geq\frac{1}{a}$, $E_0(0,0)$ is globally asymptotically stable (see Fig. \ref{fig:czhong1} a)). \\
  	$~~~~~~$	(ii) If $p<\frac{1}{a}$ and near $\frac{1}{a}$, $\Gamma_p^s(E_a)$ connects $E^*$ to $E_a$, and the extinction equilibrium $E_0(0,0)$ is globally asymptotically stable (see Fig. \ref{fig:czhong1}	 b) ).\\
 	 	 $~~~~~~$ (iii) If  $p<p_{SN}$, the orbits through any point   above $\Gamma_p^s(E_a)$  converge to $E_0$, and  the orbits through any point   below $\Gamma_p^s(E_a)$ converge to $E_1$  (see Fig. \ref{fig:czhong1} c)).
 	\end{theorem}

The remaining question is how are the dynamics of system (\ref{ODE system})  for $p$ chosen in the rest ranges, i.e., $p_{SN}<p<\frac{1}{a}$. Since larger $c$ means less steep of the $v-$nullcline, the critical points (for example $p_H$, $p^{\#}$)   may appear in different ranges, and there are the following  three different cases.

\begin{figure}[htb]\centering
 	   	 	a)\includegraphics[width=0.45 \textwidth]{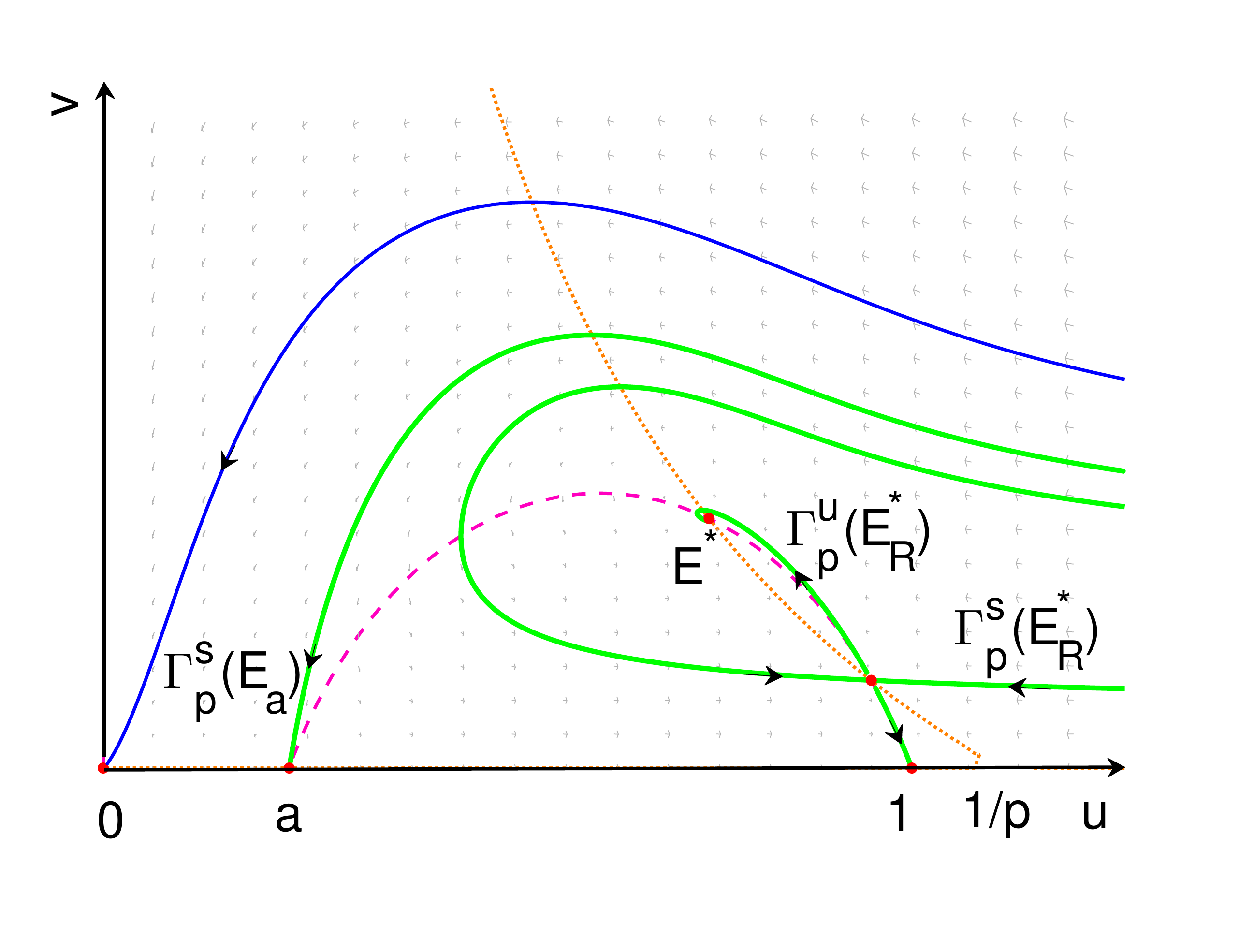}
 	   	 		b)\includegraphics[width=0.45\textwidth]{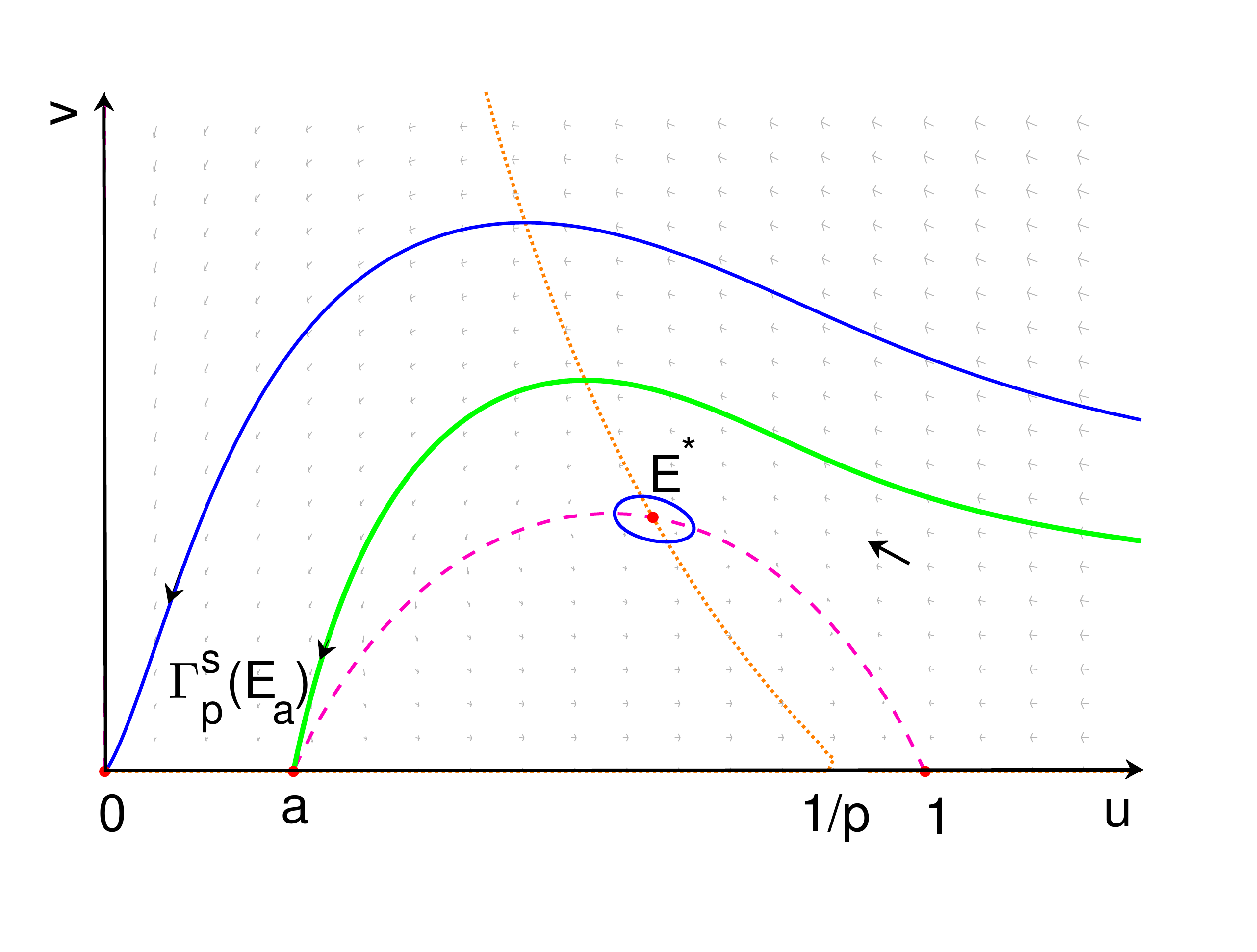}
 	   	 	c)\includegraphics[width=0.45\textwidth]{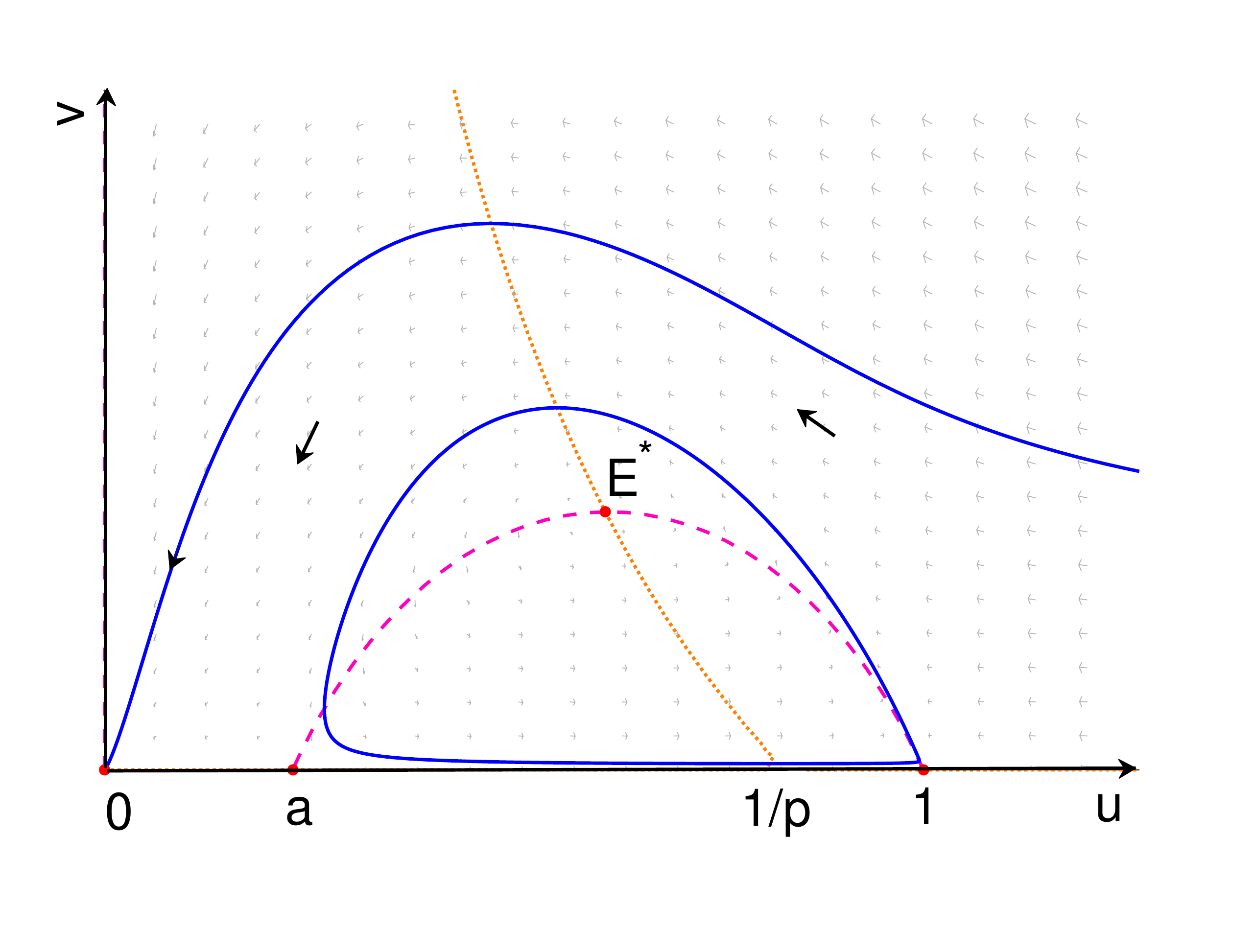}
 	   	 	d)\includegraphics[width=0.45\textwidth]{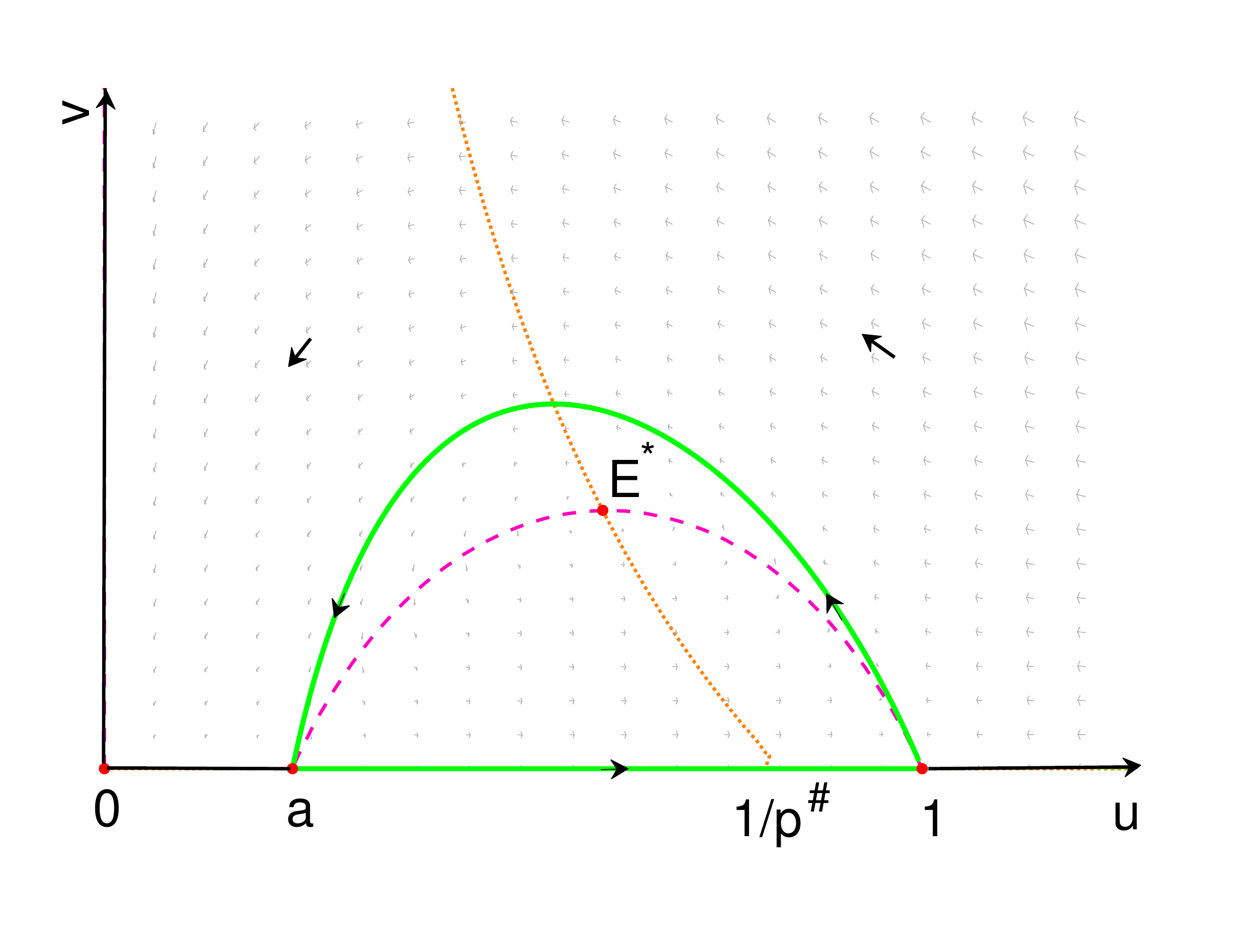}
 	   	 	\caption{ Phase portrait of system (\ref{ODE system})  for case $1$   when a) $p>p_{SN}$ and close to $p_{SN}$; b) $p>p_H$ and close to $p_H$;   c) $p_H<p<p^{\#}$  and  close to $p^{\#}$; and	d) $p=p^{\#}$.  }\label{fig:czhong2}
 	   	 		  \end{figure}

   \textbf{Case 1 } Choose suitable $c$ such that $c>\frac{1}{r(1-a)}$, $ U_1(1^+)<S_a(1^+)$, and   $p_H>1$.  From theorem \ref{heterostrong}, there is a loop of heteroclinic orbits when $p=p^{\#}$. We have the following conclusion.
   \begin{theorem} \label{case2dyn}Suppose that $c>\frac{1}{r(1-a)}$,  $ U_1(1^+)<S_a(1^+)$, and   $p_H>1$. Assume that the first Lyapunov  coefficient $a(p_H)<0$, and  every periodic orbit of system (\ref{ODE system}) is  orbitally stable,  then $p_H< p^{\#}  $.\\
   $~~~~~~$(i) If $p_{SN}<p<1$, the upward $\Gamma_p^u(E_R^*)$ connects $E_R^*$ to $E^*$. The orbits through any point above $\Gamma_p^s(E_a)$ converge to $E_0$, and the orbits   through any point inside the  stable manifold $\Gamma_p^s(E_R^*)$ converge to $E^*$. The orbits through any point below  $\Gamma_p^s(E_a)$   and exterior to  $\Gamma_p^s(E_R^*)$ converge to $E_1$ (see Fig. \ref{fig:czhong2} a) ).\\
   	   $~~~~~~$(ii) If $1< p <p_H $, $\Gamma_p^u(E_1)$ connects $E_1$ to $E^*$. The orbits through any point above $\Gamma_p^s(E_a)$ converge to $E_0$, and the orbits   through any point below $\Gamma_p^s(E_a)$ converge to $E^*$.	 	  \\
   	    $~~~~~~$(iii) When  $p_H < p <p^{\#} $,  $E^*$ is a repellor, and there is a unique limit cycle under $\Gamma_p^s(E_a)$. The orbits   through any point below $\Gamma_p^s(E_a)$ converge to the  limit cycle.  (see Fig. \ref{fig:czhong2}  b), c)).  \\
    $~~~~~~$	(iv) When $p=p^{\#} $, $\Gamma_{p^{\#}}^s(E_a)=\Gamma_{p^{\#}}^u(E_1)$, and there are two  heteroclinic orbits forming a loop  of heteroclinic orbits between $E_1$ and $E_a$ (see Fig. \ref{fig:czhong2}  d)). The orbits through any point exterior to the cycle  converge to $E_0$, and the orbits through any point interior to the cycle  converge to the cycle.\\
   	 $~~~~~~$(v) If $p^{\#}<p<\frac{1}{a}$, $\Gamma_p^s$ connects $E^*$ to $E_a$, and the extinction equilibrium $E_0(0,0)$ is globally asymptotically stable (phase portrait is similar as Fig. \ref{fig:czhong1}  b)).       	
   \end{theorem}
  \noindent Proof. From proposition \ref{manifoldstrong}, $S_a(p)$ is decreasing with $p$, combining with  $ U_1(1^+)<S_a(1^+)$,  then $\Gamma_p^s(E_a)$ enters $E_a$ from the region $\{(u,v):v>f(u)\}$   for all $p<1$.  Since $E^*$ is stable when $p<p_H$, then $E^*$ is stable for all $p_{SN}<p<1$. If there is a periodic orbit below $\Gamma_p^s(E_a)$, it must encircle $E^*$, and it is unstable from inside, which contradicts with our assumption of  the  orbital stability of periodic orbit. It means that there is no periodic orbit under $\Gamma_p^s(E_a)$ for all $p<1$. From proposition \ref{omegaalphastrong}, $E^*$ is the $\omega-$limit set of $\Omega_{22}$, and  $\Gamma_p^s(E_R^*)$ must be above $\Gamma_p^u(E_R^*)$.  The upward $ \Gamma_p^u(E_R^*) $ connects $E_R^*$ to $E^*$	
   (see Fig. \ref{fig:czhong2} a) ).   Thus, the orbits   through any point inside the  stable manifold $\Gamma_p^s(E_R^*)$ converge to $E^*$.
  The proof of (ii)-(v) is similar as that of theorem \ref{cxiaoyudynamics}.          $~~~~~~~~~~~~~~~~~~~~~~~~~~~~~~~~~~~~~~~~~~~~~~~~~~~~~~~~~~~~~~~~~~~~~~~~~~~~~~~~~~~~~~~~~~~~~~~~~~~~~~~~~~~~~~~~~~~~~~~~~~~~~~~~~~~~~~~~~~~~~~~~~~~~~~~~~~~~~\Box$
   	
\begin{figure}[htb]\centering
 	  	a)\includegraphics[width=0.3\textwidth]{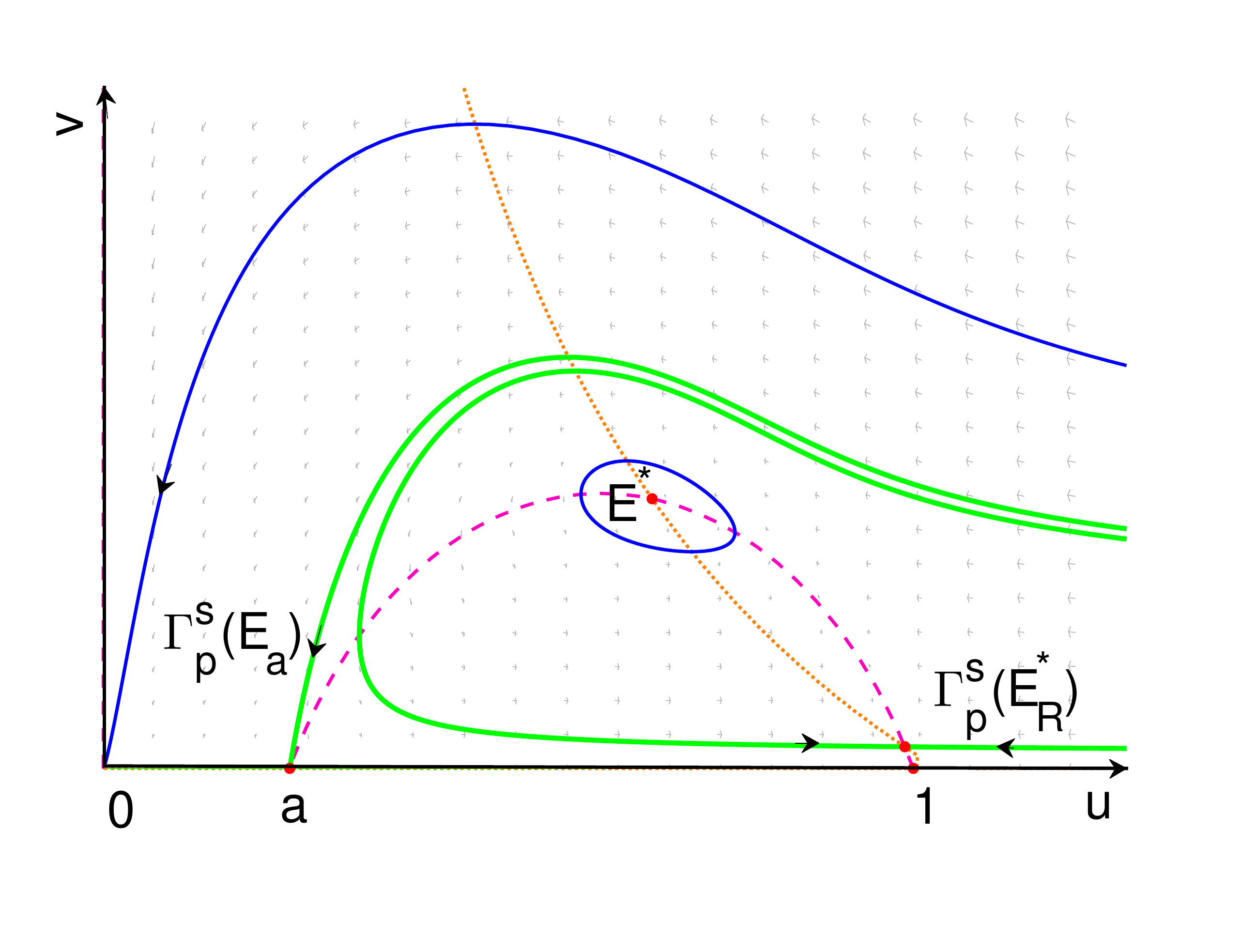}
  	      b)\includegraphics[width=0.3\textwidth]{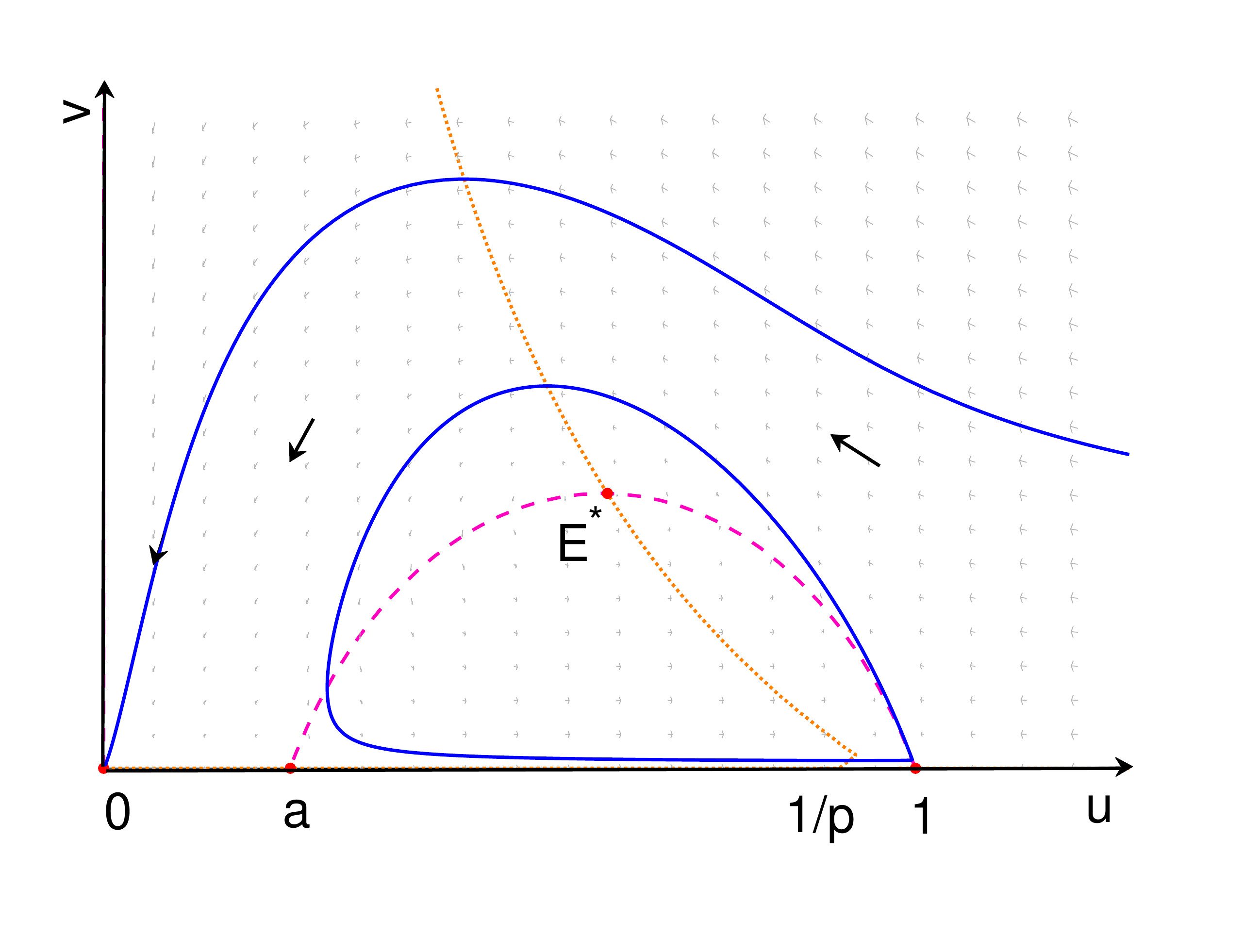}
 	   		c)\includegraphics[width=0.3\textwidth]{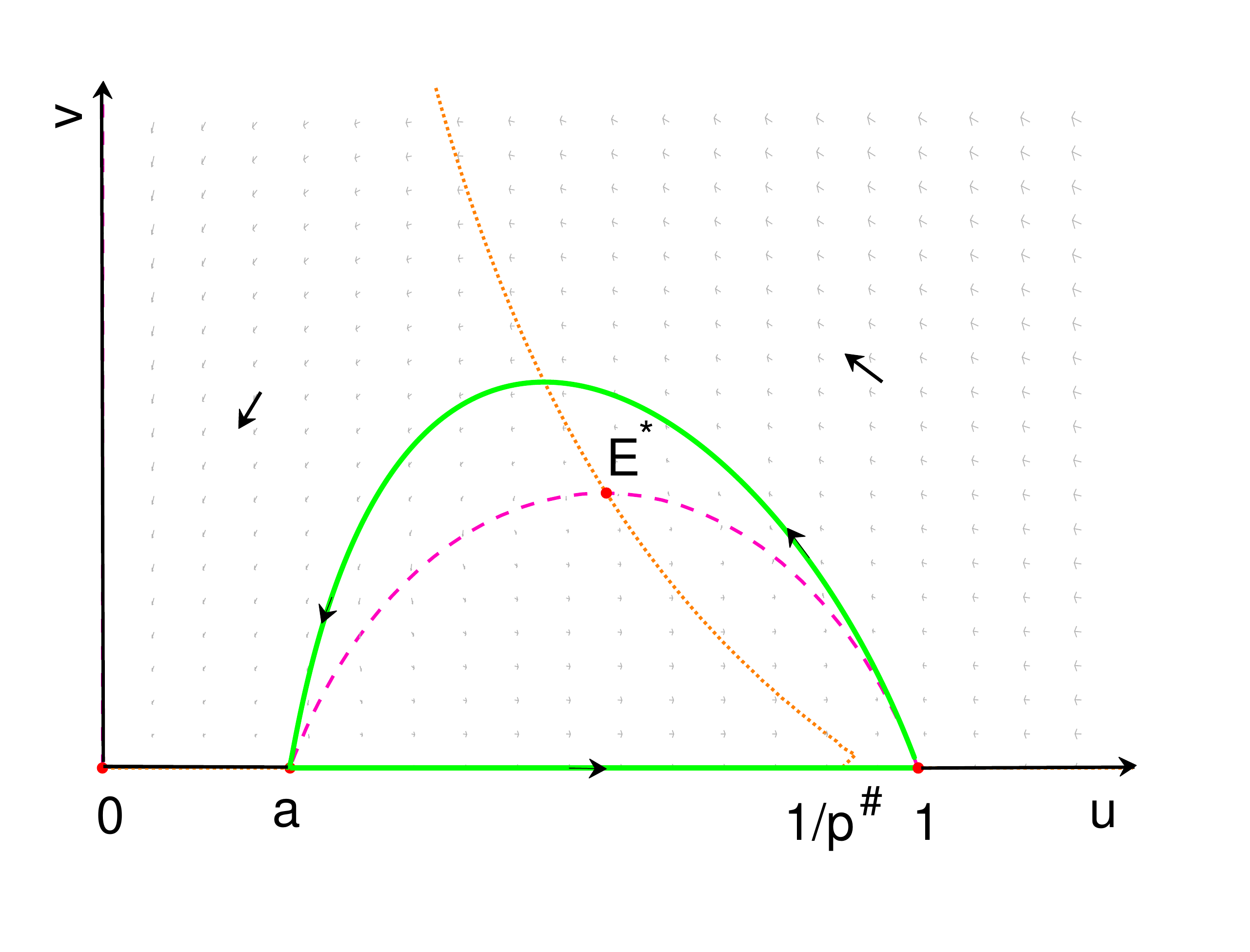}    \caption{  Phase portrait of system (\ref{ODE system})  for case $2$   when a) $p_H<p<1$ and close to $p_H$;     b) $p_H<p<p^{\#}$ and  close to $p^{\#}$;   c) $p=p^{\#}$.	 }	   	  	\label{fig:cda5}   	 		  \end{figure}

 	  \textbf{Case 2 }  Choose suitable $c$ such that   $c>\frac{1}{r(1-a)}$,
 	   	    	    $U_1(1^+)<S_a(1^+)$, and  $p_H<1$.  	From Theorem
 	    	\ref{heterostrong},  there is a loop of heteroclinic orbits when $p=p^{\#}$.
 	   \begin{theorem}    Choose suitable $c$ such that   $c>\frac{1}{r(1-a)}$,
 	    	    $U_1(1^+)<S_a(1^+)$, and  $p_H<1$. Assume that the first Lyapunov coefficient $a(p_H)<0$, and  every periodic orbit of system (\ref{ODE system}) is orbitally stable,  then $p_H< p^{\#}  $.\\
  	 	$~~~~~~$	(i) If $p_{SN}<p<p_H$, the upward $\Gamma_p^u(E_R^*)$ connects $E_R^*$ to $E^*$, and the downward $\Gamma_p^u(E_R^*)$ connects $E_R^*$ to  $E_1$. The orbits through any point above $\Gamma_p^s(E_a)$ converge to $E_0$, and the orbits   through any point inside the  stable manifold $\Gamma_p^s(E_R^*)$ converge to $E^*$. The orbits through any point below  $\Gamma_p^s(E_a)$   and exterior to  $\Gamma_p^s(E_R^*)$ converge to $E_1$  (phase portrait is similar as Fig. \ref{fig:czhong2}  a)  ). \\
  	   	  $~~~~~~$  (ii) When  $ p_H <p<1$,   there is a unique limit cycle inside the stable manifold $\Gamma_p^s(E_R^*)$ of $E_R^*$. (see Fig. \ref{fig:cda5}  a)  ).   	
  	   	  The orbits through any point interior to  $\Gamma_p^s(E_R^*)$ converge to the limit cycle. 	\\
  	  $~~~~~~$  (iii) When  $1 < p <p^{\#} $,   there is a unique limit cycle under $\Gamma_p^s$. The orbits   through any point below $\Gamma_p^s$ converge to the  limit cycle.  (see Fig. \ref{fig:cda5}  b)  ).  \\
 	   $~~~~~~$ (iv) When $p=p^{\#} $, $\Gamma_{p^{\#}}^s(E_a)=\Gamma_{p^{\#}}^u(E_1)$, there are two  heteroclinic orbits forming a loop of heteroclinic orbits from $E_1$ to $E_a$ and back to $E_1$ (see Fig. \ref{fig:cda5}  c)). The orbits through any point exterior to the cycle  converge to $E_0$, and the orbits through any point interior to the cycle  converge to the cycle.  \\	
 $~~~~~~$(v) If $p^{\#}<p<\frac{1}{a}$, $\Gamma_p^s$ connects $E^*$ to $E_a$, and the extinction equilibrium $E_0(0,0)$ is globally asymptotically stable (phase portrait is similar as Fig. \ref{fig:czhong1}  b)).   		
\end{theorem}
 \noindent  Proof. We can prove (i)   similar as in Theorem \ref{case2dyn}.
 (ii) Since $a(p_H)<0$ and $\alpha'(p_H)>0$, there is a stable periodic orbit bifurcating from Hopf bifurcation when  $p> p_H $.
 Similar as the proof in Theorem \ref{case2dyn}, we can prove that $\Gamma_p^s(E_R^*)$ must be above $\Gamma_p^u(E_R^*)$.
  From proposition \ref{omegaalphastrong}, the unique limit cycle is  the $\omega-$limit set of $\Omega_{22}$.  The proof of (iii)-(v) is similar as that of Theorem \ref{cxiaoyudynamics}.
$~~~~~~~~~~~~~~~~~~~~~~~~~~~~~~~~~~~~~~~~~~~~~~~~~~~~~~~~~~~~~~~~~~~~~~~~~~~~~~~~~\Box$

 	 	\begin{figure}[htb]\centering
 	    a)\includegraphics[width=0.3\textwidth]{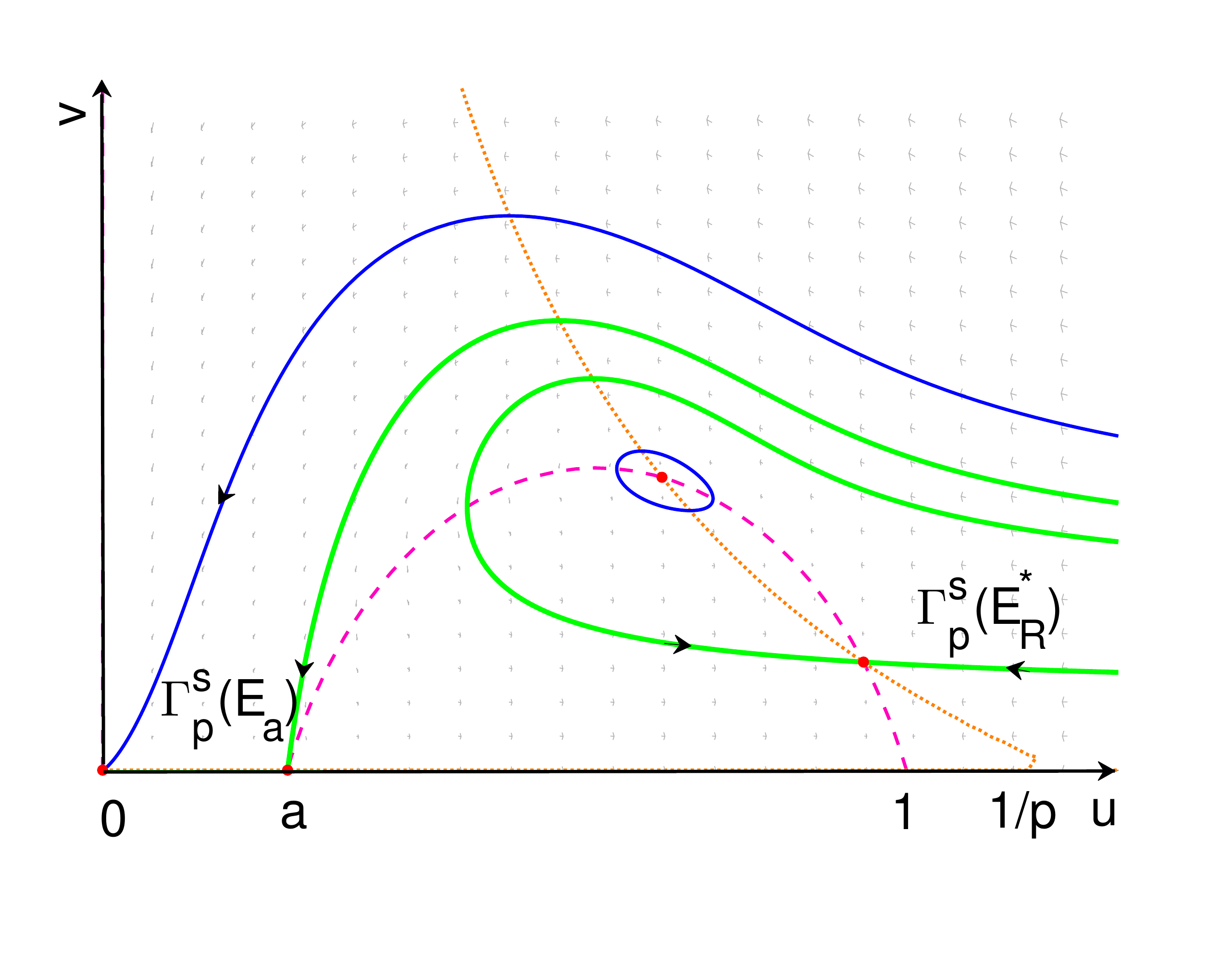}
 	   		b)\includegraphics[width=0.3\textwidth]{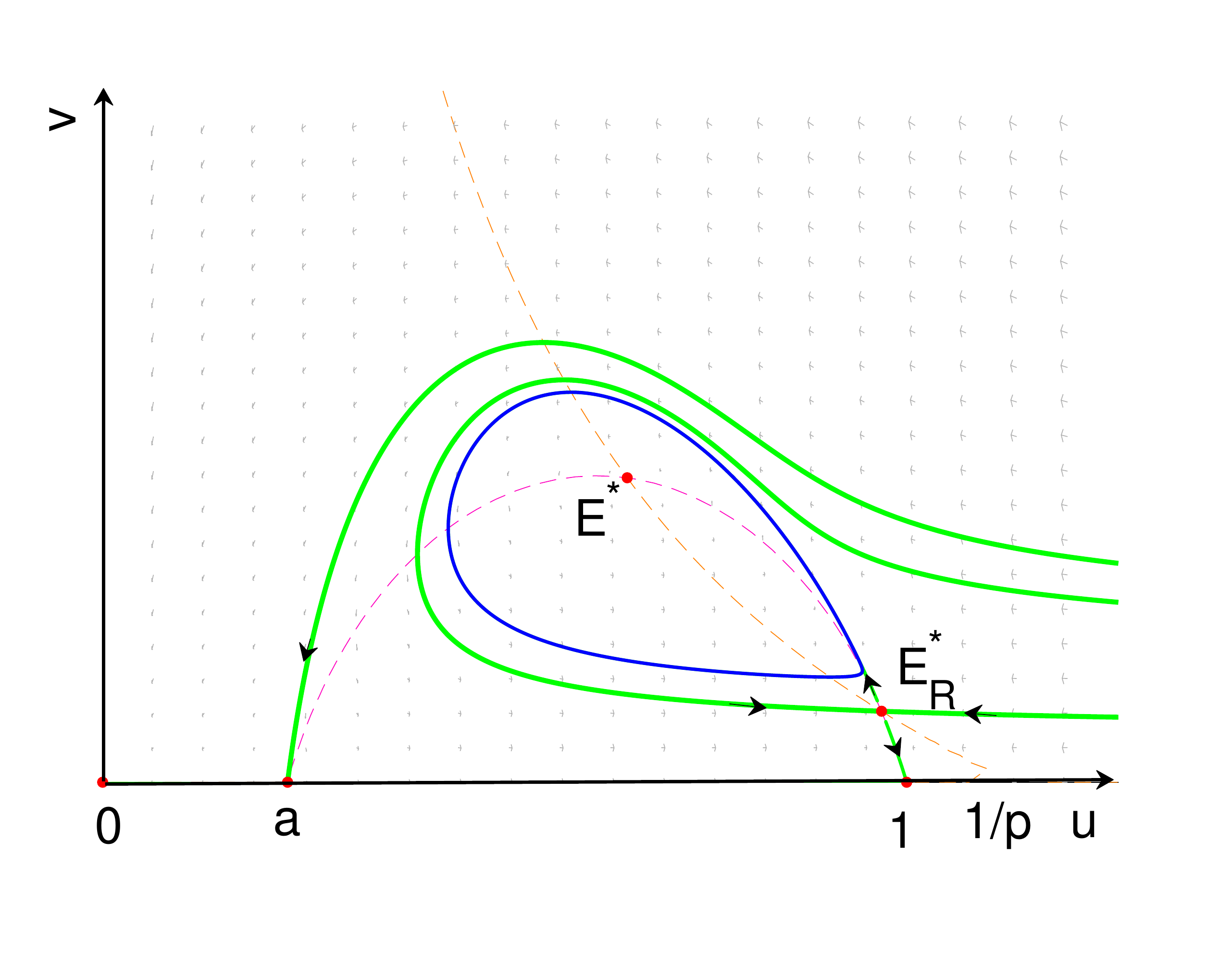}
 	   	 		c)\includegraphics[width=0.3\textwidth]{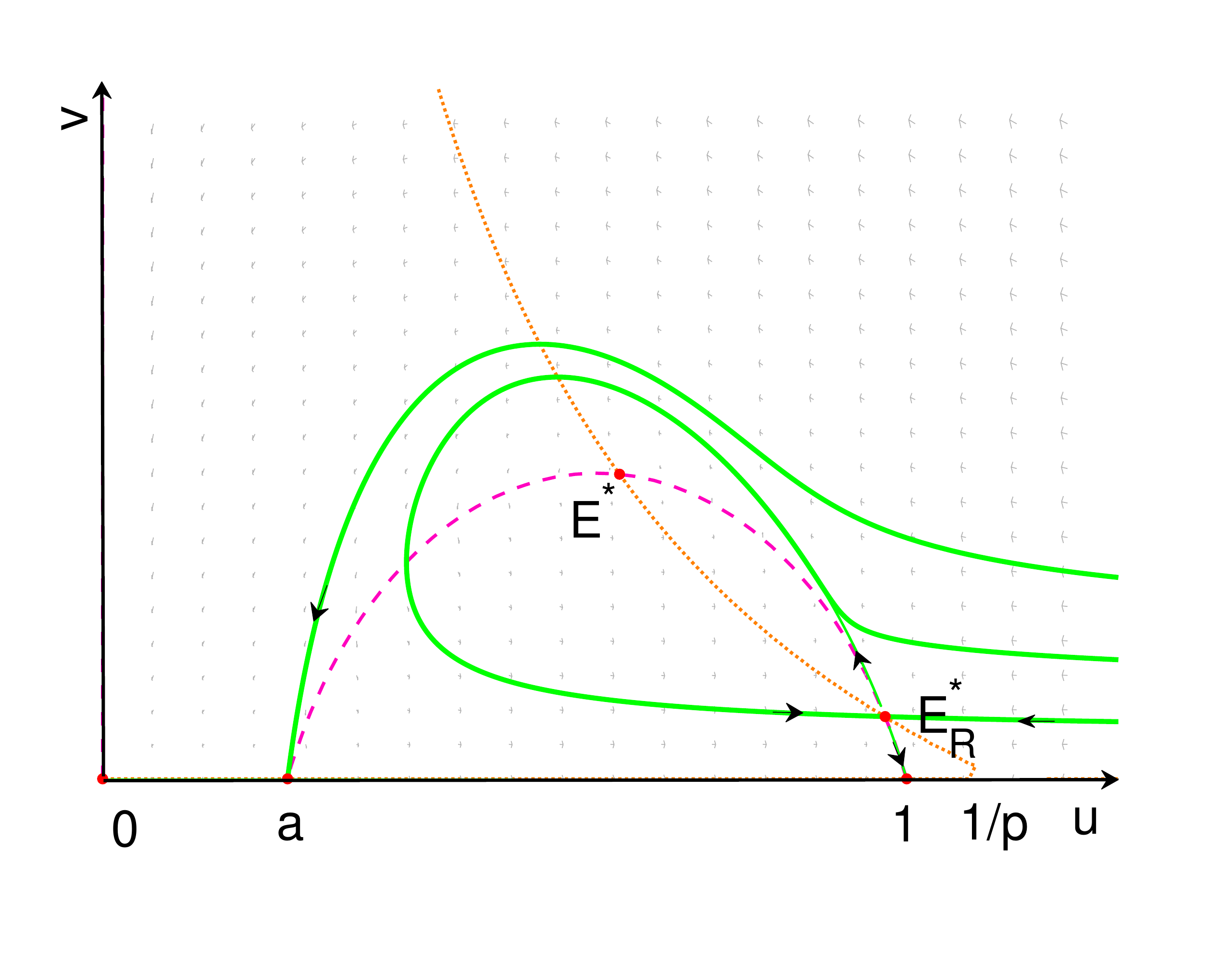}
 	   		d)\includegraphics[width=0.3\textwidth]{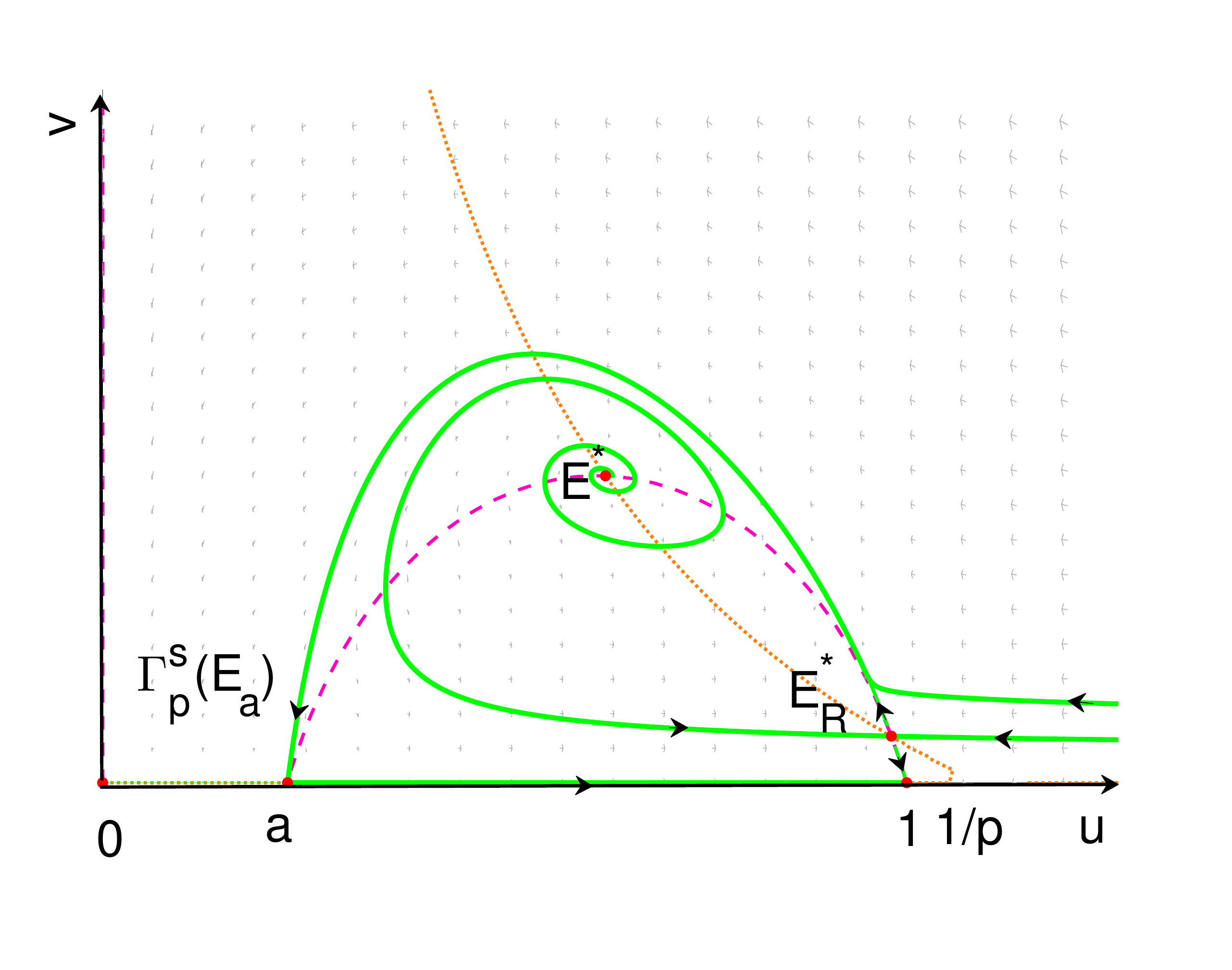}
 	   	e)\includegraphics[width=0.3\textwidth]{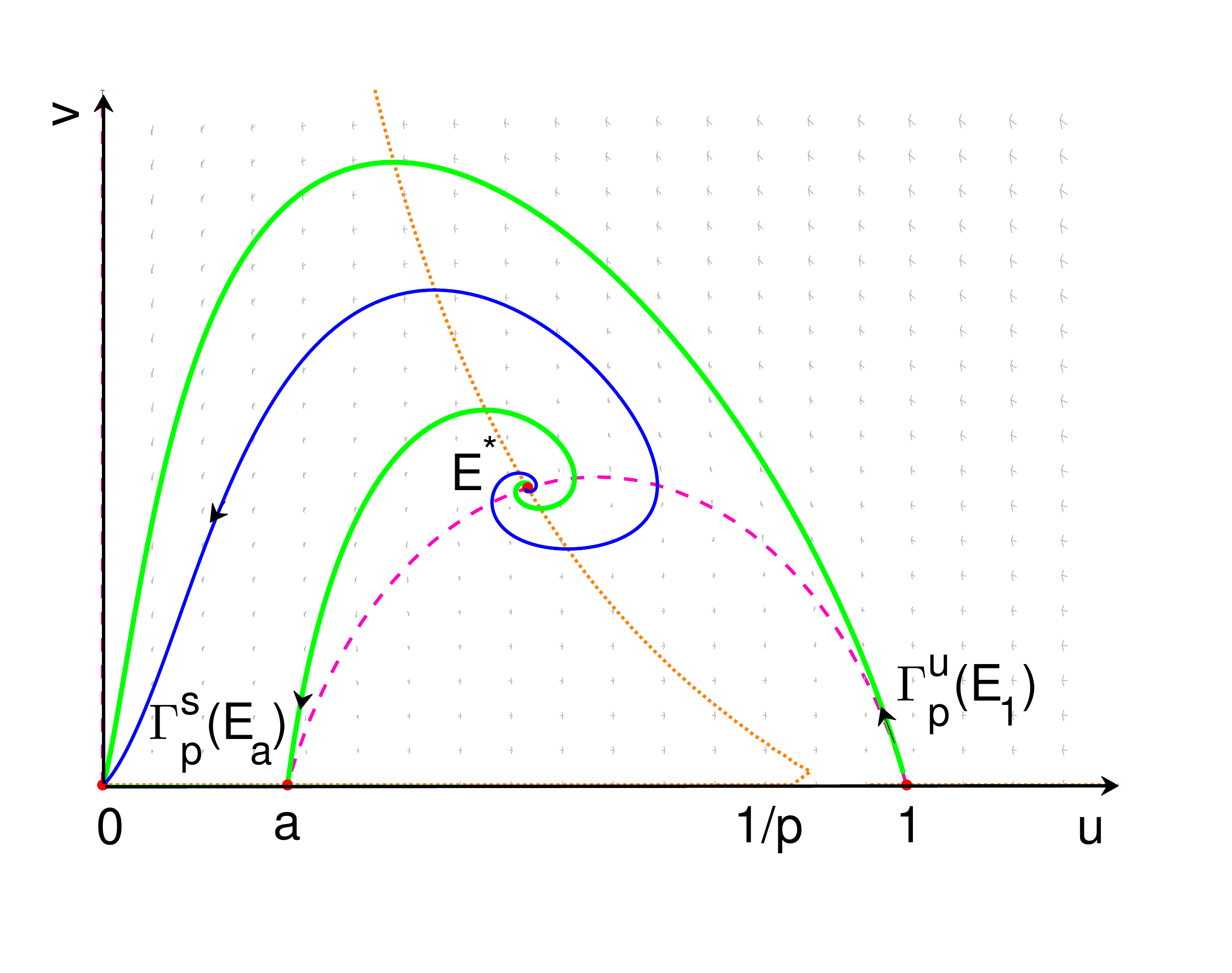}
 	   	   	 	\caption{  Phase portrait of system (\ref{ODE system}) for case $3$  when 	a) When $p>p_H$  and near $p_H$, there is a  limit cycle inside the stable manifold of $E_R^*$; b) When $p<p_{hom}$ and near $p_{hom}$, there is a limit cycle inside the stable manifold of $E_R^*$;  c) When $p=p_{hom}$,  the upward unstable manifold  and the left  stable manifold of $E_R^*$ collide,  which forms a homoclinic cycle;
 	   	 d)   The upward unstable manifold of $E_R^*$, the unstable manifold  of $E_a$ on the $u-$axis, and the downward  unstable manifold of $E_R^*$ forms a loop of heteroclinic orbits  among $E_R^*$, $E_a$ and $E_1$ when $p=\overline{p}^{\#}$;
 	   	 	e) $1<p<\frac{1}{a}$. 
 	   	  	   	  }	   	  	  \label{fig:c8}
 	   	  \end{figure}

     \textbf{Case 3 }  Choose $c$  such that $c>\frac{1}{r(1-a)}$ and   $U_1(1^+)>S_a(1^+)$. 	From Theorem
      	    	\ref{heterostrong}, there is no heteroclinic orbits from $E_1$ to $E_a$. It is easy to prove that  when $1<p<\frac{1}{a}$, $E_0$ is globally asymptotically stable (see Fig. \ref{fig:c8} e)).

      \begin{remark} 
       In this case, the location of $\Gamma_p^s(E_a)$, $\Gamma_p^s(E_R^*)$ and $\Gamma_p^u(E_R^*)$ are very complicated. We can observe the following dynamics  by numerical simulations.  (i) Loop of heteroclinic orbits.
     The downward  branch of the unstable manifold of $E_R^*$ connects  $E_R^*$ to $E_1$, and the upward connects $E_R^*$ to $E_a$, which  collides with the  the stable manifold of $E_a$.
      The  upward and downward unstable manifold of $E_R^*$, together with the unstable manifold  of $E_a$ on the $u-$ axis,  forms a loop of heteroclinic orbits among $E_R^*$, $E_a$ and $E_1$ when $p=\overline{p}^{\#}$ (see Fig. \ref{fig:c8} d)).  (ii) Homoclinic cycle. The upward unstable manifold  and the left stable manifold of $E_R^*$ collide,   which forms a homoclinic cycle. Denote the parameter $p$ as $p=p_{hom}$ (see Fig. \ref{fig:c8} c)).
      (iii) Limit cycle induced by Hopf bifurcation (see Fig. \ref{fig:c8} a), b)).
  \end{remark}

   \subsection{Numerical simulations for ODE system}
   \label{simulationsode}

   In this section, we carry out some mumerical simulations for  system (\ref{ODE system}).  We fix the parameters as   \begin{equation}\label{paraodesimu}
          a=0.23,r=1.1, m=0.31.
           \end{equation}

   \subsubsection{Numerical simulations for system (\ref{ODE system}) with weak cooperative hunting}
\begin{figure}
          \centering
       a) \includegraphics[width=0.46 \textwidth]{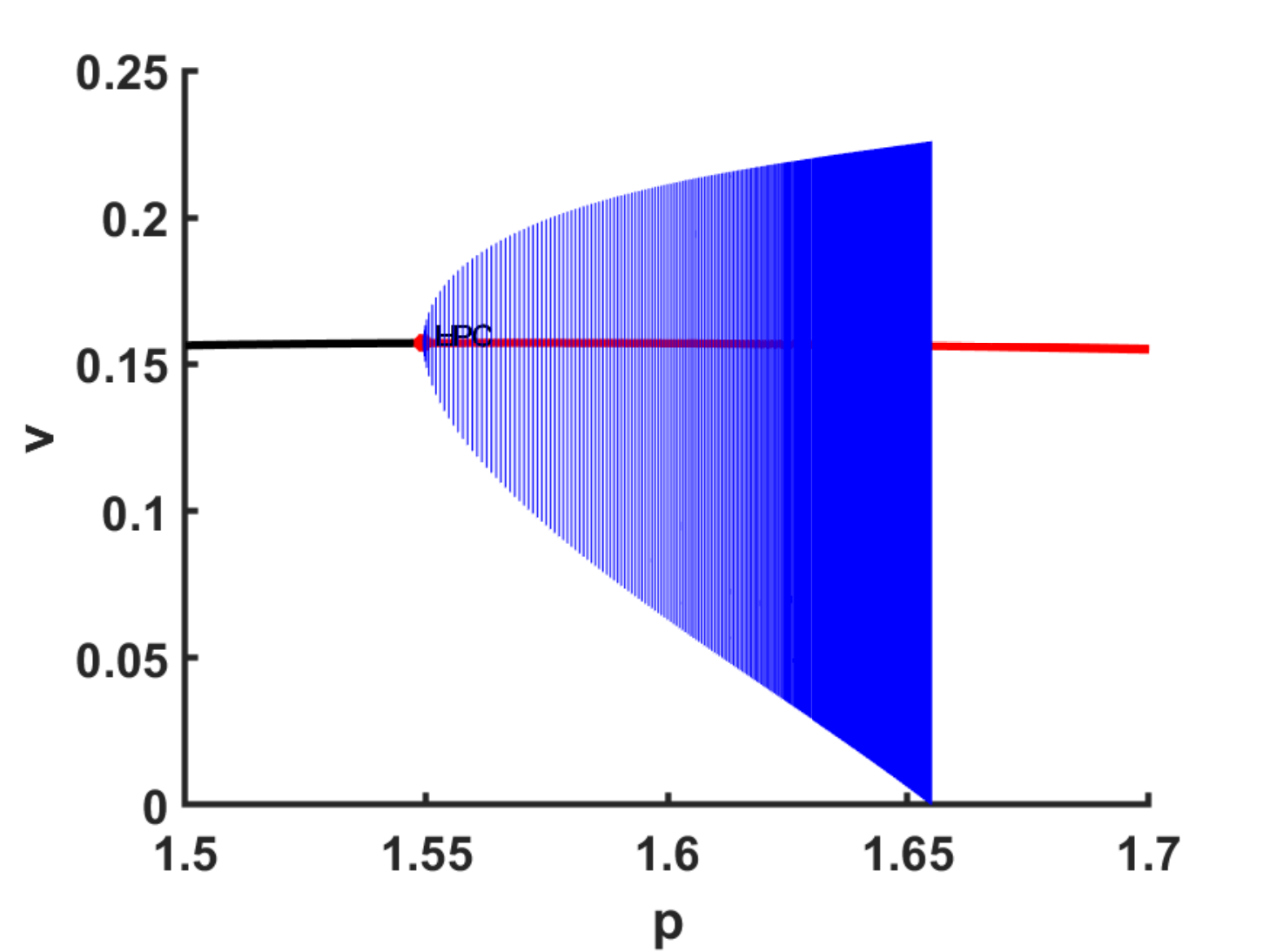}
         b)   \centering
       \includegraphics[width=0.46 \textwidth]{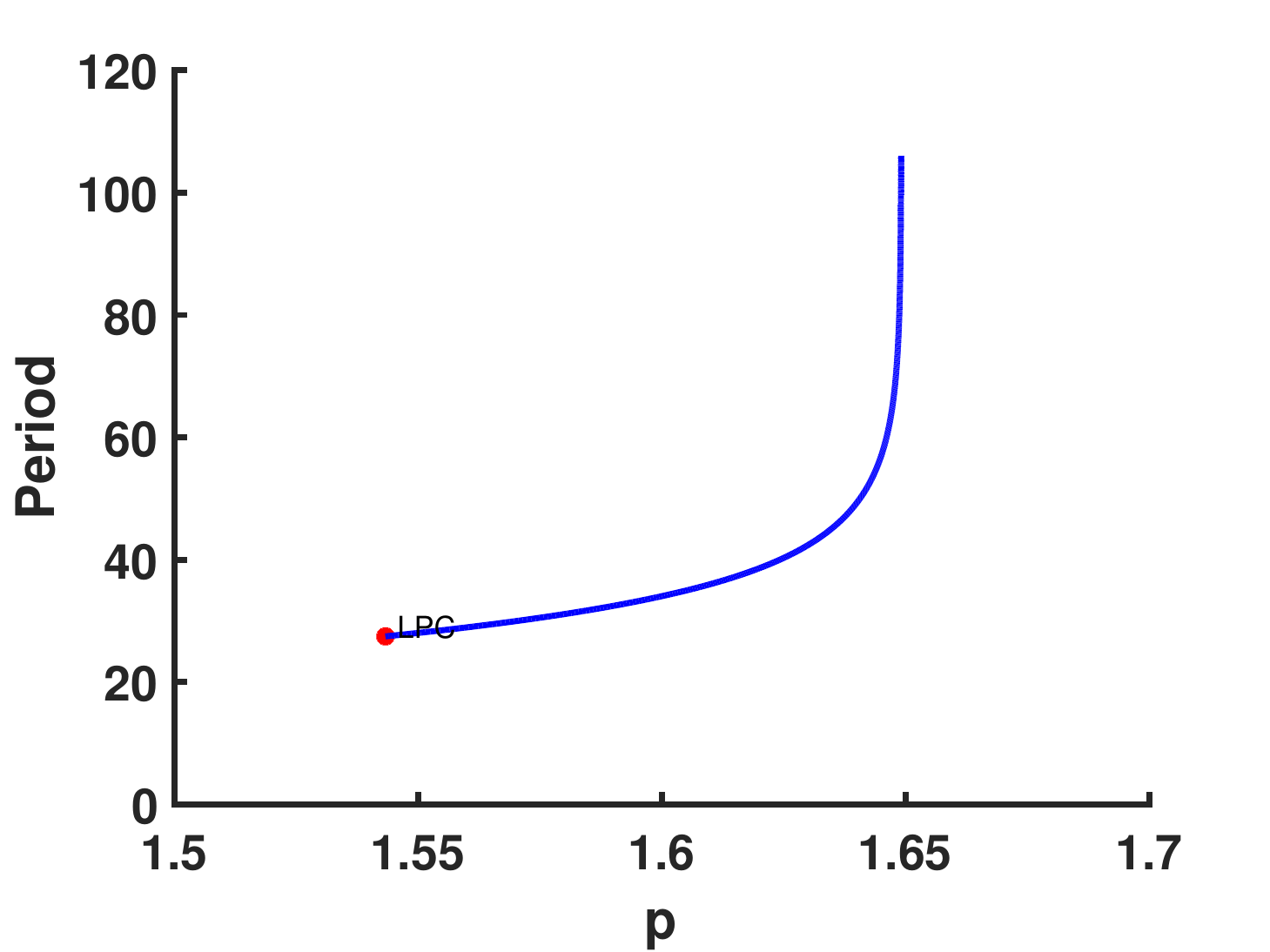}
        \caption{ The bifurcation diagram  with  $ a=0.23,m=0.31,r=1.1$  and  $c=0.25 $, and the conversion rate $p$ as the bifurcation
        parameter.  a)   The black curve is the stable steady state and the   red curve is the unstable steady state. The blue
        curves denote  the   $v-$amplitude of the periodic orbits,  which begin  at the Hopf bifurcation point $p=p_H=1.5432$. b) The period of the limit cycles.}
       \label{fig:c025bifurcationHopf}
      \end{figure}

   For  system (\ref{ODE system}),
    choose $c=0.25$  such that $c<\frac{1}{r(1-a)}$,  and vary the parameter $p$. Using the method in \cite{Dhooge A}, we can get  Hopf bifurcation point $p_H=1.5432$, and the first Lyapunov coefficient is $-1.1211$, which means that the bifurcating periodic solution is asymptotically stable, and it is bifurcating from $E^*$ as $p$ increases past $p_H$ from Theorem \ref{Hopfdirection}. We can also draw the bifurcation diagram in Fig. \ref{fig:c025bifurcationHopf}, which shows that as $p$ increases from $p_H$, the period of limit cycle is increasing, and it  tends to infinite as $p\rightarrow 1.6491$.
   Moreover, the amplitude of oscillation increases as $p$ increases, and the minimum of the predator population is close to zero when $p$ is close to $p^{\#}=1.6491$, which will increase the risk on the extinction of predator.

     To show the  complex dynamics of system (\ref{ODE system}),  different values of $p$  are chosen, listed in Table 2.  The corresponding phase portraits    for different values of $p$ have been illustrated  in Figs. \ref  {fig:ode01equipxiaoyu1} and
    	\ref{fig:ode01equi},  which is drawn by pplane8 \cite{Polking JA}.
    \begin{table}[tbp]\label{ppp}
    \label{parameterp}
    \caption{Values of parameter $p$ chosen in Figs. \ref  {fig:ode01equipxiaoyu1} and
     	\ref{fig:ode01equi}}
    \centering
    \begin{tabular}{ccccccccc}
    \hline
    Figure &  Fig. \ref  {fig:ode01equipxiaoyu1} a& Fig. \ref  {fig:ode01equipxiaoyu1} b  & Fig. \ref  {fig:ode01equipxiaoyu1} c & Fig. \ref  {fig:ode01equipxiaoyu1} d& Fig. 	 \ref{fig:ode01equi} a & Fig. 	\ref{fig:ode01equi} b  & Fig. 	\ref{fig:ode01equi} c& Fig. 	\ref{fig:ode01equi} d  \\ \hline  
    $p$ &  5.1 & 3.2&  1.2 & 0.9 & 1.55 & 1.645& 1.6491& 1.95\\       \hline
    \end{tabular}
    \end{table}

   \subsubsection{Numerical simulations for ODE system with strong cooperative hunting}
   In section \ref{section3cases}, we discuss three different cases  of dynamics when    $c$ is chosen as different values.
   To show case $1$, we fix $c=3$, and vary $p$ as a bifurcation parameter.  We get $p_H=1.1024$, and the first Lyapunov coefficient is $-2.3280$.   Choosing different values of $p$ listed in Table 3,  we have drawn the corresponding phase portraits    for different values of $p$   in Fig. \ref{fig:czhong2}.
    \begin{table}[tbp]\label{pcase12}
                \caption{Values of parameter $p$ chosen for showing case $1$ in Figs. \ref  {fig:czhong2} and \ref  {fig:cda5}}
           \centering
           \begin{tabular}{ccccc}
           \hline
           Figure &  Fig. \ref  {fig:czhong2}  a& Fig. \ref  {fig:czhong2}  b  & Fig. \ref  {fig:czhong2}  c  & Fig. \ref  {fig:czhong2}  d \\ \hline  
           $p$ & 0.96  &1.105  &  1.2& 1.2068\\       \hline
            Figure &  Fig. \ref  {fig:cda5}  a& Fig. \ref  {fig:cda5}  b  & Fig. \ref  {fig:cda5}  c &  \\ \hline  
                         $p$ &  0.965 & 1.052&1.05665& \\       \hline
           \end{tabular}
           \end{table}

   To  show  case $2$,  we choose   $c=5$, and vary $p$, we can get $p_H=0.9608$, and the first Lyapunov coefficient is $-2.9553$.   Choosing $p$ as listed in Table 3, the corresponding phase portraits   have been illustrated in Fig. \ref{fig:cda5}.

\begin{figure}
           \centering
        a) \includegraphics[width=0.46 \textwidth]{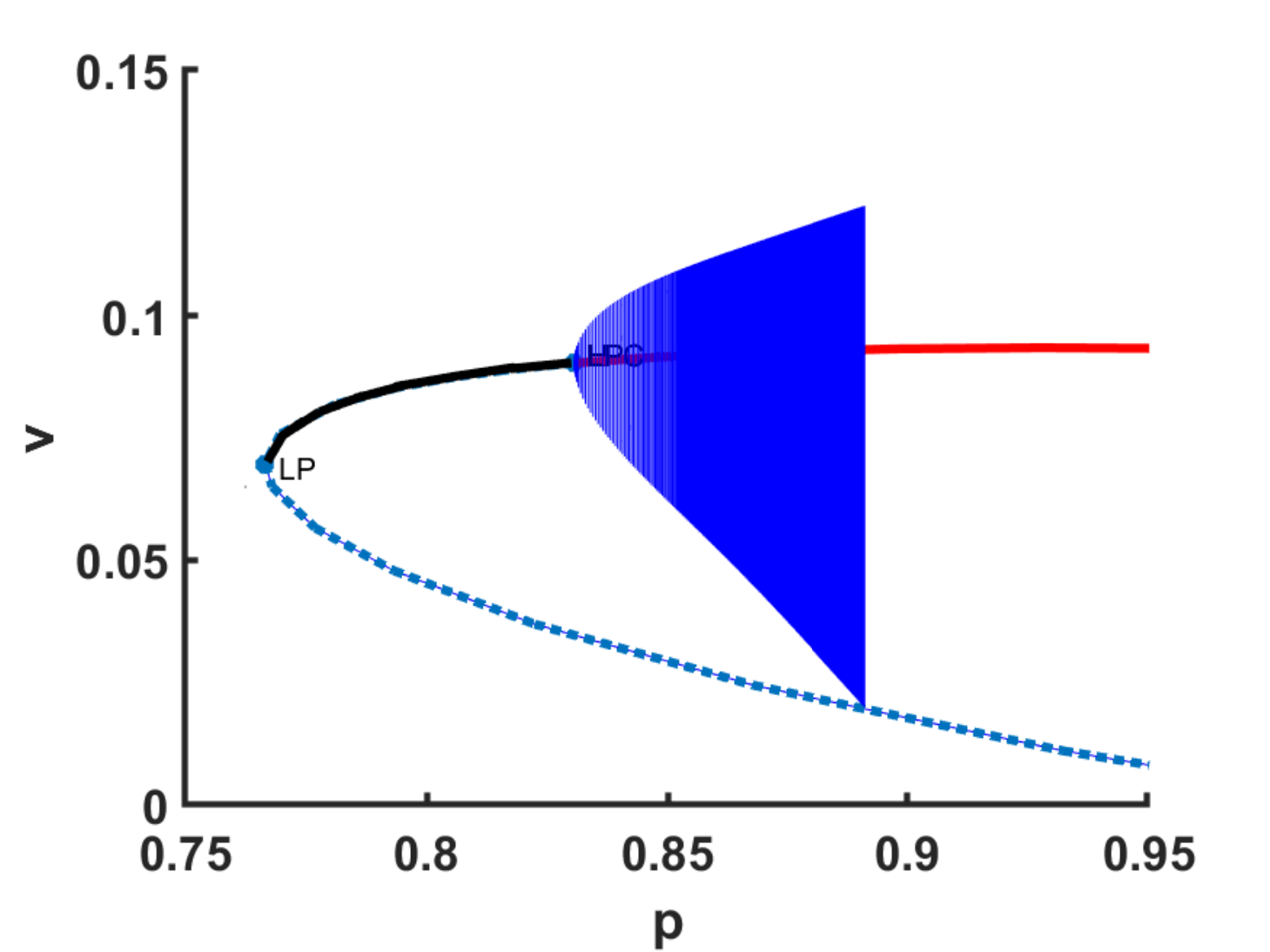}
          b)   \centering
        \includegraphics[width=0.46 \textwidth]{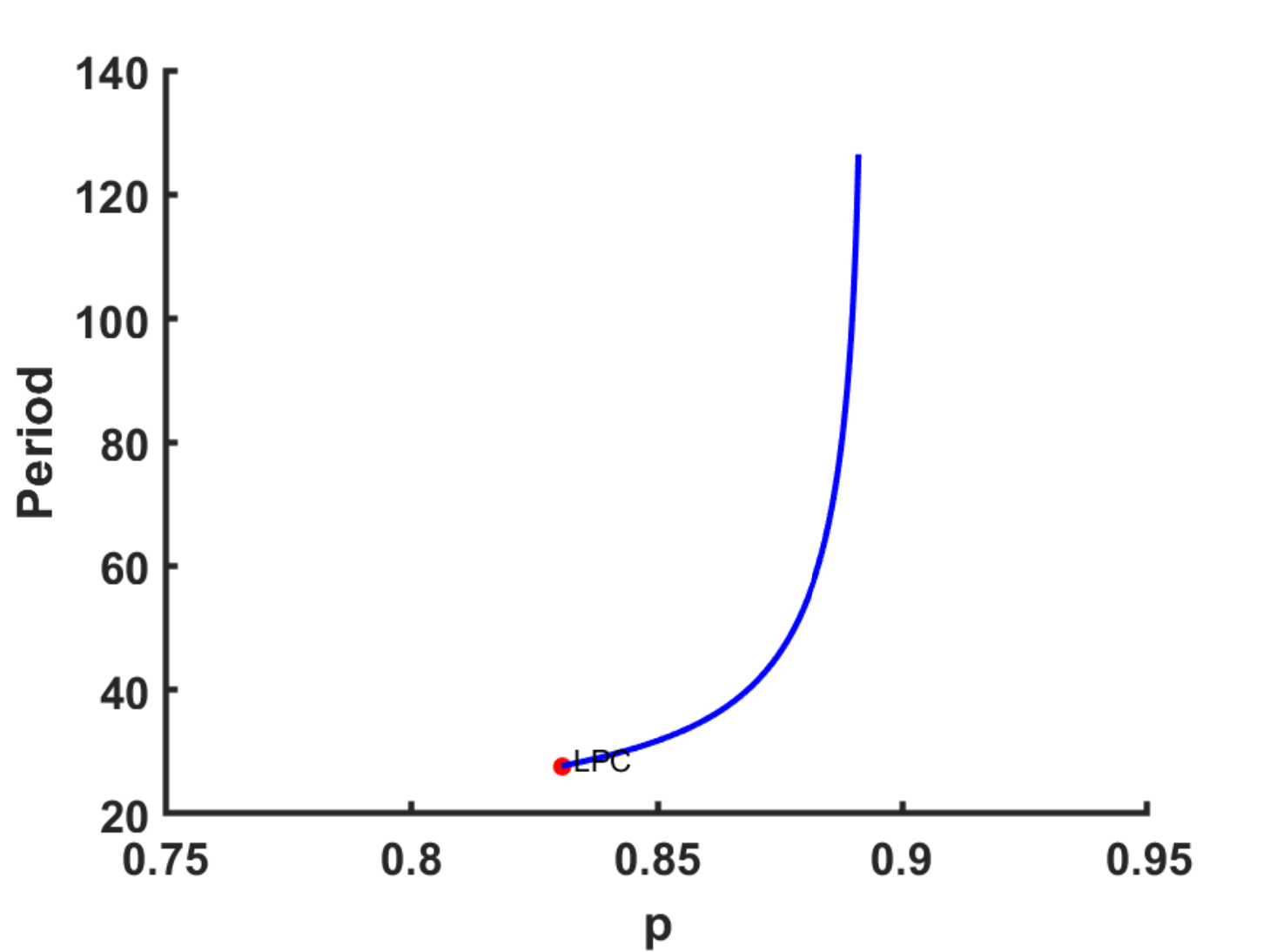}
         \caption{  The bifurcation diagram  with  $ a=0.23, m=0.31,r=1.1$  and  $c=8$, and the conversion rate $p$ as the bifurcation parameter.    a) The black curve is the stable steady state $E^*$, the   red curve is the unstable steady state  $E^*$, and the dotted blue curve is the unstable steady state  $E_R^*$. The blue  curves denote  the   $v-$amplitude of the periodic orbits,  which begin  at the Hopf bifurcation point $p=p_H=0.8306$. b) The period of the limit cycles.}
        \label{fig:c8bifurcationHopf}
       \end{figure}

   To  show  case $3$,  we choose   $c=8$, and vary $p$, we can get $p_H=0.8306$, and the first Lyapunov coefficient is $-3.732$. We can draw the bifurcation diagram  in Fig. \ref{fig:c8bifurcationHopf},  in which the curve above is the $v-$value of interior equilibrium $E^*$, and the curve below is the $v-$value of $E_R^*$.
   We can see that as $p$ increasing from $p_H$, the period of limit cycle is increasing, and it  tends to infinite as $p\rightarrow 0.8919$.
      Moreover, the amplitude of oscillation increases as $p$ increases, and the amplitude line touches the curve of $v-$value of $E_R^*$, which means that there is a cycle goes through $E_R^*$. In fact, there is a homoclinic cycle when $p=p_{hom}=0.8919$ (see  Fig. \ref  {fig:c8}  c)).
    Choosing $p$ as listed in Table 4,   the corresponding phase portraits   have been illustrated in Fig. \ref{fig:c8}.

     \begin{table}[tbp]\label{pcase3}
                   \caption{Values of parameter $p$ chosen for showing case $3$ in Fig. \ref {fig:c8}}
              \centering
              \begin{tabular}{cccccc}
              \hline
              Figure &  Fig. \ref  {fig:c8}  a& Fig. \ref {fig:c8}  b  & Fig. \ref  {fig:c8}  c &  Fig. \ref  {fig:c8}  d& Fig. \ref {fig:c8}  e    \\ \hline  
              $p$ & 0.834  & 0.88& 0.8919&  0.91564& 1.1 \\       \hline
              \end{tabular}
              \end{table}

    \section{ Diffusion-driven Turing instability and Turing-Hopf bifurcation}\label{sectionTuring}

     In the previous section, we have discussed both local and global dynamics of ODE  model (\ref{ODE system}).  In fact, the preys and predators distribute inhomogeneously in different locations, and the spatial diffusion plays an important part in the process of population evolution.  Turing instability and Turing-Hopf bifurcation induced by  diffusion have been widely investigated recently (see \cite{YongliSongTuringHopf,X. Xu,Q. An,guogaihui}).
     Taking into account the diffusion in system (\ref{ODE system}), we consider the following diffusive predator-prey system:
    \begin{equation}
    \label{diffusion}
    \left\{
    \begin{array}{l}
    \dfrac{\partial u(x,t)} {\partial t}= d_1\Delta u(x,t)+ru(x,t)[1-u(x,t)] (u(x,t)-a), \\~~~~~~~~~~~ -(1+cv(x,t))u(x,t)v(x,t),~~~~~~~~~~~~~~~~~~~~~~~~~~~~~~~~~~~~~~~~~~~~~~~~~~~~x\in (0,l\pi),\\
    \dfrac{\partial v(x,t)}{\partial t }= d_2\Delta v(x,t)+mv(x,t)\left[  p u(x,t)\left( 1+cv(x,t)\right) -1 \right]  , ~x\in (0,l\pi) \\
    \partial_x u(x,t)= 0,~~ \partial_x v(x,t) =0,   ~~~~~~~~~~~~~~~~~~~~~~~~~~~~~~~~~~~~~~~~~~~~~~~~~~~~~~~~~~~~~~~~x=0, l\pi,\\
    \end{array}
    \right.
    \end{equation}
   where  $d_1,d_2>0$ are the diffusion coefficients characterizing the rates of the spatial dispersion of the prey and predator population, respectively.

   \subsection{ Turing instability and Turing-Hopf bifurcation induced by diffusion}\label{sectiondiffusion}

    In this subsection, we consider the effect of the diffusion on the stability of the constant steady state $E^*$.   If $E^*$ is linearly stable in the absence of diffusion, and it  becomes unstable in the presence of diffusion, we call such an instability Turing instability. Since $E_R^*$ is always a saddle if it exists, thus  Turing instability can only happen near $E^*$ in both cases: weak cooperation and strong cooperation.  We first investigate the existence of Turing instability, then we consider Turing-Hopf bifurcation near $E^*$.

For Neumann boundary condition, we define the real-valued Sobolev space
          \begin{equation}
          	X=\{(u, v)^T\in H^2(0,l\pi)\times H^2(0,l\pi),~ \dfrac{\partial u}{\partial x} =\dfrac{\partial v}{\partial x}= 0, {\rm at}~ x=0,l\pi\}.
          	\end{equation}
    The linearization of  system (\ref{diffusion})  at  the constant steady state $E^*(u^*,v^*)$ is given by
    \begin{equation}\label{}
    \left( \begin{array}{c}
    \frac{\partial u}{\partial t}   \\ \frac{\partial v}{\partial t}
    \end{array}\right) = D \left( \begin{array}{ c}
         \Delta u\\\Delta v
            \end{array}\right)  +A\left( \begin{array}{ c}
     u\\v
        \end{array}\right) \stackrel{\vartriangle}{=}L\left( \begin{array}{ c}
             u\\v
                \end{array}\right),
    \end{equation}
         where\begin{equation*}
  D = \left( \begin{array}{ cc}
   d_1   &0\\0& d_2
      \end{array}\right), ~ A=\left( \begin{array}{cc}
   ru^*(1+a-2u^*) &-2cu^*v^*-u^*\\
     mp(1+cv^*)v^*& mpcu^*v^*\\
    \end{array}\right)  .
       \end{equation*}

  It is well known that the eigenvalues  of $D\Delta$ on $X$ are
 $-d_1\frac{n^2}{l^2}$ and  $-d_2\frac{n^2}{l^2}$, $k\in \mathbb{N}_0=\{0,1,2,...\}$,
  with corresponding normalized eigenfunctions $\beta_n^{(1)}$ and  $\beta_n^{(2)}$, where
  \begin{equation*}
 \beta_n^{(1)}(x)=\left( \begin{array}{c }
\gamma_n\\0
 \end{array}\right),~~  \beta_n^{(2)}(x)=\left( \begin{array}{c }
 0\\ \gamma_n
  \end{array}\right),~~ \gamma_n(x)=\dfrac{\cos\frac{n}{l}x}{\parallel\cos\frac{n}{l}x\parallel_{L^2}}=\left\lbrace \begin{array}{ll}
  \sqrt{\frac{1}{l\pi}},& n=0,\\
  \sqrt{\frac{2}{l\pi}}\cos\frac{n}{l}x, & n\geq 1.
  \end{array}
  \right.\end{equation*}
 Applying the general theory about elliptic operators, we know that $\beta_n^{(1)}$ and  $\beta_n^{(2)}$   form an orthonormal basis for $X$.

From straightward calculation,  we obtain the characteidtic equations
    \begin{equation}
    \label{charactertau120}
    \Delta_n=\lambda^2-T_n\lambda+J_n=0,      ~~~~~~~n=0,1,2,...,
    \end{equation}
    where
    \begin{equation}\label{TnJn}
    \begin{aligned}
    &T_n
    =-\left( d_1\frac{n^2}{l^2}+d_2\frac{n^2}{l^2}+ru^*(2u^*-a-1)-mpcu^*v^*\right) =- d_1\frac{n^2}{l^2}-d_2\frac{n^2}{l^2}+{\rm tr}J_{E^*},\\
    &J_n
    =d_1d_2\frac{n^4}{l^4}+ru^*(2u^*-1-a)d_2\frac{n^2}{l^2}-mpcu^*v^*d_1\frac{n^2}{l^2}\\&~~~~~~+ru^*(a+1-2u^*)mpcu^*v^*+(2cu^*v^*+u^*)mv^*p(1+cv^*)\\&~~~~~~=d_1d_2\frac{n^4}{l^4}-ru^*(1+a-2u^*)d_2\frac{n^2}{l^2}-mpcu^*v^*d_1\frac{n^2}{l^2}+{\rm det}J_{E^*}.
    \end{aligned}
    \end{equation}

     From the previous section, in the absence of diffusion, i.e., $d_1=d_2=0$, $E^*$ is asymptotically stable if $1<p<p_H$ in the case of $c<\frac{1}{r(1-a)}$ (if $p_{SN}<p<p_H$ in the case of $c>\frac{1}{r(1-a)}$).  If there exists an $n\in\mathbb{N}$, such that  $\Delta_n=0$ has roots with positive real part when $p<p_H$, Turing instability occurs.
     Since ${\rm tr}J_{E^*}<0$ when $ p<p_H$,  and  ${\rm det}J_{E^*}>0$ always satisfies,  from (\ref{TnJn}), we have $T_n<0 $ for any $n$ when $p<p_H$.   Therefore, when $ p<p_H$,  the signs of the real parts of roots of (\ref{charactertau120}) are determined by the signs of $J_n$.

     In fact, the curve of $J_n=0$ is a hyperbola on $d_1-d_2$ plane,  whose horizontal and vertical asymptotes are   $d_2=\frac{mpcu^*v^*l^2}{n^2}$ and $d_1=-\frac{ru^*(2u^*-1-a)l^2}{n^2}$. We only need to consider the right branch of the  hyperbola,  which   intersects   with the $d_1-$axis at $d_1=\frac{{\rm det}J_{E^*}l^2}{mpcu^*v^*n^2}$. Thus, with respect to $n$, both the horizontal asymptote of the right branch of the  hyperbola and   the intersection with the $d_1-$axis are decreasing.


     Denote
     \begin{equation}\label{d1d2}
     d_2^T(n,d_1)\stackrel{\vartriangle}{=}\dfrac{mpcu^*v^*d_1\frac{n^2}{l^2}-{\rm det}J_{E^*}}{d_1\frac{n^4}{l^4}+ru^*(2u^*-1-a)\frac{n^2}{l^2}}.
     \end{equation}
      Noting that $u^*>\frac{a+1}{2}$ when $p\leq p_H$, we have $J_n<0$ if $d_2<d_2^T(n,d_1)$, and $J_n>0$ if $d_2>d_2^T(n,d_1)$.       

         \begin{theorem}\label{cen0den0}           If  $c<\frac{1}{r(1-a)}$, suppose  $1<p<p_H$ (If  $c>\frac{1}{r(1-a)}$, suppose  $p_{SN}<p<p_H$). \\
               The constant steady state  $E^*$    is locally asymptotically stable if $d_2>d_2^T(n,d_1)$ for all $n\in\mathbb{N}$, while it is Turing unstable if  $d_2<d_2^T(n,d_1)$ for some $n\in{\mathbb{{N}}}$ satisfying $n >l\sqrt{\frac{{\rm det}J_{E^*}}{mpcu^*v^*d_1}}$.  \end{theorem}

     On the $d_1-d_2$ plane,  we call  the boundary curve of stable region of steady state $E^*$  the Turing bifurcation curve $l^T$, which  is formed by a sequence of curve segments $l^T_n$ $(n=1,2,...)$
         \begin{equation}\label{th}
      d_2=d_2^T(n,d_1),~~ {\rm for}~ d_{1,n}<d_1\leq d_{1,n-1}, ~
          \end{equation}
          where $d_{1,n}$ is the intersection of $d_2=d_2^T(n+1,d_1)$ and $d_2=d_2^T(n,d_1)$, and $d_{1,0}$ can be infinite.

When $p=p_H$, $T_n<0$ for all $n\in\mathbb{N}$, similar as the analysis in the case of $p<p_H$, we can also get a Turing bifurcation curve, and we have the following conclusion.
 \begin{remark}
  Assume, in the $p-d_2$ plane, $d_2=d_2^T(n,d_1)$ and $p=p_H$ intersects at a point $TH$, then the characteristic equation has a pair of purely imaginary roots and a zero root at $TH$. According to the general theory of \cite{Guckenheimer}, a Turing-Hopf bifurcation may appear. However, rigorous derivation for the normal form is impossible due to the implicit form of $E^*$ in this paper. We shall give some numerical illustrations to show the dynamics of the system near the Turing-Hopf bifurcation point.
 \end{remark}


\subsection{Numerical simulations  for  diffusive system }

        In this section, we carry out some numerical simulations for  diffusive system (\ref{diffusion}).  We fix the parameters as   \begin{equation}\label{paradiffsimu}
                a=0.23,r=1.1, m=0.31,l=2.
                 \end{equation}

\begin{figure}
           \centering
    a) \includegraphics[width=0.45\textwidth]{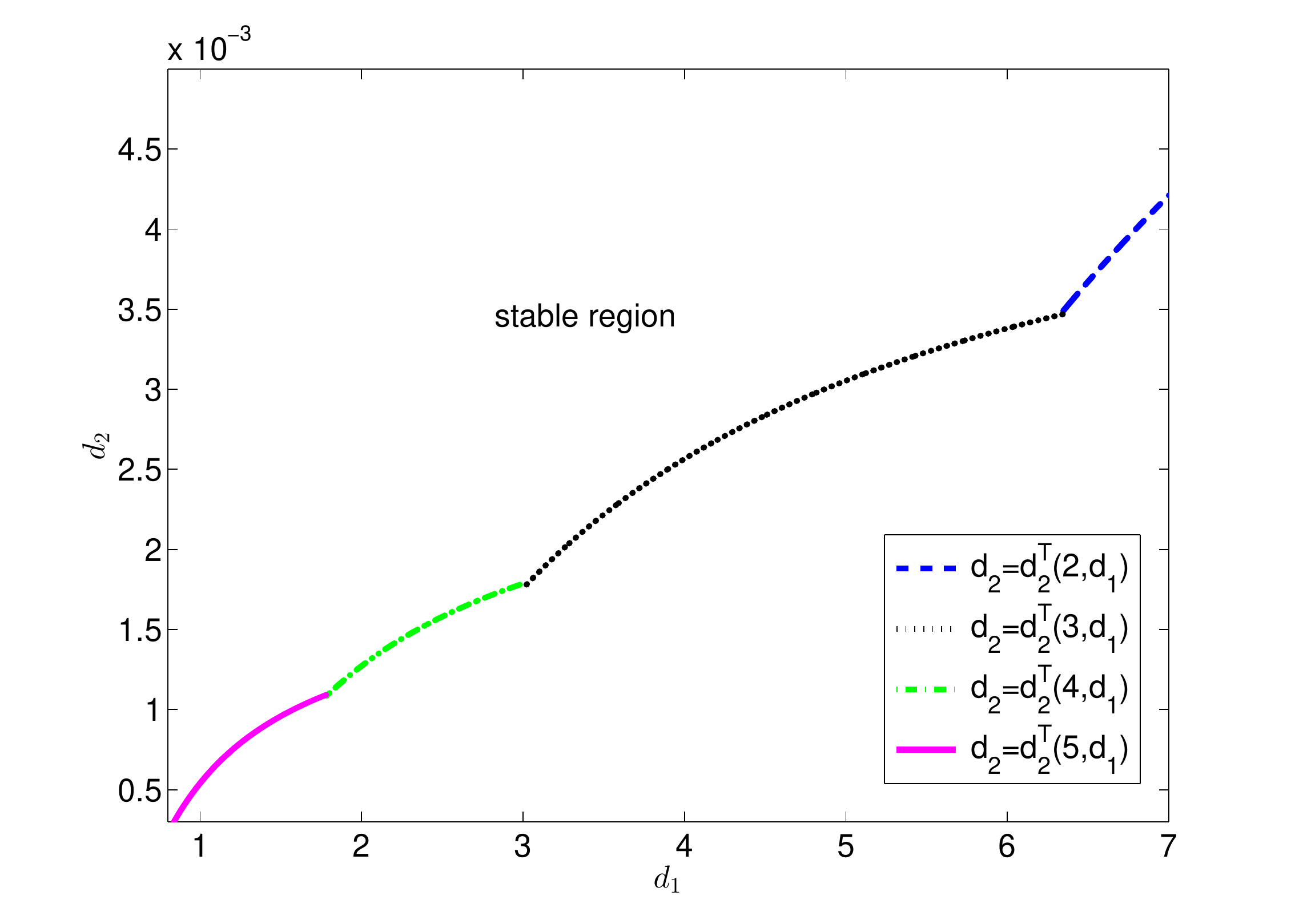}
      b) \includegraphics[width=0.45\textwidth]{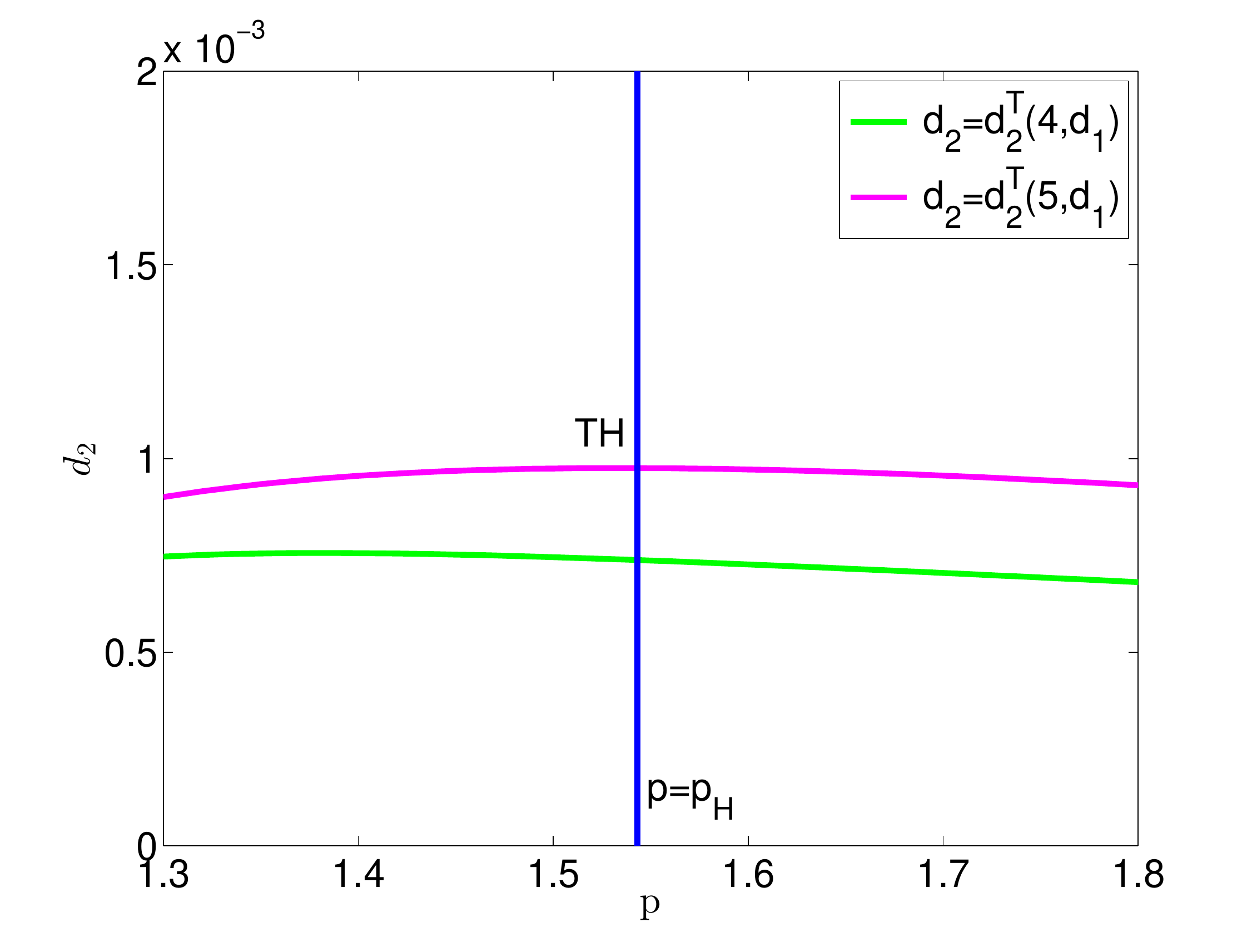}
         \caption{ a) Partial Turing bifurcation curve $l^T$ of system (\ref{diffusion})  with the parameters $  a=0.23,r=1.1,m=0.31,l=2,c=0.25<\frac{1}{r(1-a)},p=1.4<p_H $. b) When $a=0.23,r=1.1,m=0.31,l=2,c=0.25,d_1=1.496 $,  Turing bifurcation curve $d_2=d_2(5,d_1)$ intersects line $p=p_H$ at $TH~ (1.5432,0.0009)$, which is a Turing-Hopf bifurcation point.
         }
          \label{fig:turingc025d}
           \end{figure}

\begin{figure}
          \centering
   a) \includegraphics[width=0.45\textwidth]{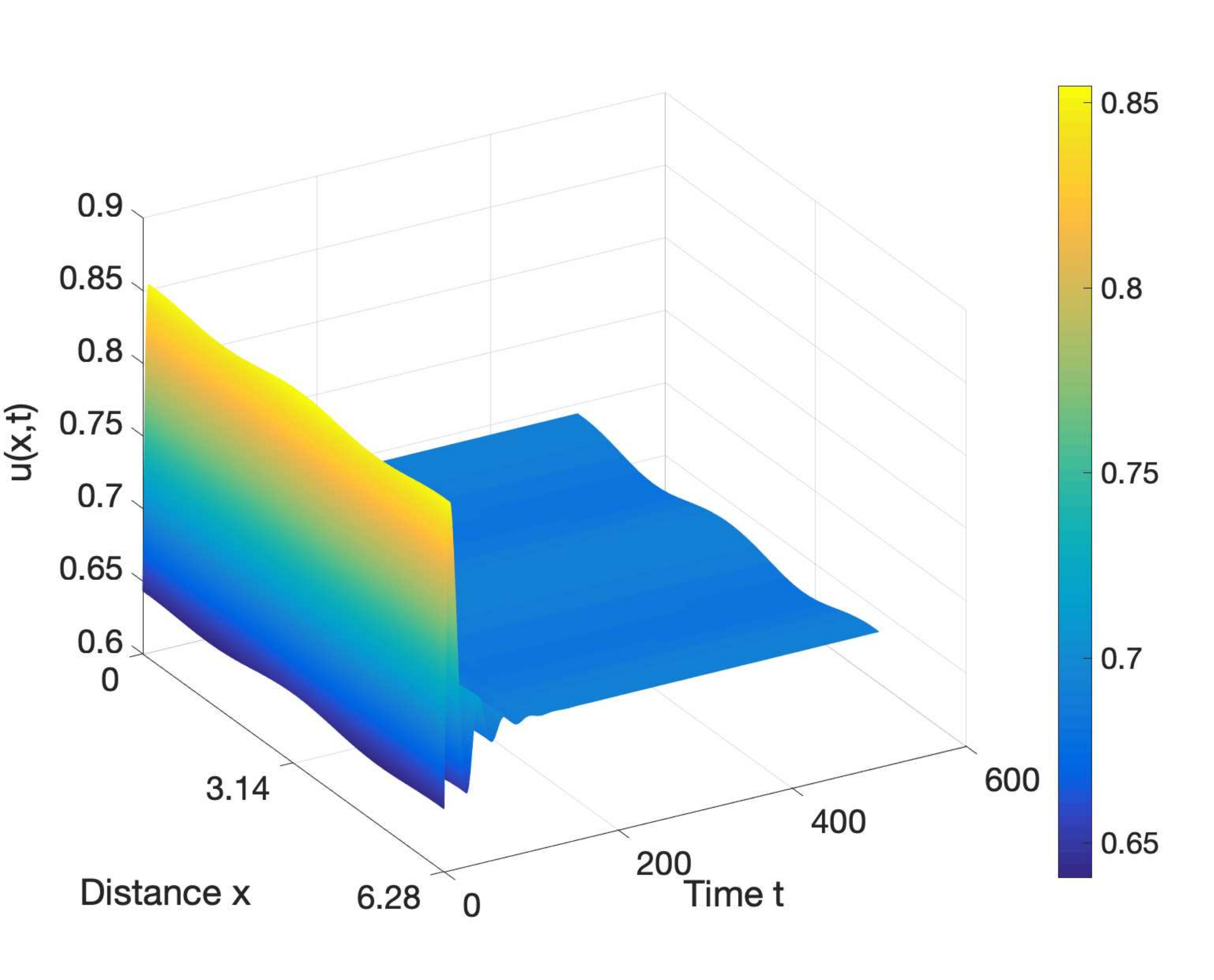}
       b)   \centering
   \includegraphics[width=0.45\textwidth]{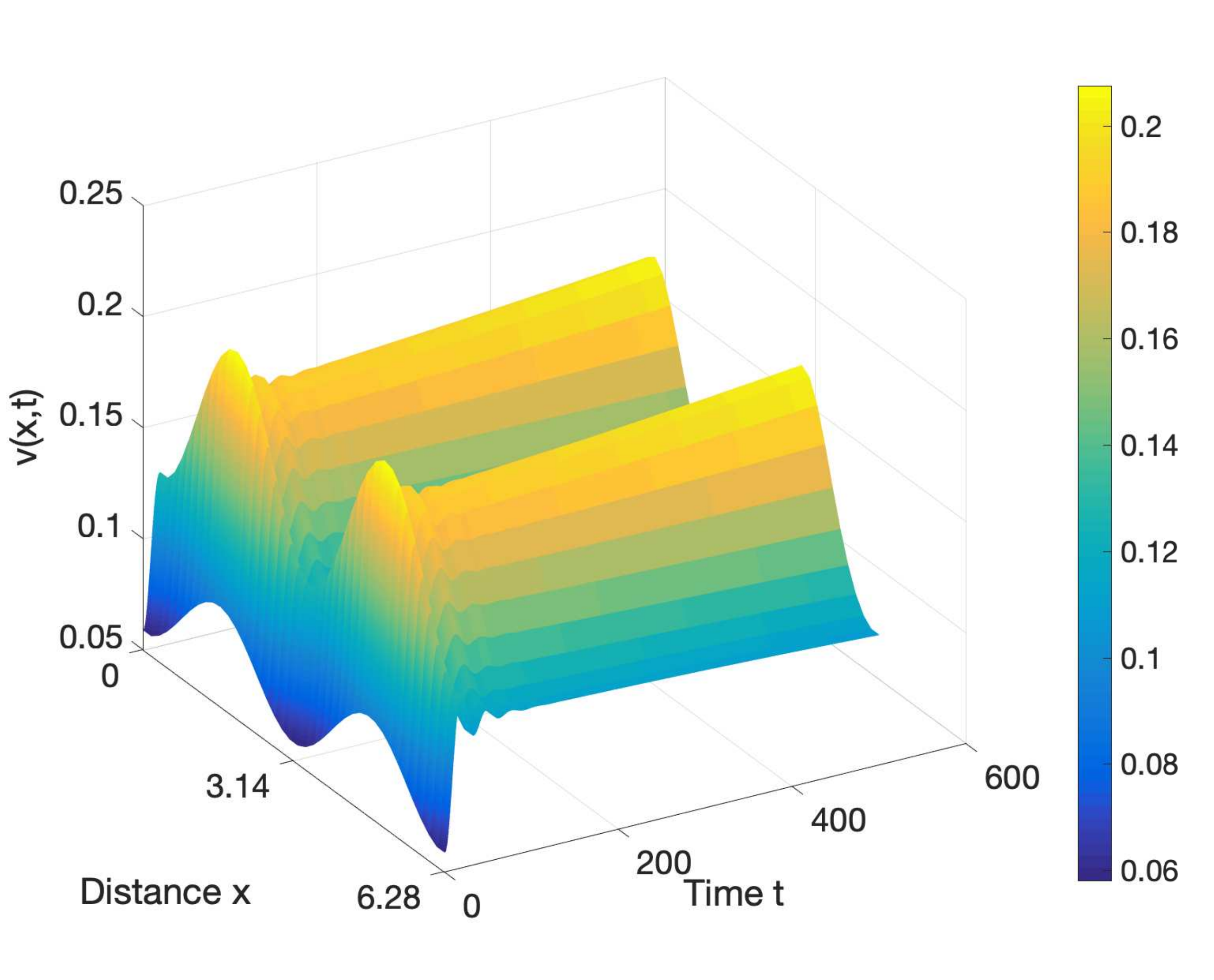}
  \caption{a) The Prey pattern  and  b) the predator pattern for system (\ref{diffusion}) with the parameters $  a=0.23,r=1.1,m=0.31,l=2,,c=0.25,p=1.4,d_1=1.496,d_2=0.000688 $.  The patterns indicate  instability induced by diffusion for the prey and predator.  Initial conditions are $u(x,0)=0.6+0.1\cos 2x$, $v(x,0)=0.08-0.02\cos 2x$.}
         \label{fig:turingc025}
          \end{figure}

   To illustrate the dynamics of (\ref{diffusion}) in the case of  weak cooperation, we choose $c=0.25$ such that $c<\frac{1}{r(1-a)}$. Let $p=1.4<p_H=1.5432$, from the  previous section, $E^*(0.6882,0.1516)$ is stable when $d_1=d_2=0$. On $d_1-d_2$ plane, we can draw curves determined by (\ref{d1d2}) for $n\in\mathbb{N}$,  and we only illustrate four curves for $n=2,3,4,5$. Correspondingly,  from  (\ref{th}), we can get four curve segments, forming partial Turing bifurcation curve  $l^T$ (see Fig. \ref{fig:turingc025d} a)).  From Theorem \ref{cen0den0},  if we choose $d_1$ and $d_2$ in the region above the curve,  the constant steady state $E^*$ is stable.  If we choose $d_1=1.496$ and $d_2=0.000688$ in the region under $l^T$,   the steady state $E^*$ of   the diffusive system is Turing unstable, and there is  a spatially inhomogeneous steady state   (see Fig. \ref{fig:turingc025}).

\begin{figure}
                 \centering
                  a) \includegraphics[width=0.45 \textwidth]{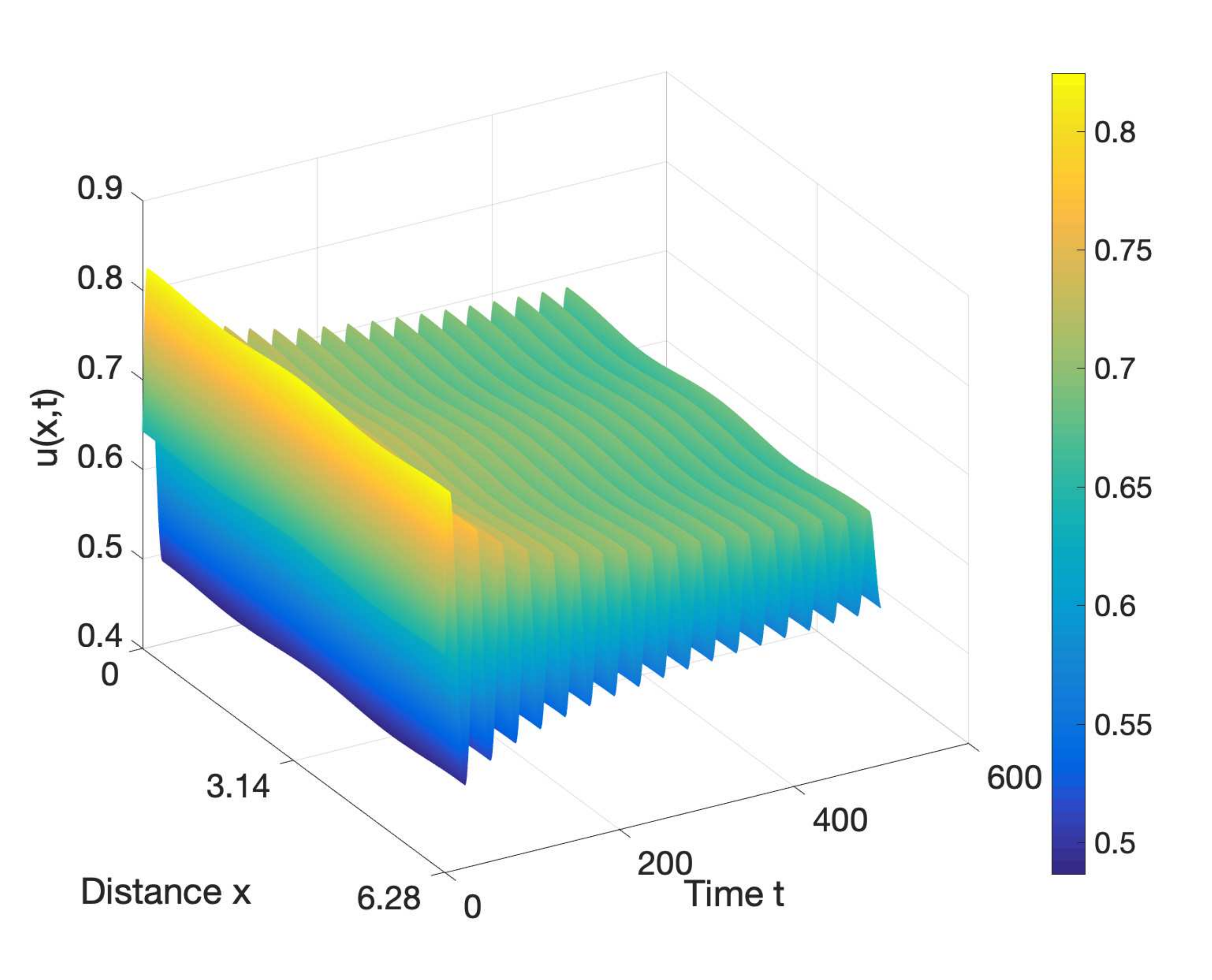}              b)   \centering
               \includegraphics[width=0.45 \textwidth]{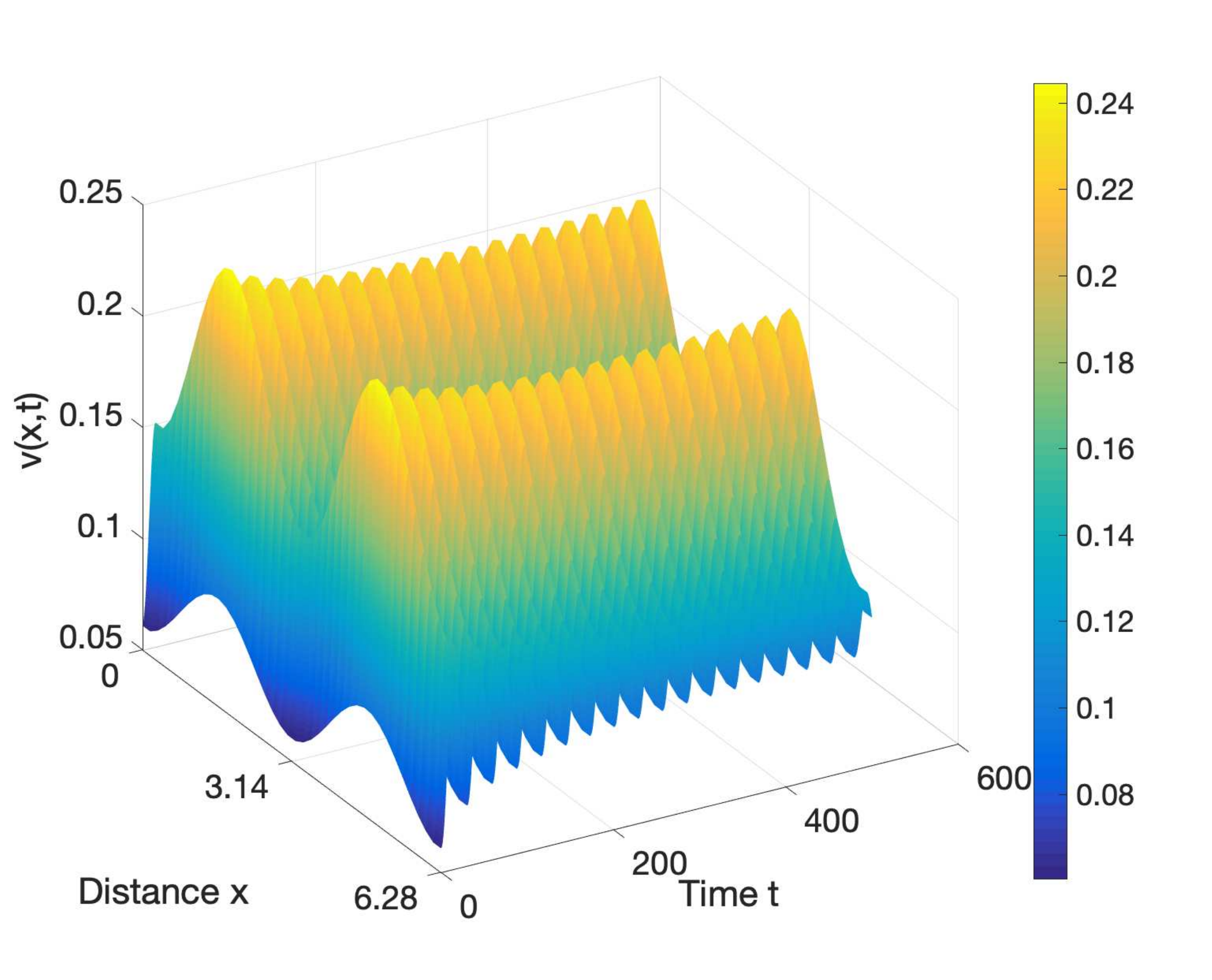}
              c) \includegraphics[width=0.45 \textwidth]{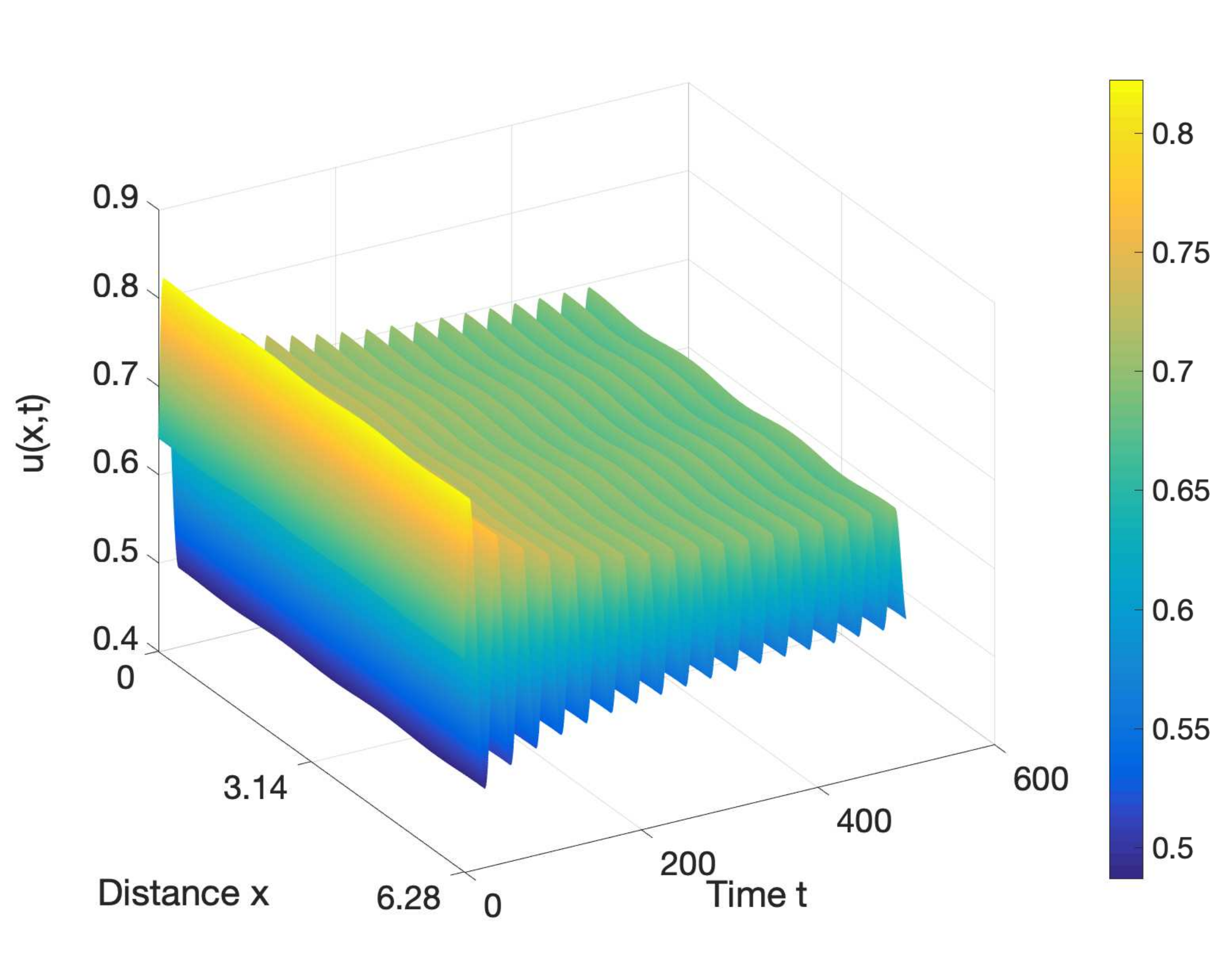}
             d)               \includegraphics[width=0.45 \textwidth]{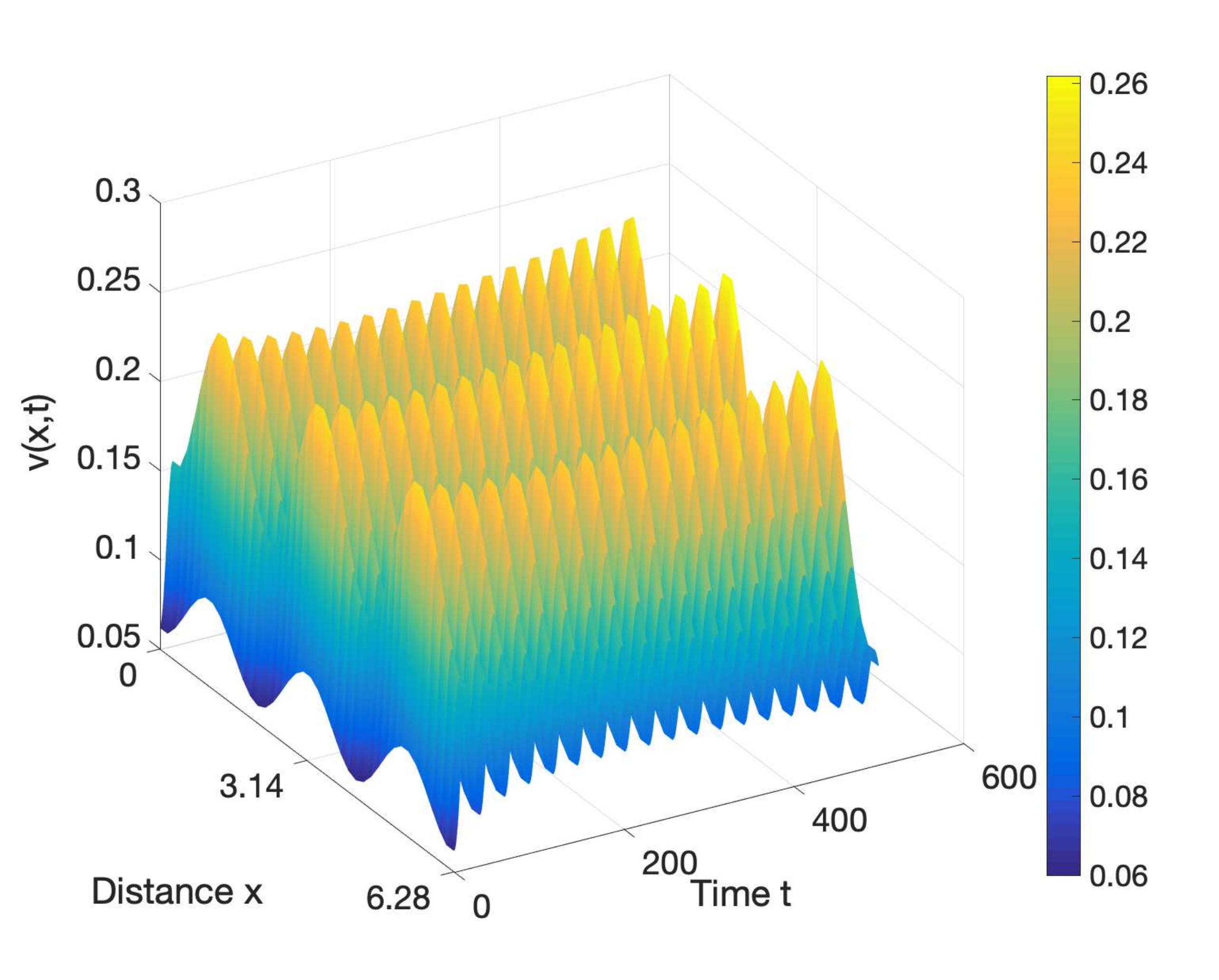}
        \caption{ Two spatially inhomogeneous periodic solutions coexist when we choose $ a=0.23, m=0.31,r=1.1$, $l=2$,  $c=0.25$, $p=p_H=1.5432$, $d_1=1.496$, and $d_2=0.00081$ in  system (\ref{diffusion}).   Initial conditions are $u(x,0)=0.6+0.1\cos 2x$, $v(x,0)=0.08-0.02\cos 2x$ for a) b), and $u(x,0)=0.6+0.1\cos 3x$, $v(x,0)=0.08-0.02\cos 3x$ for c) d). }
              \label{fig:turinghopfc025}
             \end{figure}

   If we fix $d_1=1.496$, from  (\ref{d1d2}), we can draw   Turing bifurcation curve  on $p-d_2$ plane  (see Fig. \ref{fig:turingc025d} b)),  and $d_2=d_2^T(5,d_1)$ intersects line $p=p_H$ at $TH~ (1.5432,0.0009)$.
     Choosing $p=1.5432,d_2=0.00081 $  near $TH$, two stable spatially inhomogeneous  periodic solutions  coexist when we choose two different initial values (see Fig. \ref{fig:turinghopfc025}).

\begin{figure}
           \centering
    a) \includegraphics[width=0.45\textwidth]{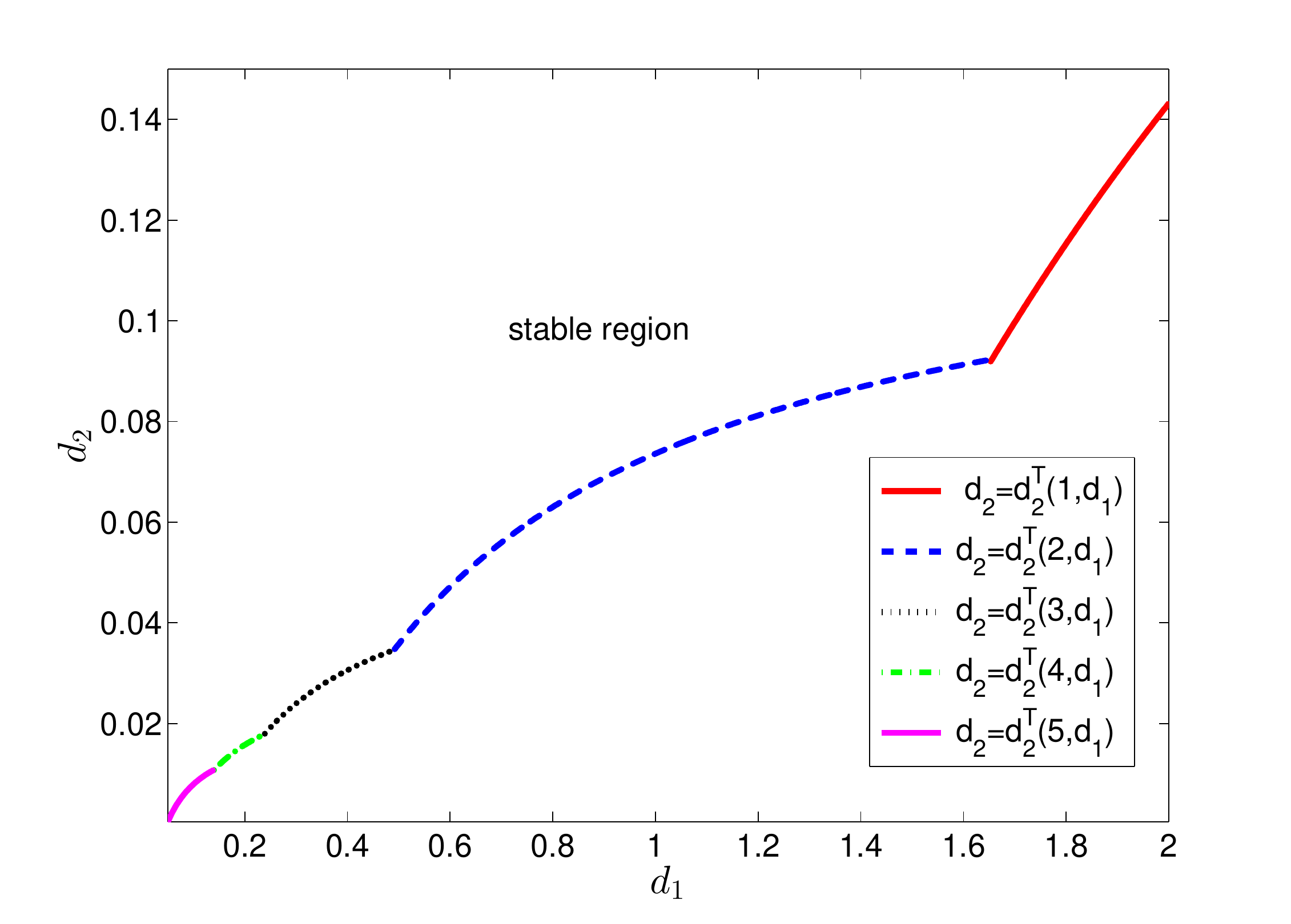}
      b) \includegraphics[width=0.45\textwidth]{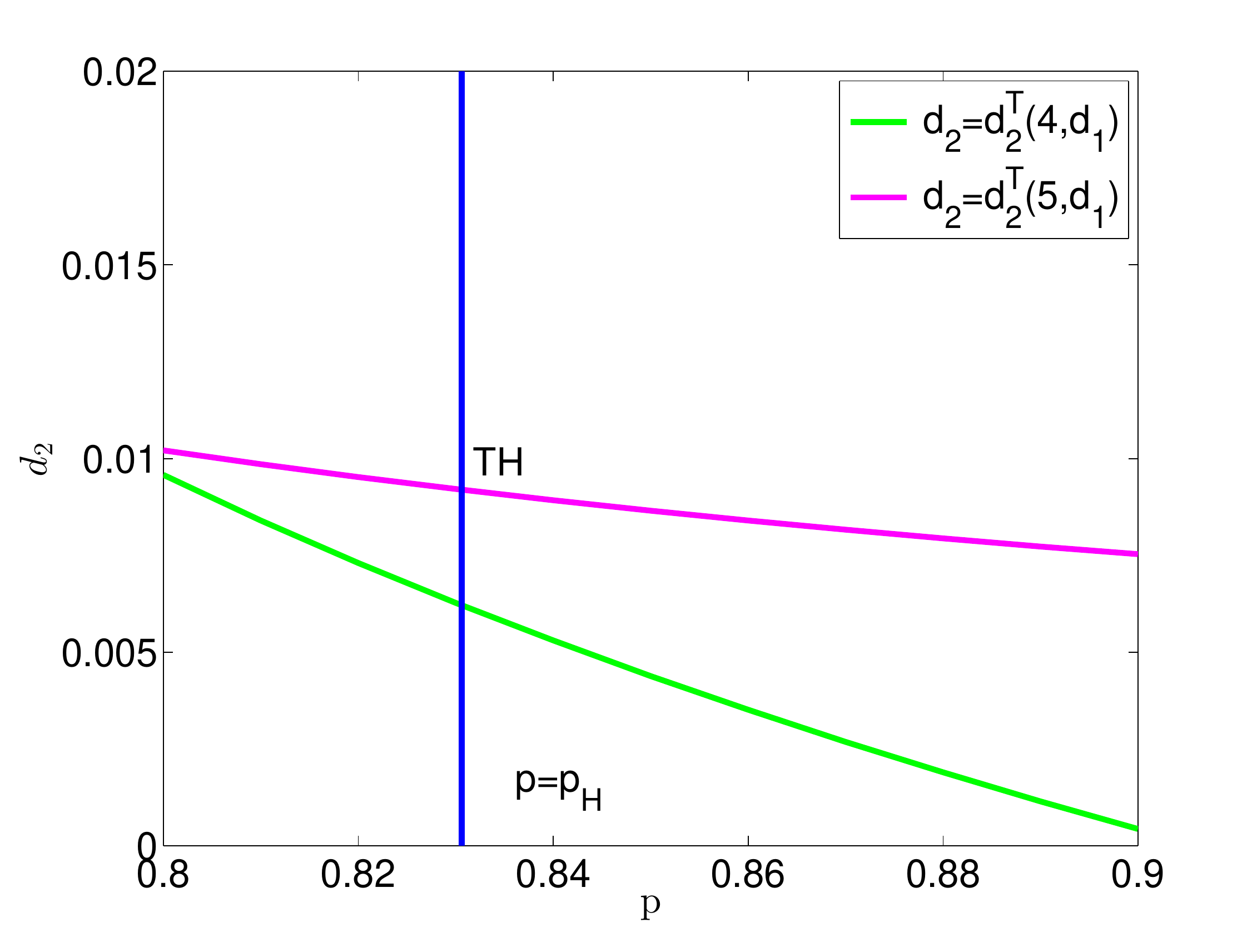}
         \caption{a) Partial Turing bifurcation curve $l^T$ of system (\ref{diffusion})  with the parameters $  a=0.23,r=1.1,m=0.31,l=2,c=8,p=0.8<p_H $. b)
         b) When $a=0.23,r=1.1,m=0.31,l=2,c=8,d_1=0.13 $,  Turing bifurcation curve $d_2=d_2(5,d_1)$ intersects line $p=p_H$ at $TH~ (0.8306,0.009)$, which is a Turing-Hopf bifurcation point.
         }
          \label{fig:turingc8d}
           \end{figure}

  To illustrate the dynamics of (\ref{diffusion}) in the case of  strong cooperation, we choose $c=8$ such that $c>\frac{1}{r(1-a)}$. Let $p=0.8<p_H=0.8306$. From (\ref{d1d2}) and (\ref{th}),  we can draw partial Turing bifurcation curve $l^T$ on $d_1-d_2$ plane (see Fig. \ref{fig:turingc8d} a)). The constant steady state $E^*( 0.7392,0.0864)$ is stable when we choose $d_1$ and $d_2$ in the region above $l^T$.
   If we choose $d_1=0.13$ and $d_2=0.006$ in the region below $l^T$,  from Theorem \ref{cen0den0}, the steady state $E^*$ of   the diffusive system is Turing unstable.  There are three spatially inhomogeneous steady states coexist  when we choose three different initial values.

\begin{figure}
                   \centering
                a) \includegraphics[width=0.45 \textwidth]{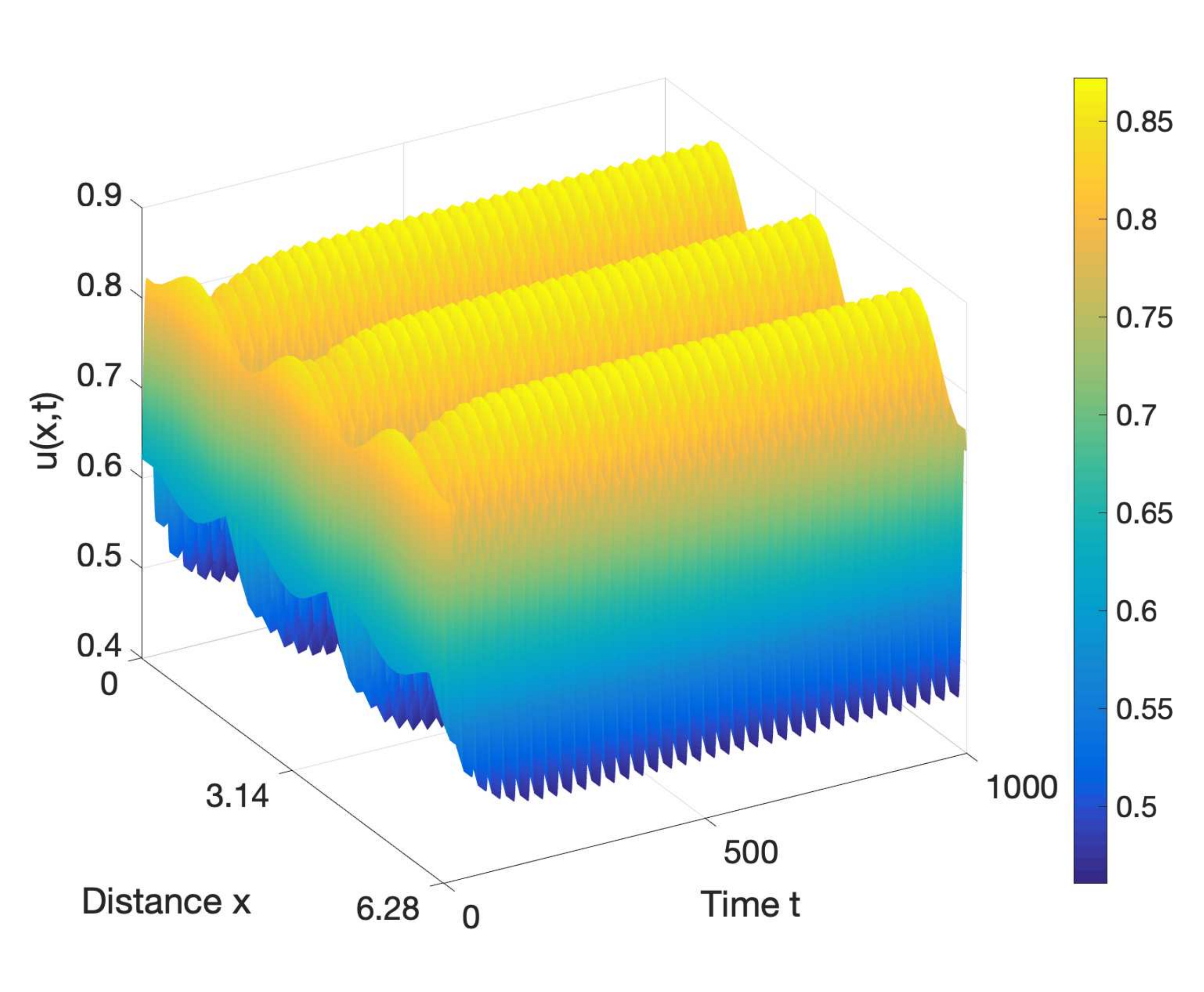}
                  b)   \centering
                \includegraphics[width=0.45 \textwidth]{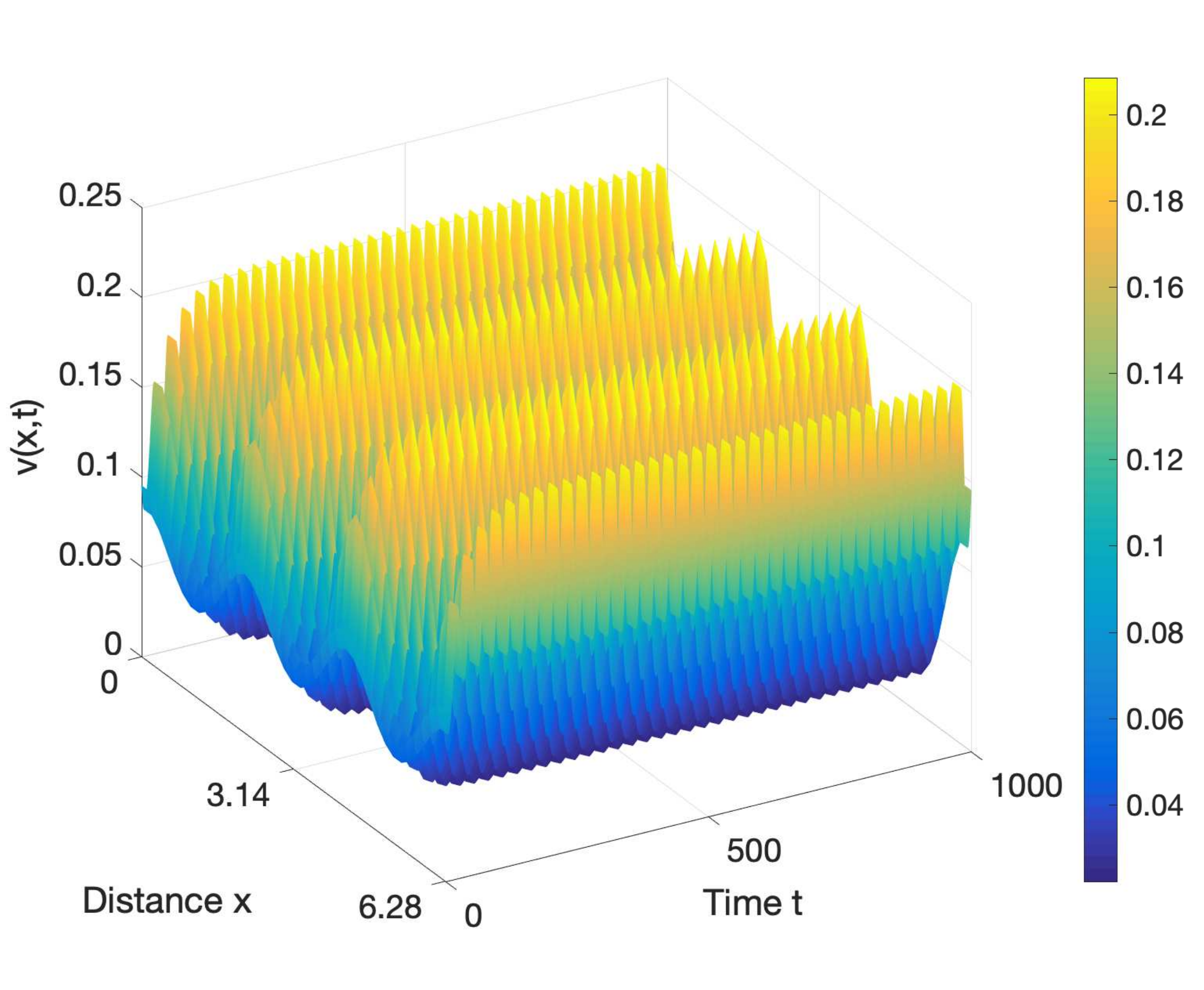}
                           c) \includegraphics[width=0.45 \textwidth]{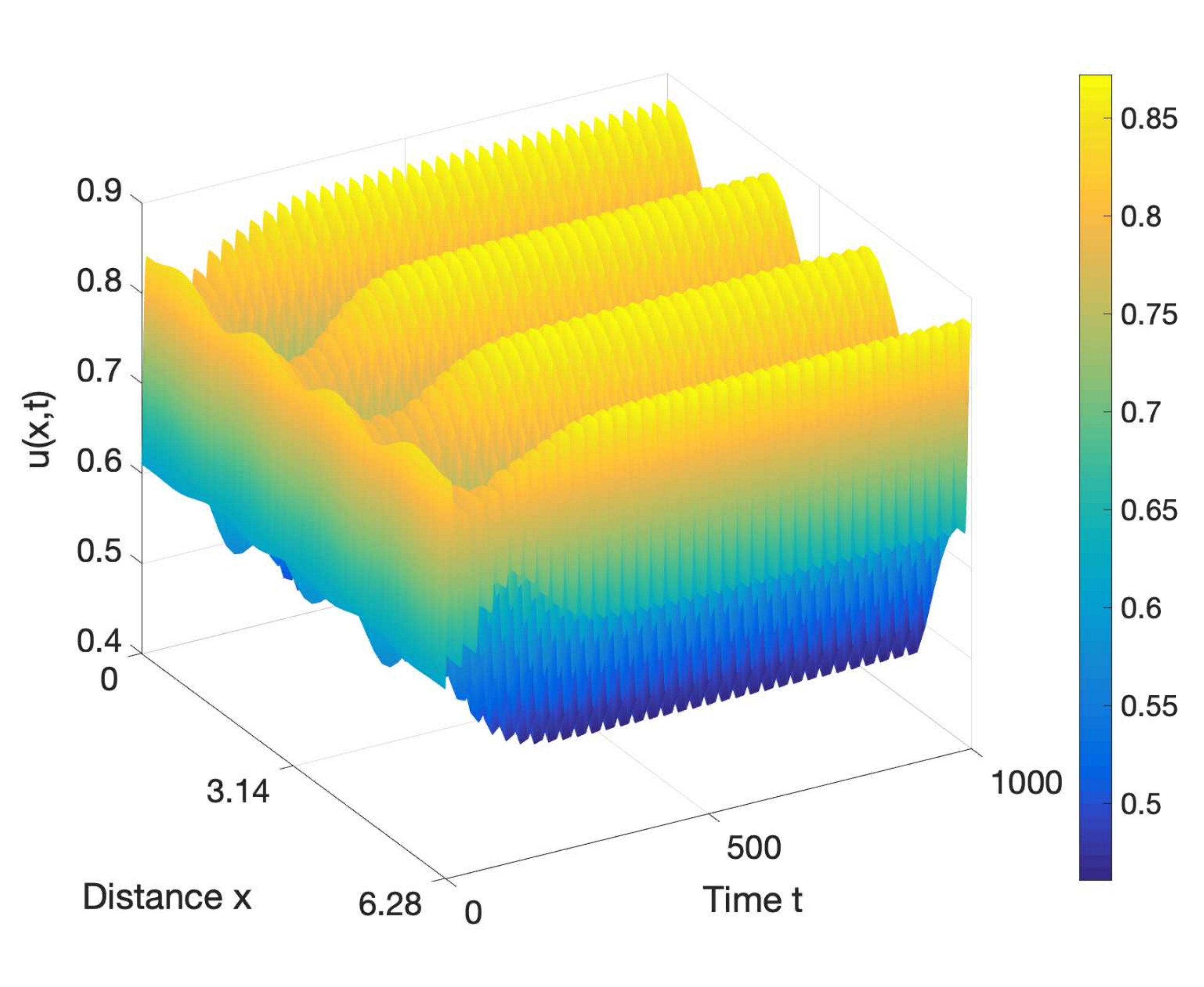}
                                d)   \centering
                              \includegraphics[width=0.45 \textwidth]{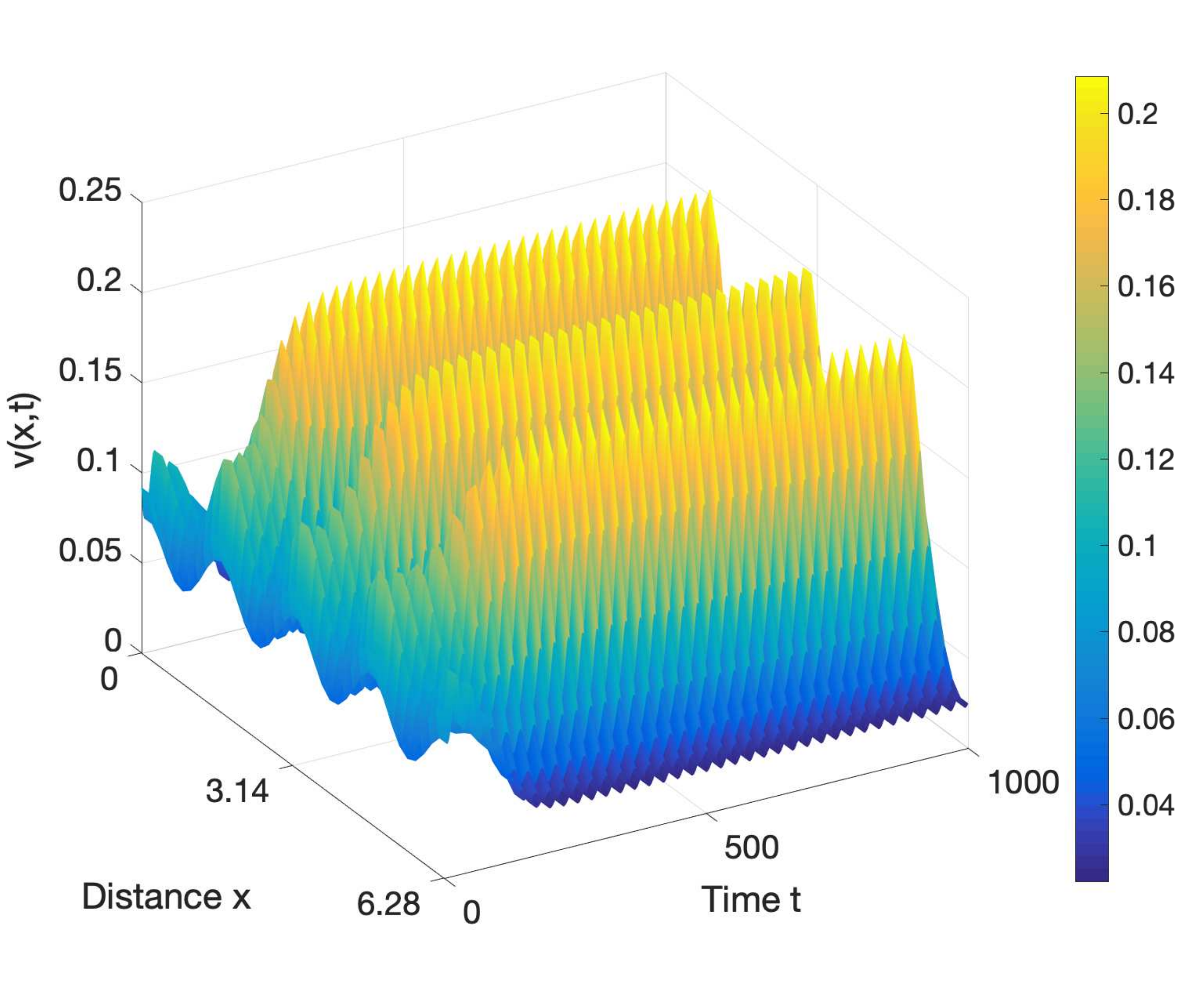}
                 \caption{ Two stable spatially inhomogeneous periodic solutions of system (\ref{diffusion})  coexist   with $ a=0.23, m=0.31,r=1.1$, $l=2$, $c=8$, $p=p_H=0.8306$, $d_1=0.13$, and $d_2=0.0065$.   Initial values are $u(x,0)=0.6+0.1\cos 4x$, $v(x,0)=0.08-0.02\cos 4x$ for a) b), and $u(x,0)=0.6+0.1\cos 3x$, $v(x,0)=0.08-0.02\cos 3x$ for c)  d). }
                \label{fig:turinghopfc8}
               \end{figure}

 If we fix $d_1=0.13$, from  (\ref{d1d2}), we can draw   Turing bifurcation curve  on $p-d_2$ plane  (see Fig. \ref{fig:turingc8d} b)),  and $d_2=d_2^T(5,d_1)$ intersects line $p=p_H$ at $TH~ (0.8306,0.009)$.
     Choosing $p=0.8306,d_2=0.0065 $  near $TH$,  we can illustrate two stable spatially inhomogeneous  periodic solutions  coexisting when we choose two different initial values (see Fig. \ref{fig:turinghopfc8}).

    \section{The dynamics of diffusive system with two delays}\label{sectiondelays}

  It has been widely  accepted that time delays   have very complex effect  on  the dynamics of a system, for example, some delays can destroy  the  stability of equilibria and  induce various oscillations  and periodic solutions.  There are time delays  in almost every  process of population interaction,  so it is more realistic to introduce time delays when we model the interaction of predator and prey.
    For example,   after the predator consuming the prey,   the reproduction of predator   is not  instantaneous  but   taking time for the  transition from prey  biomass into predator biomass.  We call this kind of time delay as a gestation delay.    There is an extensive literature about the studies of the dynamics of predator-prey
   models with the effect of time delay due to gestation of  the predator (see, for example, \cite{Kuang Ygestation,Gopalsamy Kgestation,Martin Agestation,Xia Jgestation} and references  cited therein).    Mature delay of the prey  has  also been   considered \cite {Wright E,RM May,XuCdelay,Yandelay}.
 However, to the best of our knowledge,  there is very little literatures on  delayed predator-prey system
  with Allee effect \cite{JankovicAl]eedelay}.

We consider the following diffusive system with two delays
\begin{equation}
\label{diffusion predatorless}
\left\{
\begin{array}{l}
\dfrac{\partial u(x,t)} {\partial t}= d_1\Delta u(x,t)+ru(x,t)[1-u(x,t-\tau_1)] (u(x,t)-a) \\~~~~~~~~~~~~~~~~ -(1+cv(x,t))u(x,t)v(x,t),~~~~~~~~~~~~~~~~~~~~~~~~~~~~~~~~~~~~~~~~~~~~~~~~~x\in (0,l\pi),\\
\dfrac{\partial v(x,t)}{\partial t }= d_2\Delta v(x,t)+mv(x,t)\left[  p u(x,t-\tau_2)\left( 1+cv(x,t-\tau_2)\right) -1 \right]  , ~x\in (0,l\pi),~ \\
\partial_v u(x,t)= 0,~~ \partial_v v(x,t) =0,   ~~~~~~~~~~~~~~~~~~~~~~~~~~~~~~~~~~~~~~~~~~~~~~~~~~~~~~~~~~~~~~~~x=0,l\pi.\\
\end{array}
\right.
\end{equation}
Here $\tau_1$ is   the time
  delay due to the maturation of the prey,
$\tau_2$  is  time delay due to gestation of  the predator, and we define $\tau={\rm max}\{\tau_1,\tau_2\}$.

Systems with multiple delays have attracted much attention \cite{K. L. Cooke,H. I. Freedman,J. Wei and S. Ruan,Xu J,twodelays,Huang C}.   Generally, delay may induce Hopf bifurcation, and if Hopf bifurcation curves     intersect,  double Hopf bifurcation may arise.  To figure out the effect of delay on the dynamics of systems,   Hopf bifurcation and double Hopf bifurcation induced by delay  have been investigated \cite{Yudouble,Dujia1,Ma double,Campell double}.
However,  in most of  literatures we mentioned above, the system  is  reduced into a system with one delay.
In fact, systems with multiple delays conform to reality better than one single delay.  Recently, the research on the dynamics and bifurcation analysis of system with  two simultaneously varying delays   are of great interest to scholars. In \cite{twodelays,Ge Jtwodelay}, the double Hopf bifurcation induced by two delays are studied by different methods. Through the analysis of double Hopf bifurcation, we can classify the topological structures of various  bifurcating solutions. By the classification, the dynamics in the neighbourhood of the double Hopf bifurcation point  in system can be obtained completely.

In the previous section, we have found the conditions of stability and Turing instability of constant steady state $E^*$ of  system (\ref{diffusion predatorless}) when $\tau_1=\tau_2=0$. We know  that  if  $ p<p_H$ and $d_2>d_2^T(n,d_1)$ for all $n\in\mathbb{N}$, the   steady state $E^*$ is locally asymptotically stable. In this section, we investigate the diffusive system (\ref{diffusion predatorless}) under the assumption  $ p<p_H$ and $d_2>d_2^T(n,d_1)$, and  focus on the effect of two delays on the dynamics of the diffusive system near $E^*$.

\subsection{Hopf and double Hopf bifurcation induced by two delays }\label{sectionhopddouble}
In this subsection, we investigate  the existence of Hopf bifurcation induced by two delays by the method of stability switching curves given in Ref. \cite{Lin X}, and give the condition of the existence of double Hopf bifurcation.

The linearization of  system (\ref{diffusion predatorless})  at   the  steady state $E^*(u^*,v^*)$ is given by
\begin{equation}\label{E2}
 \frac{\partial }{\partial t} \left( \begin{array}{l}
     u(x,t) \\
       v(x,t) \\
           \end{array}\right)
     =
    (D\Delta      +
     A) \left(  \begin{array}{l}
           u(x,t)\\
            v(x,t)\\
                     \end{array}
          \right)
    + B \left(  \begin{array}{l}
                u(x,t-\tau_1)\\
                 v(x,t-\tau_1)\\
                               \end{array}
               \right)
    +C\left(  \begin{array}{l}
            u(x,t-\tau_2)\\
              v(x,t-\tau_2)\\
              \end{array}
     \right),
\end{equation}
where
  \begin{equation}\label{ABC}
  \begin{array}{l}
     A=\left( \begin{array}{cc}
    ru^*(1-u^*)& -2cu^*v^*-u^*\\
    0& 0\\
       \end{array}\right),~
  B=\left( \begin{array}{cc}
      -ru^*(u^*-a)& 0\\
     0 & 0\\
         \end{array}\right),\\
    C=\left( \begin{array}{cc}
     0& 0\\
   mp(1+cv^*)v^* &mpcu^*v^*\\
          \end{array}\right),
   D={\rm diag} (d_1,d_2).
   \end{array}
  \end{equation}
 $u(x,t)$ and $v(x,t)$ satisfy the homogeneous Neumann boundary condition.

From Wu \cite{JWu}, the corresponding characteristic equation of Eq. (\ref{E2}) is
\begin{equation}\label{charaubar}
{\rm det}\left(\lambda I-M_n-A -Be^{-\lambda \tau_1}-Ce^{-\lambda \tau_2} \right) =0,
\end{equation}
where $I$ is a $2\times2$ identity matix, $M_n= -n^2/l^2{\rm diag}(d_1, d_2)$, $n\in \mathbb{N}_0=\{0,1,2,...\}$.
Eq. (\ref{charaubar}) can be written in the following form:
\begin{equation}
\label{character}
D_n(\lambda;\tau_1,\tau_2)= P_{0,n}(\lambda)+P_{1,n}(\lambda)e^{-\lambda\tau_1}+P_{2,n}(\lambda)e^{-\lambda\tau_2}+P_{3,n}(\lambda)e^{-\lambda(\tau_1+\tau_2)}=0,
\end{equation}
where
\begin{equation}\label{p1p2p3}
\begin{array}{l}
P_{0,n}(\lambda)=(\lambda+d_1\frac{n^2}{l^2}-a_{11})(\lambda+d_2\frac{n^2}{l^2}),
\\P_{1,n}(\lambda)=-b_{11}(\lambda+d_2\frac{n^2}{l^2}),
\\P_{2,n}(\lambda)=-c_{22}(\lambda+d_1\frac{n^2}{l^2}-a_{11})-a_{12}c_{21},
\\P_{3,n}(\lambda)=b_{11}c_{22}.
\end{array}
\end{equation}

The characteristic equation with the form of Eq. (\ref{character}) has been investigated by  Lin and Wang \cite{Lin X}.  They derived an explicit expression for the stability switching curves, on which there is a pair of purely imaginary roots for Eq. (\ref{character}).
Moreover, they gave a criterion to determine the crossing directions, i.e., on which side of the stability switching curve there are two more characteristic roots with positive real parts.  Using this method, we can find all the stability switching curves  in the $(\tau_1, \tau_2)$ plane, and determine the crossing directions. We leave the details in the   \ref{sectionStability switching curves} and  \ref{sectioncrossing}.

Moreover, we have the following  Hopf bifurcation theorem with two parameters.
          \begin{theorem}\label{bifur}
For each $j\in\{1,2,\cdots,N\}$, $\mathcal{T}_n^j$, defined by (\ref{Tk}),  is a Hopf bifurcation curve in the following sense:
for any $p\in \mathcal{T}_n^j$ and for any smooth curve $\Gamma$ intersecting with $\mathcal{T}_n^j$ transversely at $p$, we define the tangent of $\Gamma$ at $p$ by $\overrightarrow{l}$. If $
\frac{\partial {\rm Re}\lambda}{\partial \overrightarrow{l}}\mid_p\neq 0
$, and the other eigenvalues of (\ref{character}) at $p$ have non-zero real parts, then system (\ref{diffusion predatorless}) undergoes a Hopf bifurcation at $p$ when parameters $(\tau_1,\tau_2)$ cross $\mathcal{T}_n^j$ at $p$ along $\Gamma$.
\end{theorem}
The theorem can be proved  by a similar  method in \cite{twodelays}.

\begin{remark}
	Suppose  that $\mathcal{T}_{k_1}^{j_1}$ and  $\mathcal{T}_{k_2}^{j_2}$ intersect at a point $(\tau_1,\tau_2)$, with the corresponding values of $\omega$ being   $\omega_{j_1,k_1}\in \Omega_{j_1,k_1}$ and $\omega_{j_2,k_2}\in\Omega_{j_2,k_2}$.  Then there are two pairs of purely imaginary roots of (\ref{character}) at the intersection, which are $\pm i  \omega_{j_1,k_1}$ and $ \pm i \omega_{j_2,k_2}$, denoted by $\pm i  \omega_1$ and $\pm i  \omega_2$ for convenience. Thus, system (\ref{diffusion predatorless}) may undergo double Hopf bifurcations at the intersection of  two stability switching curves near the constant steady state   $E^*$.
\end{remark}

\subsection{Normal form on the center manifold  for  double Hopf bifurcation}
\label{normal form}
In this subsection, applying  the normal form method of partial functional differential equations \cite{Faria}, we derive   normal form  of double Hopf bifurcation taking two delays as bifurcation parameters.  Then, we can classify the topological structures of various  bifurcating solutions, and get  the dynamics in the neighbourhood of the double Hopf bifurcation point  in system (\ref{diffusion predatorless}).

For the Neumann boundary condition, we define the real-valued Sobolev space
$$
  	X=\{(u,v)^T\in H^2(0,l\pi)\times H^2(0,l\pi) \arrowvert   \dfrac{\partial u}{\partial x} =\dfrac{\partial v}{\partial x}= 0, x=0,l\pi\},
$$
  and the abstract space
  $\mathcal{C} = C([-1, 0], X)$.
 Define the complexification space of the real-valued Hilbert space  $X$ by
$$
 X_{\mathbb{C}}:=X\oplus iX=\{U_1+iU_2:U_1,U_2\in X\},
$$
   with  the general complex-value   $L^2$ inner product
$$
\langle U,V\rangle=\int_0^{l\pi}(\overline{u}_1v_1+\overline{u}_2v_2)dx,
~~\mbox{for}~~U=(u_1,u_2)^T, V=(v_1,v_2)^T\in X_{\mathbb{C}}.
$$
Let $\mathscr{C}:=C([-1,0],X_{\mathbb{C}})$ denotes the phase space with the sup norm, and write $u^t\in\mathscr{C}$ for $u^t(\theta)=u(t+\theta)$, $-1\leq \theta\leq 0$.

Without loss of generality, we  assume  $\tau_1<\tau_2$ in this section. Denote the double Hopf bifurcation point by $(\tau_1^*,\tau_2^*)$.   Let $\overline{u}(x,t)=u(x,\tau_2t)-u^*,\overline{v}(x,t)=v(x,\tau_2t)-v^*$,  and set  $\sigma_1=\tau_1-\tau_1^*$ and  $\sigma_2=\tau_2-\tau_2^*$  as two bifurcation parameters.  Drop the bars, and denote $U(t)=(u(t),v(t))^T$,  then system (\ref{diffusion predatorless})  can be written as
 \begin{equation}
     \label{dudt}
   \dfrac{dU(t)}{dt}=D(\tau_1^*+\sigma_1,\tau_2^*+\sigma_2)\Delta U(t)+L(\tau_1^*+\sigma_1,\tau_2^*+\sigma_2)(U^t)+F(\tau_1^*+\sigma_1,\tau_2^*+\sigma_2,U^t),
          \end{equation}
where
 \begin{equation*}
 \begin{array}{l}
     D(\tau_1^*+\sigma_1,\tau_2^*+\sigma_2)=(\tau_2^*+\sigma_2)D=\tau_2^*D+\sigma_2 D,\\
    L(\tau_1^*+\sigma_1,\tau_2^*+\sigma_2)(\phi)
    =(\tau_2^*+\sigma_2)\left[A\phi(0)+B\phi(-\frac{\tau_1^*+\sigma_1}{\tau_2^*+\sigma_2})+C\phi(-1))\right],\\
    F(\tau_1^*+\sigma_1,\tau_2^*+\sigma_2,\phi)=   (\tau_2^*+\sigma_2) (F_1,
   F_2)^T,
      \end{array}
         \end{equation*}
with $A,~B$, $C$ and $D$ being defined in (\ref{ABC}), and for $\varphi\in \mathcal{C}$
\begin{equation*}
      \begin{array}{l}
F_1=r(1-u^*)\phi_1^2(0)-(1+2cv^*)\phi_1(0)\phi_2(0)+(ra-2ru^*)\phi_1(0)\phi_1(-\frac{\tau_1^*+\sigma_1}{\tau_2^*+\sigma_2})\\
~~~~~~~~~~-cu^*\phi_2^2(0) -r\phi_1^2(0)\phi_1(-\frac{\tau_1^*+\sigma_1}{\tau_2^*+\sigma_2})-c\phi_1(0)\phi_2^2(0),
\end{array}
\end{equation*}
and
\begin{equation*}
 \begin{array}{l}
F_2=mp(1+cv^*)\phi_2(0)\phi_1(-1) +mpc[u^*\phi_2(0) +v^*\phi_1(-1)+\phi_2(0)\phi_1(-1)]\phi_2(-1).
\end{array}
\end{equation*}

Consider the linearized system of (\ref{dudt})
   \begin{equation}
     \label{dudtlinear}
   \dfrac{dU(t)}{dt}=\tau_2^*D\Delta U(t)+\tau_2^*(AU^t(0)+BU^t(-\tau_1^*/\tau_2^*)+CU^t(-1))\stackrel{\vartriangle}{=}D_0\Delta U(t)+L_0(U^t).
          \end{equation}
We know that the  normalized eigenfunctions  of $D\Delta$ on $X$ corresponding to the eigenvalues $-d_1\frac{n^2}{l^2}$ and $-d_2\frac{n^2}{l^2}$( $n\in \mathbb{N}_0=\{0,1,2,\cdots\}$) are
$$
\beta_n^{1}(x)=\gamma_n(x)(1,0)^T~ \mbox{and}~~\beta_n^{2}(x)=\gamma_n(x)(0,1)^T,
$$
respectively, where
 $\gamma_n(x)=\dfrac{\cos\frac{n}{l}x}{\parallel\cos\frac{n}{l}x\parallel_{L^2}}$.  Define
$\mathscr{B}_n:={\rm span} \left\lbrace \langle v(\cdot),\beta_n^{j}\rangle\beta_n^{j}~\arrowvert  ~v\in \mathscr{C},j=1,2\right\rbrace$, satisfying $L(\mathscr{B}_n)(\tau_1,\tau_2)\subset {\rm span}\{\beta_n^{1},\beta_n^{2}\}$.
Denote
 $\left\langle v(\cdot),\beta_n \right\rangle=\left(
\langle v(\cdot),\beta_n^{1} \rangle,
\langle v(\cdot),\beta_n^{2} \rangle
\right)^T. $

Write system (\ref{dudt})  as
          \begin{equation}
               \label{dudt2}
             \dfrac{dU(t)}{dt}=D_0\Delta U(t)+L_0(U^t)+\widetilde{F} (\sigma,U^t),
                    \end{equation}
where $$\begin{aligned}
\widetilde{F} (\sigma,U^t)&=\sigma_2(D\Delta U^t(0)+AU^t(0)+BU^t(-\tau_1^*/\tau_2^*)+CU^t(-1))\\&+(\tau_2^*+\sigma_2)B(U^t(-(\tau_1^*+\sigma_1)/(\tau_2^*+\sigma_2))-U^t(-\tau_1^*/\tau_2^*))+F(\tau_1^*+\sigma_1,\tau_2^*+\sigma_2,U^t).
\end{aligned}$$
System  (\ref{dudt2}) can be rewritten as an abstract ordinary differential equation on $\mathscr{C}$  \cite{Faria}
\begin{equation}\label{ode}
\frac{d}{dt}U^t=\mathcal{A} U^t+X_0\widetilde{F} (\sigma,U^t),
\end{equation}
where $\mathcal{A}$ is the infinitesimal generator of the $C_0$-semigroup of solution maps of  the linear equation (\ref{dudt}), defined by
\begin{equation}
  \label{A}
\mathcal{A}\varphi=\dot{\varphi}+X_0[D_0\Delta\varphi(0)+L_0(\varphi)-\dot{\varphi}(0)],
 \end{equation}
  with ${\rm dom}(\mathcal{A})=\{\varphi\in\mathscr{C}:\dot{\varphi}\in\mathscr{C},\varphi(0)\in {\rm dom}(\Delta)\}$,  and $X_0$ is given by
  $X_0(\theta)=0$ for  $\theta\in[-1,0)$ and $X_0(0)=I$. Clearly,  $\mathcal{A}:\mathscr{C}_0^1\cap\mathscr{C}\rightarrow \mathscr{C}.$

   Then on $\mathscr{B}_n$, the linear equation
  $$
  \frac{d}{dt}U(t)=D_0\Delta U(t)+L_0(U^t)
  $$
 is equivalent to the retarded functional differential equation on $\mathbb{C}^2$:
   \begin{equation}
   \label{RFDE}
   \dot{z}(t)=-\frac{n^2}{l^2}D_0z(t)+L_0z^t.
  \end{equation}
  By the Riesz representation theorem, there exists a matrix whose components are bounded variation functions $\eta_k\in BV([-1,0],\mathbb{R}^{2\times 2})$ such that
   \begin{equation*}
   -\frac{n_k^2}{l^2}D_0\varphi(0)+L_0(\varphi)=\int_{-1}^0d\eta_k(\theta)\varphi(\theta), \varphi\in \mathscr{C}.
   \end{equation*}
   Let $A_k$ ($k=1,2$) denote the infinitesimal generator of the semigroup generated by (\ref{RFDE}),  and $A_k^*$ denote the formal adjoint of $A_k$ under the bilinear form
   \begin{equation*}
   (\alpha,\beta)_k=\alpha(0)\beta(0)-\int_{-1}^0\int_0^\theta\alpha(\xi-\theta)d\eta_k(\theta)\beta(\xi)d\xi.
   \end{equation*}

The calculations of normal form are very long, so we leave them in  supplementary materials.
 Based on the derivation in   supplementary materials,  the normal form truncated to the third order on the center manifold for double Hopf bifurcation is obtained.
Making the  polar  coordinate transformation,
then we obtain the following system corresponding to    the truncated  normal form
 \begin{equation}\label{normalformcylin}
 \begin{aligned}
 &\dot{\rho}_1=r_1(\nu_1+r_1^2+br_2^2),\\
 &\dot{\rho}_2=r_2(\nu_2+cr_1^2+dr_2^2),
 \end{aligned}
 \end{equation}
where
$$
\nu_1=\epsilon_1({\rm Re}B_{11}\sigma_1+{\rm Re}B_{21}\sigma_2),
 \nu_2=\epsilon_1({\rm Re}B_{13}\sigma_1+{\rm Re}B_{23}\sigma_2),   b=\frac{\epsilon_1\epsilon_2{\rm Re}B_{1011}}{{\rm Re}B_{0021}}, c=\frac{{\rm Re}B_{1110}}{{\rm Re}B_{2100}},d=\epsilon_1\epsilon_2,
$$
and $B_{11}, B_{21},B_{13},B_{23},B_{1011},B_{0021},B_{1110}$ and $B_{2100}$ are the coefficients of the normal form obtained in supplementary materials.   From chapter 7.5 in Ref. \cite{Guckenheimer},  there are twelve distinct kinds of unfoldings for Eq. (\ref{normalformcylin}).

\subsection{Numerical simulations for diffusive system with two delays}\label{section numerical delays}

 In this section, we carry out some numerical simulations for diffusive system (\ref{diffusion predatorless}). Fix
               \begin{equation*}
               \label{paraode}
               a=0.23,r=1.1,m=0.31, l=2,
                 \end{equation*}
          which is the same as in (\ref{paradiffsimu}).
    Fix $c=0.25$ such that $c<\frac{1}{r(1-a)}$,    and choose $p=1.2<p_H$, $d_1=0.3$ and $d_2=0.4$, such that $d_2>d_2^T(1,d_1)$ for all $n\in\mathbb{N}$ .  From Theorem \ref{cen0den0}, one  can get the unique constant  steady state $E^*( 0.80945,0.11797)$  is locally asymptotically stable.

\begin{figure}
            a)   \centering
       \includegraphics[width=0.44\textwidth,height=0.32\textwidth]{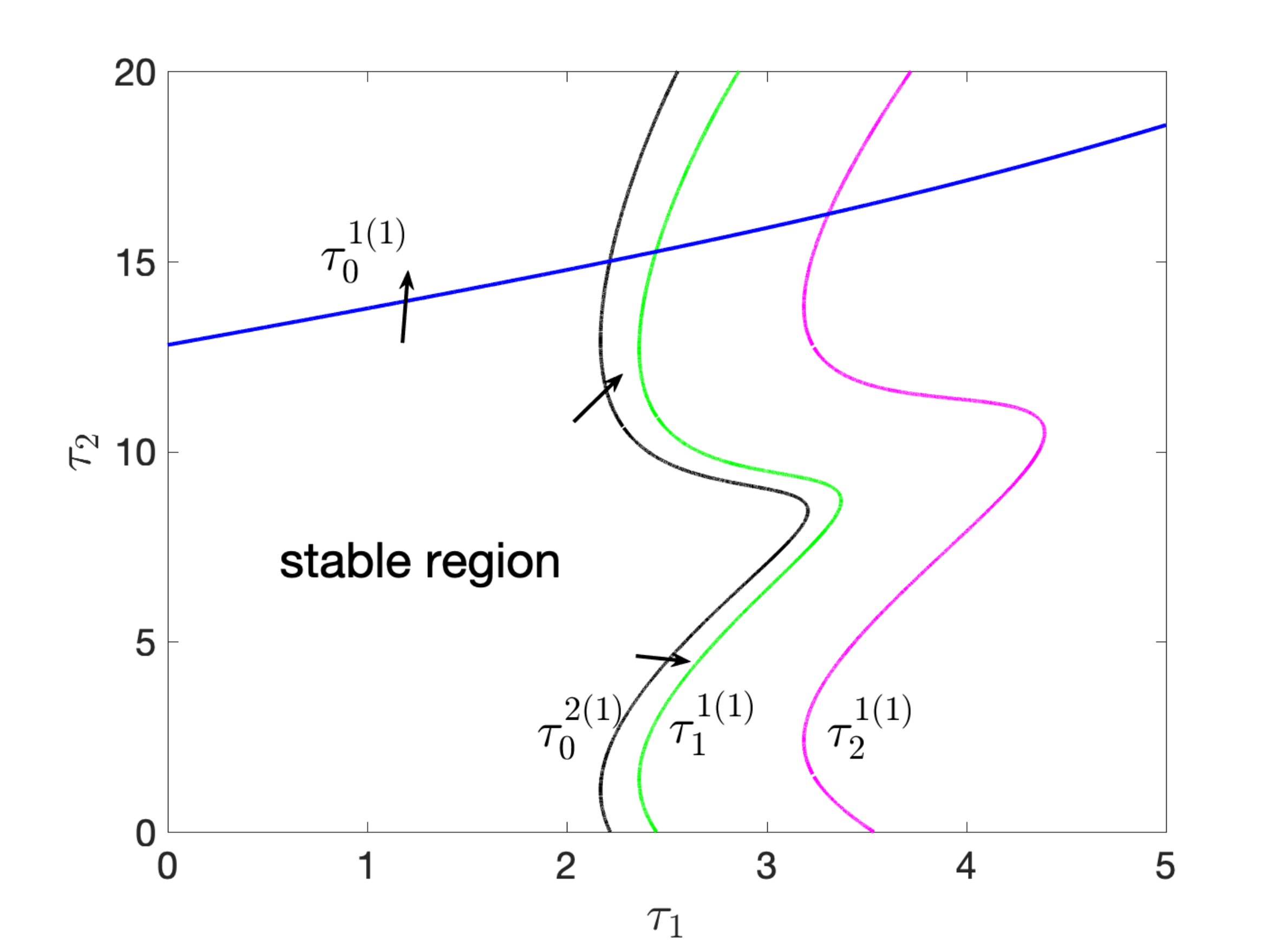}
       b) \includegraphics[width=0.44\textwidth,height=0.32\textwidth]{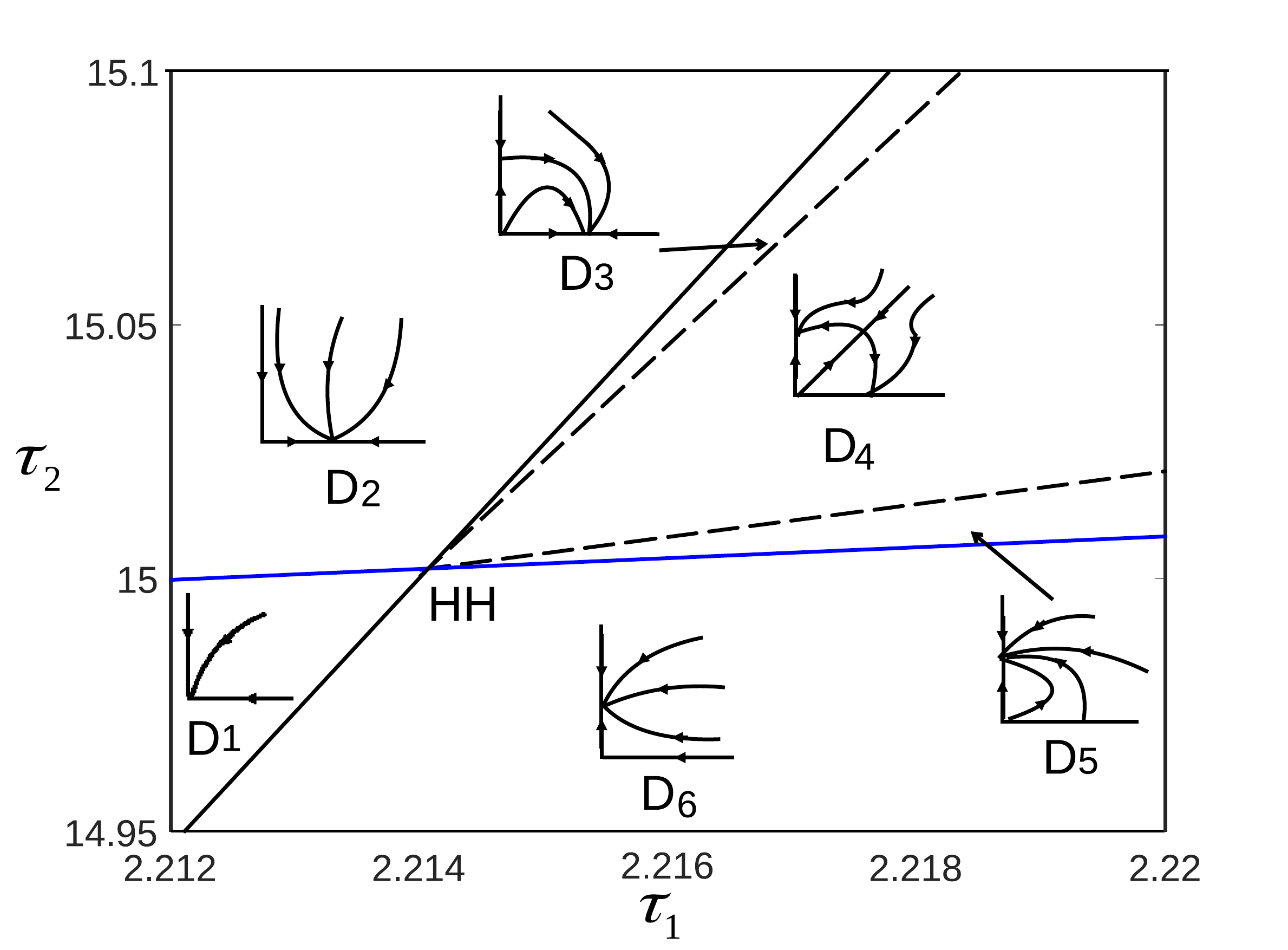}
   \caption  {a) The left-most curve and the lowest curve of $\mathcal{T}_0$ intersect at  $(\tau_1,\tau_2)=(2.2140,15.0019)$,
   	which is a double Hopf bifurcation point on the $\tau_1-\tau_2$ plane.  Crossing directions are marked by arrows.  b) The complete bifurcation sets near double Hopf bifurcation point HH.}
    \label{fig:tau1tau2}
   \end{figure}

Now we  illustrate the effect of two delays on the dynamics (\ref{diffusion predatorless}).  Following the steps in   \ref{sectionStability switching curves} and \ref{sectioncrossing}, we can draw all the stability curves on $(\tau_1,\tau_2)$ plane, and decide the crossing direction, which are shown in Fig. \ref{fig:F0T0}-\ref{fig:F3T3}  in   \ref{appendixsimulations}.
  Combining all the stability switching curves  together,  and zooming in the part of $(\tau_1,\tau_2)\in [0,5]\times[0,20]$, we get the Hopf bifurcation curves shown in Fig. \ref{fig:tau1tau2} a). Consider   the  bottom left region bounded by left-most curve of $\mathcal{T}^1_0$ and the lowest curve of $\mathcal{T}^2_0$ (see Fig. \ref{fig:tau1tau2} a)), which are   both part of $\mathcal{T}_0$. Since the crossing directions of the two switching curves (the black line and blue line) are all pointing outside of the region, the constant steady state $E^*$ is stable in the bottom left region.  Moreover,   the two stability switching curves  intersect at the point $( 2.21407,15.0019)$, which is the double Hopf bifurcation point, denoted by HH. For HH,  $\omega_1=0.1078$, $\omega_2=0.4920$. Using the normal form derivation process given in supplementary materials, we have
    $B_{11}=-0.0512$,
  $B_{21} = 0.0478$,
  $B_{13 }= 1.7736$,
  $B_{23} = -0.0706$,
  $B_{2100} = -0.7294$,
  $B_{1011}= -1.5826$,
  $B_{0021} = -0.2410$,
  $B_{1110} = -10.8681$.      Furthermore, we have the  normal form (\ref{normalformcylin}) with $\epsilon_1 =-1, \epsilon_2 =   -1, b =   6.5658,c = 14.8998, d =       1$, and  $d-bc =      -96.8306$.
 According to  chapter 7.5 in Ref. \cite{Guckenheimer}, case Ib arises, and   we have the bifurcation set near HH showing in Fig. \ref{fig:tau1tau2} b).  In region  $D_1$,  the  positive equilibrium is  asymptotically stable. In  region $D_2$ or  $D_6$, there is a stable periodic solution.  When the parameters cross  into the region D4,   there are two stable momogeneous periodic solutions coexisting in D4.

\section{Concluding remarks}
 The first part of this paper is mainly about    the  dynamics of the ODE predator-prey model (\ref{ODE system}) with  Allee effect in prey and cooperative hunting in predator.
  We start with  the ODE system in the case of weak cooperation.   Taking $p$, the convertion rate, as the  bifurcation parameter,  we prove that there is a loop of heteroclinic orbits between $E_1$ and $E_a$  formed by the intersection of the  stable manifold of $E_a$ and the unstable manifold of $E_1$  when $p=p^{\#}$ ($p^{\#}>1$).  It shows that $p^{\#}$ is an important threshold, which distinguishes the globally stability for $E_0$ and the existence of a  separatrix curve. To be specific, when $p>p^{\#}$, $E_0$ is globally stable, which means the extinction of both species, and we call the phenomenon overexploitation.  When $p<p^{\#}$, there is a  separatrix curve $\Gamma_p^s$, which is determined by the stable manifold of $E_a$.   High initial predator density will lead to extinction for both species,   and conversely low  initial predator density  will approach either a steady state or a periodic oscillation.  When $p<1$,  $E_1(1,0)$ is locally stable, which means that the prey achieves its  carrying capacity and the predator is extinct. When $1<p<p_H$, the stable state becomes $E^*$, which implies the coexistence of both species. When $p$ crosses $p_H$, the stable state switches into a limit cycle.   The period of the limit cycle increases   with $p$, and tends to $\infty$ as $p\rightarrow p^{\#}$.    The existence of a limit cycle with a large amplitude increases the risk of the extinction for the predator.

     In the case of strong cooperation,   three different cases are considered. We prove the  existence  of a loop of heteroclinic orbits between $E_1$ and $E_a$ when $p=p^{\#}$ in the first and second  cases and nonexistence of such a loop of heteroclinic orbits  in the third case.  We also exhibit  the existence of loop of heteroclinic orbits among $E_R^*$, $E_a$ and $E_1$  formed by  the unstable manifold of $E_R^*$ and the unstable manifold  $E_a$  when $p=\overline{p}^{\#}$,  and a homoclinic cycle when $p=p_{hom}$ by numerical simulations in case 3.
   When $p>p^{\#}$ in case 1 and case 2 ($p>\overline{p}^{\#}$ in case 3), both species distinct.  When $1<p<p^{\#}$ in case 1 and case 2,  there is a  separatrix curve $\Gamma_p^s(E_a)$, and there are two stable states coexisting, which are extinction for both species and coexistence or oscillation. When $p<p^{\#}$ (or $p<\overline{p}^{\#}$) and $ p_{SN}<p<1$, there are two separatrix curves, the stable manifold of $E_a$ and the stable manifold of $E_R^*$,  separating  the first quadrant into three parts. If the initial population is above the stable manifold of $E_a$, both species are extinct; if the initial population is inside the stable manifold of $E_R^*$, the stable state can be $E^*$ or limit cycle;  if the initial population is below the stable manifold of $E_a$ and outside the stable manifold of the stable manifold of $E_R^*$,  the prey achieves its  carrying capacity and the predator is extinct.
    When $p<p_{SN}$, there is one  separatrix curve $\Gamma_p^s(E_a)$,   and there are two stable states coexisting.  The species with initial population above the separatrix finally become extinct, however only the prey survive when  initial population is below it.

There is a unique equilibrium for weak cooperation when $1<p<\frac{1}{a}$.  It is shown that there are at most two stable states coexist for different initial values separated by the stable manifold of $E_a$.  There are two interior equilibria $E^*$ and $E_R^*$ in the case of strong cooperation when $ p_{SN}<p<1$.   The equilibrium $E_R^*$ brings one more   separatrix curve than the case of weak cooperation,  thus there may be three stable states coexisting.

  In the second part, we consider the corresponding diffusive system, and focus on the effect of diffusion on the dynamics. Taking diffusion coefficients $d_1$ and $d_2$ as bifurcation parameters,  we give the conditions of existence of Turing instability and Turing-Hopf bifurcation.    We illustrate complex dynamics of the diffusive system, including the existence of spatially inhomogeneous steady state, coexistence of two spatially inhomogeneous periodic solutions.

   The third part of this paper is about the diffusive system with two delays. We discuss  the joint effect of two delays on the dynamical behavior of the diffusive system. Applying the method of stability switching curves, we find all the stability switching curves at which the characteristic roots are purely imaginary.  Combining the geometrical and analytic method, we can decide the  crossing direction of the characteristic roots as long as we confirm the positive direction for stability switching curves. Then we can get the condition of existence of Hopf bifurcation. By searching the intersection of  stability switching curves near the stable region, we get the double Hopf bifurcation point. Through the calculation of normal form of the system,  we get  the corresponding unfolding system and the bifurcation set, from which we can figure out the complete dynamics near the bifurcation point on the $(\tau_1,\tau_2)$ parametric plane. We theoretically prove and illustrate that near the bifurcation point, there are the phenomena of the stability of positive equilibrium, stable periodic solutions and   coexistence of two periodic solutions.

 \appendix

 \section{Stability switching curves}
 \label{sectionStability switching curves}

 Applying  the method proposed in  Lin and Wang \cite{Lin X} to find the stability switching curves, we need to verify the assumptions (i)-(iv) in Ref.  \cite{Lin X}.
  \begin{itemize}

 \item[(i)] Finite number of characteristic roots on $\mathbb{C}_+=\{\lambda\in\mathbb{C}:{\rm Re}\lambda>0\}$ under the condition $${\rm deg} (P_{0,n}(\lambda))\geq {\rm max}\{{\rm deg}(P_{1,n}(\lambda)),{\rm deg}(P_{2,n}(\lambda)),{\rm deg}(P_{3,n}(\lambda))\}.$$

 \item[ (ii)] $P_{0,n}(0)+P_{1,n}(0)+P_{2,n}(0)+P_{3,n}(0)\neq 0$.
 \item[(iii)]
 $P_{0,n}(\lambda),P_{1,n}(\lambda),P_{2,n}(\lambda),P_{3,n}(\lambda)$ are coprime polynomials.
  \item[(iv)] $\lim\limits_{\lambda\rightarrow\infty} \left(\left| \frac{P_{1,n}(\lambda)}{P_{0,n}(\lambda)}\right|+  \left| \frac{P_{2,n}(\lambda)}{P_{0,n}(\lambda)}\right|+\left| \frac{P_{3,n}(\lambda)}{P_{0,n}(\lambda)}\right|\right) <1$.
 \end{itemize}
 In fact, condition  (i) follows from  Ref. \cite{J. Hale}. In fact   $P_{0,n}(0)+P_{1,n}(0)+P_{2,n}(0)+P_{3,n}(0)=J_n$, which is determined in (\ref{TnJn}).   Since we assume $p<p_H$ and $d_2>d_2^T(n,d_1)$ for all $n\in\mathbb{N}$, from  Theorem \ref{cen0den0},
  the constant steady state  $E^*$    is locally asymptotically stable, and $J_n>0$ for all $n\in\mathbb{N}$.  Thus,
  condition (ii) satisfies. From the expressions of  $P_{0,n}(\lambda),P_{1,n}(\lambda),P_{2,n}(\lambda),P_{3,n}(\lambda)$ in (\ref{p1p2p3}), condition (iii) is  obviously satisfied. From (\ref{p1p2p3}), we know
  $$
  \lim\limits_{\lambda\rightarrow\infty} \left(\left| \frac{P_{1,n}(\lambda)}{P_{0,n}(\lambda)}\right|\right) =\lim\limits_{\lambda\rightarrow\infty} \left(  \left| \frac{P_{2,n}(\lambda)}{P_{0,n}(\lambda)}\right|\right) =\lim\limits_{\lambda\rightarrow\infty} \left( \left| \frac{P_{3,n}(\lambda)}{P_{0,n}(\lambda)}\right|\right) =0,
  $$
which means that (iv) are  satisfied.

From \cite{K.L. Cooke,Ruan wei},
 we have the following lemma.
    \begin{lemma}
 As the delays $(\tau_1,\tau_2)$ vary continuously in $\mathbb{R}_+^2$, the number of zeros (counting multiplicity) of $D_n(\lambda; \tau_1,\tau_2)$ on $\mathbb{C}_+$ can change only if a zero appears on or cross the imaginary axis.
 \end{lemma}

From condition (ii), $\lambda=0$ is not a   root of (\ref{character}). Now we  are in the position of seeking  all the points $(\tau_1,\tau_2)$ such that $D_n(\lambda; \tau_1,\tau_2)=0$ has at least one root on the imaginary axis,  on which the stability of equilibrium $E^*$ may  switch.    Substituting $\lambda =i\omega~(\omega>0)$ into (\ref{character}),
  \begin{equation*}
  \label{characteriomega}
   (P_{0,n}(i\omega)+P_{1,n}(i\omega)e^{-i\omega\tau_1})+(P_{2,n}(i\omega)+P_{3,n}(i\omega)e^{-i\omega\tau_1})e^{-i\omega\tau_2}=0.
  \end{equation*}
 Due to  $|e^{-i\omega\tau_2}|=1$, we get that
  \begin{equation*}
  |P_{0,n}(i\omega)+P_{1,n}(i\omega)e^{-i\omega\tau_1}|=|P_{2,n}(i\omega)+P_{3,n}(i\omega)e^{-i\omega\tau_1}|,
  \end{equation*}
  Thus,  we have the following equality
  \begin{equation}\label{omegatau1}
  |P_{0,n}(i\omega)|^2+|P_{1,n}(i\omega)|^2-|P_{2,n}(i\omega)|^2-|P_{3,n}(i\omega)|^2=2A_{1,n}(\omega)\cos(\omega\tau_1)-2B_{1,n}(\omega)\sin(\omega\tau_1),
  \end{equation}
  where
  \begin{equation*}
  \begin{aligned}
  A_{1,n}(\omega)={\rm Re}(P_{2,n}(i\omega)\overline{P}_{3,n}(i\omega)-P_{0,n}(i\omega)\overline{P}_{1,n}(i\omega))=\sqrt{A_{1,n}(\omega)^2+B_{1,n}(\omega)^2}\cos(\varphi_{1,n}(\omega)),\\
   B_{1,n}(\omega)={\rm Im}(P_{2,n}(i\omega)\overline{P}_{3,n}(i\omega)-P_{0,n}(i\omega)\overline{P}_{1,n}(i\omega))=\sqrt{A_{1,n}(\omega)^2+B_{1,n}(\omega)^2}\sin(\varphi_{1,n}(\omega)),\\
  \end{aligned}
  \end{equation*}
if $A_{1,n}(\omega)^2+B_{1,n}(\omega)^2>0$, with $\varphi_{1,n}(\omega)=\angle\{P_{2,n}(i\omega)\overline{P}_{3,n}(i\omega)-P_{0,n}(i\omega)\overline{P}_{1,n}(i\omega)\}\in(-\pi,\pi]$. Thus, we can write Eq. (\ref{omegatau1})  as
   \begin{equation}\label{omegatau1xin}
    |P_{0,n}(i\omega)|^2+|P_{1,n}(i\omega)|^2-|P_{2,n}(i\omega)|^2-|P_{3,n}(i\omega)|^2 =2\sqrt{A_{1,n}(\omega)^2+B_{1,n}(\omega)^2}\cos(\varphi_{1,n}(\omega)+\omega\tau_1).
    \end{equation}

Denote
\begin{equation}\label{conditiontau1}
 \Sigma^1_n=\left\lbrace  \omega\in\mathbb{R}_+:\left(  (P_{0,n}(i\omega)|^2+|P_{1,n}(i\omega)|^2-|P_{2,n}(i\omega)|^2-|P_{3,n}(i\omega)|^2\right)^2 \leq 4(A_{1,n}(\omega)^2+B_{1,n}(\omega)^2)\right\rbrace .
 \end{equation}
 We can easily know that there is  $\tau_1\in \mathbb{R}_+$ satisfying Eq. (\ref{omegatau1xin}) if and only if  $\omega\in\Sigma^1_n$.
  In fact, (\ref{conditiontau1}) also includes the case $A_{1,n}^2(\omega)+B_{1,n}^2(\omega)=0$.

Define
  \begin{equation*}
 \theta_{1,n}(\omega)=\arccos\left( \dfrac{ |P_{0,n}(i\omega)|^2+|P_{1,n}(i\omega)|^2-|P_{2,n}(i\omega)|^2-|P_{3,n}(i\omega)|^2}{2\sqrt{A_{1,n}(\omega)^2+B_{1,n}(\omega)^2}}\right) ,~~~~~\theta_{1,n}\in[0,\pi],
  \end{equation*}
  which leads to
  \begin{equation}\label{tau1}
  \tau_{1,j_1,n}^{\pm}(\omega)=\dfrac{\pm\theta_{1,n}(\omega)-\varphi_{1,n}(\omega)+2j_1\pi}{\omega},~~~j_1\in\mathbb{Z}.
  \end{equation}
 Using the same method, we can get  the corresponding results on the other delay $\tau_2$,
   \begin{equation}\label{tau2}
    \tau_{2,j_2,n}^{\pm}(\omega)=\dfrac{\pm\theta_{2,n}(\omega)-\varphi_{2,n}(\omega)+2j_2\pi}{\omega},~~~j_2\in\mathbb{Z},
    \end{equation}
   where
   \begin{equation*}
    \begin{array}{c}
     \theta_{2,n}(\omega)=\arccos\left( \dfrac{|P_{0,n}(i\omega)|^2-|P_{1,n}(i\omega)|^2+|P_{2,n}(i\omega)|^2-|P_{3,n}(i\omega)|^2}{2\sqrt{A_{2,n}(\omega)^2+B_{2,n}(\omega)^2}}\right) ,~~~~~\theta_{2,n}\in[0,\pi],\\
    A_{2,n}(\omega)
    =2\sqrt{A_{2,n}(\omega)^2+B_{2,n}(\omega)^2}\cos(\varphi_{2,n}(\omega)),\\
     B_{2,n}(\omega)
     =2\sqrt{A_{2,n}(\omega)^2+B_{2,n}(\omega)^2}\sin(\varphi_{2,n}(\omega)),\\
    \end{array}
    \end{equation*}
and the condition on $\omega$ is as follows
 \begin{equation}\label{conditiontau2}
  \Sigma^2_n=\left\lbrace    \omega\in\mathbb{R}_+:\ \left(  |P_{0,n}(i\omega)|^2-|P_{1,n}(i\omega)|^2+|P_{2,n}(i\omega)|^2-|P_{3,n}(i\omega)|^2\right)^2 \leq 4(A_{2,n}(\omega)^2+B_{2,n}(\omega)^2)\right\rbrace .
 \end{equation}
 In fact, it is  easy to show that the   inequality in  (\ref{conditiontau1}) is equivalent  to the one in  (\ref{conditiontau2}) by squaring both sides,  thus,  $\Sigma^1_n=\Sigma^2_n$, and we denote both of them as $\Omega_n$.

We call the set
\begin{equation*}\begin{array}{r}
\Omega_n=\Bigg\{     \omega\in \mathbb{R}_+:F_n(\omega)\stackrel{\vartriangle}{=} ( |P_{0,n}(i\omega)|^2+|P_{1,n}(i\omega)|^2-|P_{2,n}(i\omega)|^2-|P_{3,n}(i\omega)|^2)^2\\ \left.-4(A_{1,n}(\omega)^2+B_{1,n}(\omega)^2)\leq 0  \Bigg\}\right.
\end{array}
\end{equation*}
 the crossing set of characteristic equation $D_n(\lambda;\tau_1,\tau_2)=0$.

  Obviously, $F_n(\omega)=0$  has a finite number of roots on $\mathbb{R}_+$.  If $F_n(0)>0$,  denote the roots of  $F_n(\omega)=0$ as
$$
0<a_{1,n}<b_{1,n}\leq a_{2,n}<b_{2,n}<\cdots\leq a_{m,n}<b_{m,n}<+\infty,
$$
then we have  $\Omega_n=\bigcup\limits_{k=1}^m\Omega_{k,n},~~ \Omega_{k,n}=[a_{k,n},b_{k,n}].$
 If $F_n(0)\leq 0$, denote the roots of  $F_n(\omega)=0$ as
$$
 0<b_{1,n}\leq a_{2,n}<b_{2,n}<\cdots\leq a_{N,n}<b_{N,n}<+\infty,
$$
then we have
$$
\Omega_n=\bigcup\limits_{k=1}^N\Omega_{k,n},~~ \Omega_{1,n}=\left( 0,b_{1,n}\right] , \Omega_{k,n}=[a_{k,n},b_{k,n}]~~  (k\geq 2).
$$

In fact, we can  confirm that   when $\tau_1=\tau_{1,j_1,n}^+(\omega)$, we have $\tau_2=\tau_{2,j_2,n}^-(\omega)$, and when $\tau_1=\tau_{1,j_1,n}^-(\omega)$, we have $\tau_2=\tau_{2,j_2,n}^+(\omega)$.   On two ends of $\Omega_{j,n}$, we have $F_n(a_{j,n})=F_n(b_{j,n})=0$,  and thus
$
\theta_{i,n}(a_{j,n})=\delta_i^a\pi,~~\theta_{i,n}(b_{j,n})=\delta_i^b\pi,
$
where $\delta_i^a,\delta_i^b=0,1,i=1,2$. From (\ref{tau1}) and (\ref{tau2}), we can easily verify  that
\begin{equation}\label{connect}
\begin{aligned}
(\tau_{1,j_1,n}^{+j}(a_{j,n}),\tau_{2,j_2,n}^{-j}(a_{j,n}))=(\tau_{1,j_1+\delta_1^a,n}^{-j}(a_{j,n}),\tau_{2,j_2-\delta_2^a,n}^{+j}(a_{j,n})),\\
(\tau_{1,j_1,n}^{+j}(b_{j,n}),\tau_{2,j_2,n}^{-j}(b_{j,n}))=(\tau_{1,j_1+\delta_1^b,n}^{-j}(b_{j,n}),\tau_{2,j_2-\delta_2^b,n}^{+j}(b_{j,n})).\\
\end{aligned}
\end{equation}

Denote
\begin{equation}\label{Tzfk}
\begin{aligned}
\mathcal{T}_{j_1,j_2,n}^{\pm j}&=\left\lbrace \left( \tau_{1,j_1,n}^{\pm}(\omega),\tau_{2,j_2,n}^{\mp}(\omega)\right):\omega\in\Omega_{j,n} \right\rbrace \\&=\left\lbrace \left( \dfrac{\pm\theta_{1,n}(\omega)-\varphi_{1,n}(\omega)+2j_1\pi}{\omega},\dfrac{\mp\theta_{2,n}(\omega)-\varphi_{2,n}(\omega)+2j_2\pi}{\omega}\right):\omega\in\Omega_{j,n} \right\rbrace.
\end{aligned}
\end{equation}
Thus, for the stability switching curves corresponding to $\Omega_{j,n}$, $\mathcal{T}_{j_1,j_2,n}^{+j}$ is  connected to  $\mathcal{T}_{j_1+\delta_1^a,j_2-\delta_2^a,n}^{-j}$ at one end $a_{j,n}$,  and connected to  $\mathcal{T}_{j_1+\delta_1^b,j_2-\delta_2^b,n}^{-j}$ at the other end $b_{j,n}$.

Denote
\begin{equation}\label{Tk}
\mathcal{T}^{j}_n=\bigcup_{j_1=-\infty}^{\infty}\bigcup_{j_2=-\infty}^{\infty}(\mathcal{T}_{j_1,j_2,n}^{+j}\cup\mathcal{T}_{j_1,j_2,n}^{-j})\cap \mathbb{R}_+^2,
\end{equation}
and
\begin{equation}\label{huatn}\mathcal{T}_n=\bigcup_{j=1}^N\mathcal{T}^j_n.\end{equation}

 \begin{definition}
  Any $(\tau_1,\tau_2)\in \mathcal{T}_n$ is called a crossing point, which makes $D_n(\lambda;\tau_1,\tau_2)=0$ have at least one purely imaginary root $i\omega$ with $\omega$ belongs to the crossing set $\Omega_n$.  The set $\mathcal{T}_n$ is called stability switching curves.
 \end{definition}

\section{Crossing directions}\label{sectioncrossing}
When $(\tau_1,\tau_2)$ varies and crosses the stability switching curves from one side to the other, the number of characteristic roots with positive real part may increase. we call it the crossing direction of stability switching curves.

In order to describe the crossing direction clearly, we need to specify the positive direction for  stability switching curves $\mathcal{T}_{j_1,j_2,n}^{\pm j}$.
 We call the direction of the curve corresponding to increasing $\omega\in \Omega_{j,n}$  the positive direction, i.e. from $(\tau_{1,j_1,n}^{\pm j}(a_{j,n}),\tau_{2,j_2,n}^{\mp j}(a_{j,n}))$ to $(\tau_{1,j_1,n}^{\pm j }(b_{j,n}),\tau_{2,j_2,n}^{\mp j}(b_{j,n}))$.  From the fact that $\mathcal{T}_{j_1,j_2,n}^{+j}$ is  connected to  $\mathcal{T}_{j_1+\delta_1^a,j_2-\delta_2^a,n}^{-j}$ at  $a_{j,n}$ as we have mentioned  in the previous subsection,   the positive direction of the two curves are opposite.

\begin{figure}[h]
 	\centering
 	\includegraphics[width=0.6\textwidth,height=0.35\textwidth]{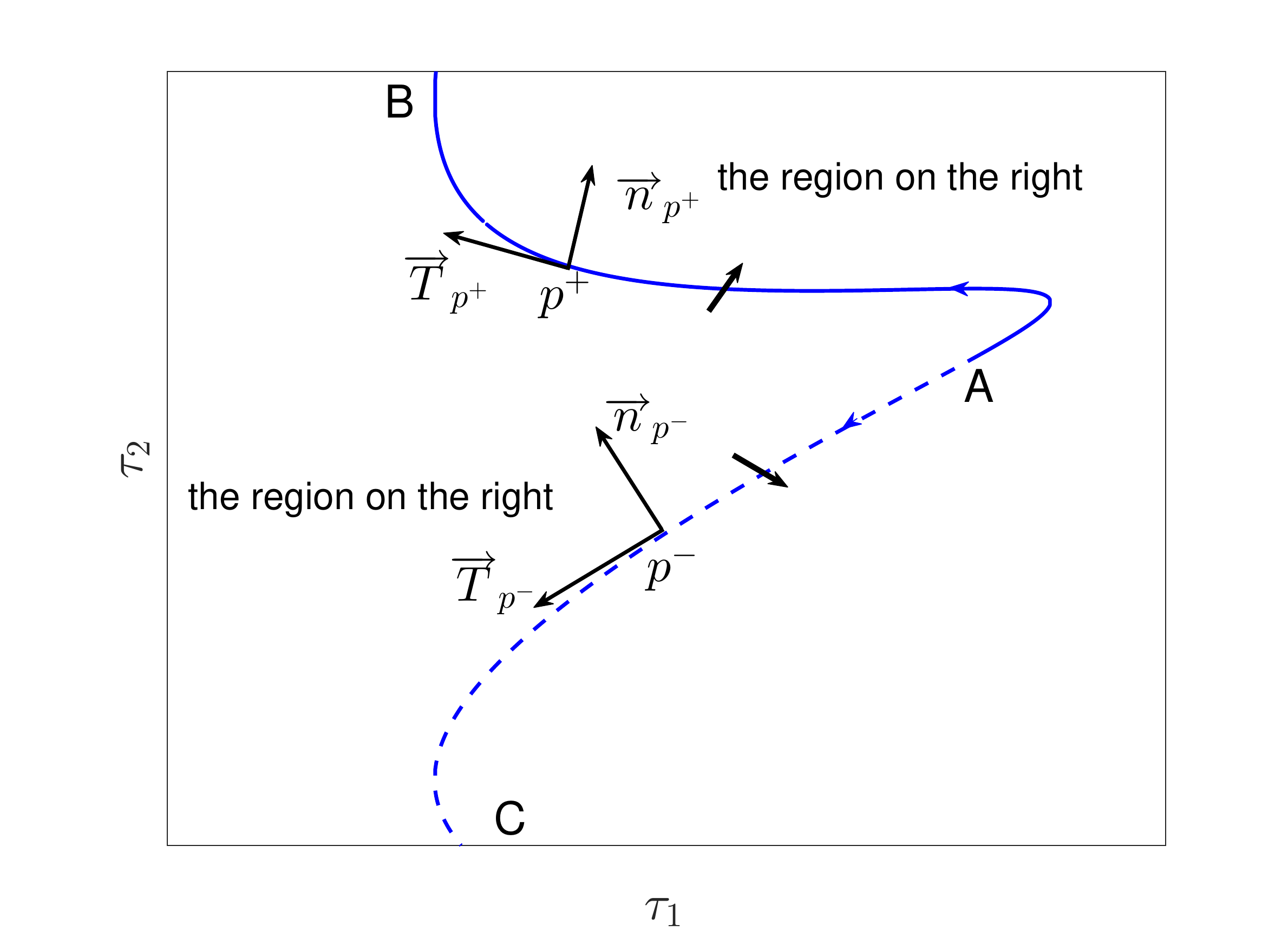}
 	\caption  {The positive direction and crossing direction of a part of stability switching curves corresponding to $\Omega_{j,n}=[a_{j,n},b_{j,n}]$.  }
 	\label{fig:neighbor}
 \end{figure}

 To make our  expression clear, we draw a schematic diagram of a part of stability switching curves corresponding to $\Omega_{j,n}=[a_{j,n},b_{j,n}]$ (see Fig. \ref{fig:neighbor}).
 The solid curve stands for $\mathcal{T}_{j_1,j_2,n}^{+ j}$, with two ends $A(\tau_{1,j_1,n}^{+j}(a_{j,n}),\tau_{2,j_2,n}^{-j}(a_{j,n}))$, and $B(\tau_{1,j_1,n}^{+j}(b_{j,n}),\tau_{2,j_2,n}^{-j}(b_{j,n}))$.   The dashed curve denotes $\mathcal{T}_{j_1+\delta_1^a,j_2-\delta_2^a,n}^{- j}$, which is   connected to $\mathcal{T}_{j_1,j_2,n}^{+ j}$ at $A$ corresponding to $a_{j,n}$, with the positive direction from $A$ to $C$.

 Call the region on the right-hand (left-hand) side as we move in the positive directions of the curve  the region on the right (left).  Since the tangent vector of $\mathcal{T}_{j_1,j_2,n}^{\pm j}$ at $p^{\pm}(\tau_{1,j_1,n}^{\pm},\tau_{2,j_2,n}^{\mp})$ along the positive direction is $(\frac{\partial \tau_1}{\partial \omega},\frac{\partial \tau_2}{\partial \omega})\mid_{p^{\pm}}\stackrel{\vartriangle}{=}\overrightarrow{T}_{p^{\pm}}$, the normal vector of $\mathcal{T}_{j_1,j_2,n}^{\pm j}$  pointing to the right region is $(\frac{\partial \tau_2}{\partial \omega},-\frac{\partial \tau_1}{\partial \omega})\mid_{p^{\pm}}\stackrel{\vartriangle}{=}\overrightarrow{n}_{p^{\pm}}$ (see Fig. \ref{fig:neighbor}). On the other hand, when a pair of complex characteristic roots cross the imaginary axis to the right half plane, $(\tau_1,\tau_2)$ moves along the direction
 $(\frac{\partial \tau_1}{\partial \sigma},\frac{\partial \tau_2}{\partial \sigma})\mid_{p^{\pm}}$.  The inner product of these two vectors  is
 \begin{equation}\label{delta}
 \delta(\omega)\mid_{p^{\pm}}:=\frac{\partial \tau_1}{\partial \sigma}\frac{\partial \tau_2}{\partial \omega}-\frac{\partial \tau_2}{\partial \sigma}\frac{\partial \tau_1}{\partial \omega}\mid_{p^{\pm}},
 \end{equation}
 the sign of which can decide the  crossing direction of the   characteristic roots.

Denote $\lambda=\sigma+i\omega$. Write  (\ref{character}) as $${\rm Re}D_n(\lambda;\tau_1,\tau_2)+i{\rm Im}D_n(\lambda;\tau_1,\tau_2)=0.$$ From the implicit function theorem, $\tau_1$, $\tau_2$ can be expressed as functions of $\sigma$ and $\omega$, and
 if $
{\rm det}\left(\begin{array}{cc}
R_1& R_2\\I_1&I_2
\end{array}\right)=R_1I_2-R_2I_1\neq 0,
$
we have
\begin{equation}\label{juzhen}
\Delta(\omega):=\left( \begin{array}{cc}
\frac{\partial \tau_1}{\partial \sigma}&\frac{\partial \tau_1}{\partial \omega}\\\frac{\partial \tau_2}{\partial \sigma}&\frac{\partial \tau_2}{\partial \omega}
\end{array}\right) \arrowvert_{\sigma=0,\omega\in \Omega_n}=-\left( \begin{array}{cc}
R_1& R_2\\I_1&I_2
\end{array}\right)^{-1}\left( \begin{array}{cc}
R_0&-I_0\\I_0&R_0
\end{array}\right),
\end{equation}
where
\begin{equation*}
\begin{aligned}
&\dfrac{\partial {\rm Re} D_n(\lambda;\tau_1,\tau_2)}{\partial \sigma}|_{\lambda=i\omega}
=\dfrac{\partial {\rm Im} D_n(\lambda;\tau_1,\tau_2)}{\partial \omega}|_{\lambda=i\omega}=R_0,
\end{aligned}
\end{equation*}
\begin{equation*}
\dfrac{\partial {\rm Re} D_n(\lambda;\tau_1,\tau_2)}{\partial \omega}|_{\lambda=i\omega}=-\dfrac{\partial {\rm Im} D_n(\lambda;\tau_1,\tau_2)}{\partial \sigma}|_{\lambda=i\omega}=-I_0,
\end{equation*}
\begin{equation*}
\begin{aligned}
&\dfrac{\partial {\rm Re} D_n(\lambda;\tau_1,\tau_2)}{\partial \tau_l}|_{\lambda=i\omega} 
=R_l,
\dfrac{\partial {\rm Im} D_n(\lambda;\tau_1,\tau_2)}{\partial \tau_l}|_{\lambda=i\omega} 
=I_l,
\end{aligned}
\end{equation*}
with $l=1,2$.
Here   we do not consider the case that $i\omega$ is the multiple root of $D_n(\lambda;\tau_1,\tau_2)=0$,  i.e.,  $\frac{d D_n(\lambda;\tau_1,\tau_2)}{d \lambda}\mid_{\lambda=i\omega}=R_0+iI_0\neq 0$, thus    we have
    $
{\rm det}\left(\begin{array}{cc}
-R_0& I_0\\-I_0&-R_0
\end{array}\right)=R_0^2+I_0^2>0. $
Therefore,  the sign of $\delta(\omega)={\rm det} \Delta(\omega)$  is decided by $R_1I_2-R_2I_1$.
We can verify that
$
R_1I_2-R_2I_1\mid_{p^{\pm}} 
=\pm\omega^2 |P_{2,n}\overline{P}_{3,n}-P_{0,n}\overline{P}_{1,n}|\sin \theta_{1,n}$.
From  $\theta_{1,n}( \mathring{\Omega}_{j,n})\subset (0,\pi)$, we have
$$\begin{array}{l}
\delta(\omega\in \mathring{\Omega}_{j,n})\mid_{p^+}>0,~~\forall p^+\in\mathcal{T}_{j_1,j_2,n}^{+ j},~ {\rm and } ~
\delta(\omega\in \mathring{\Omega}_{j,n})\mid_{p^-}<0,~~\forall p^-\in\mathcal{T}_{j_1,j_2,n}^{- j},
\end{array}$$
where $\mathring{\Omega}_{j,n}$ denotes the interior of $\Omega_{j,n}$.

\begin{lemma}\label{direction}
For any $j=1,2,\cdots,N$, we have
\begin{equation*}
\delta(\omega\in \mathring{\Omega}_{j,n}) >0 ,~~\forall (\tau_1(\omega),\tau_2(\omega))\in\mathcal{T}_{j_1,j_2,n}^{+ j}, ~{\rm and~ }
\delta(\omega\in \mathring{\Omega}_{j,n})<0,~~\forall (\tau_1(\omega),\tau_2(\omega))\in\mathcal{T}_{j_1,j_2,n}^{- j}.
\end{equation*}
\end{lemma}
Thus, we can get  a further conclusion which is more  intuitive.
\begin{remark}
The region on the right  of $\mathcal{T}_{j_1,j_2,n}^{+j}$   has two more   characteristic roots with positive real parts, and the region on the right  of   $\mathcal{T}_{j_1,j_2,n}^{- j}$ has two  less characteristic roots with positive real parts.
\end{remark}



 \section{ The stability switching curves in  Fig. \ref{fig:tau1tau2} in section \ref{section numerical delays} } \label{appendixsimulations}
 Following the steps in   \ref{sectionStability switching curves} and \ref{sectioncrossing}, we can draw all the stability curves on $(\tau_1,\tau_2)$ plane, and decide the crossing direction.

\begin{figure}
     \centering
  a) \includegraphics[width=0.3 \textwidth]{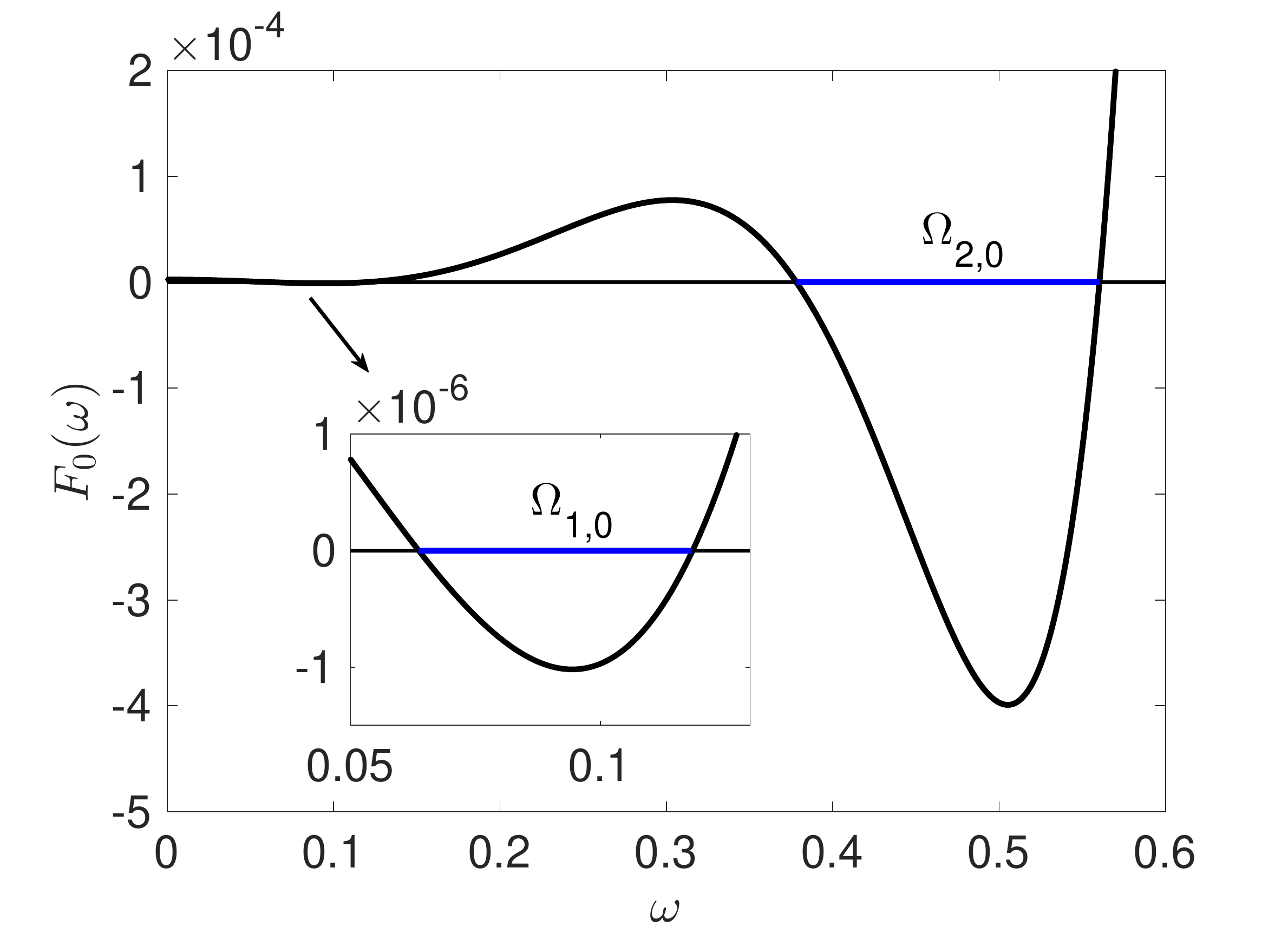}
    b)   \centering
  \includegraphics[width=0.3 \textwidth]{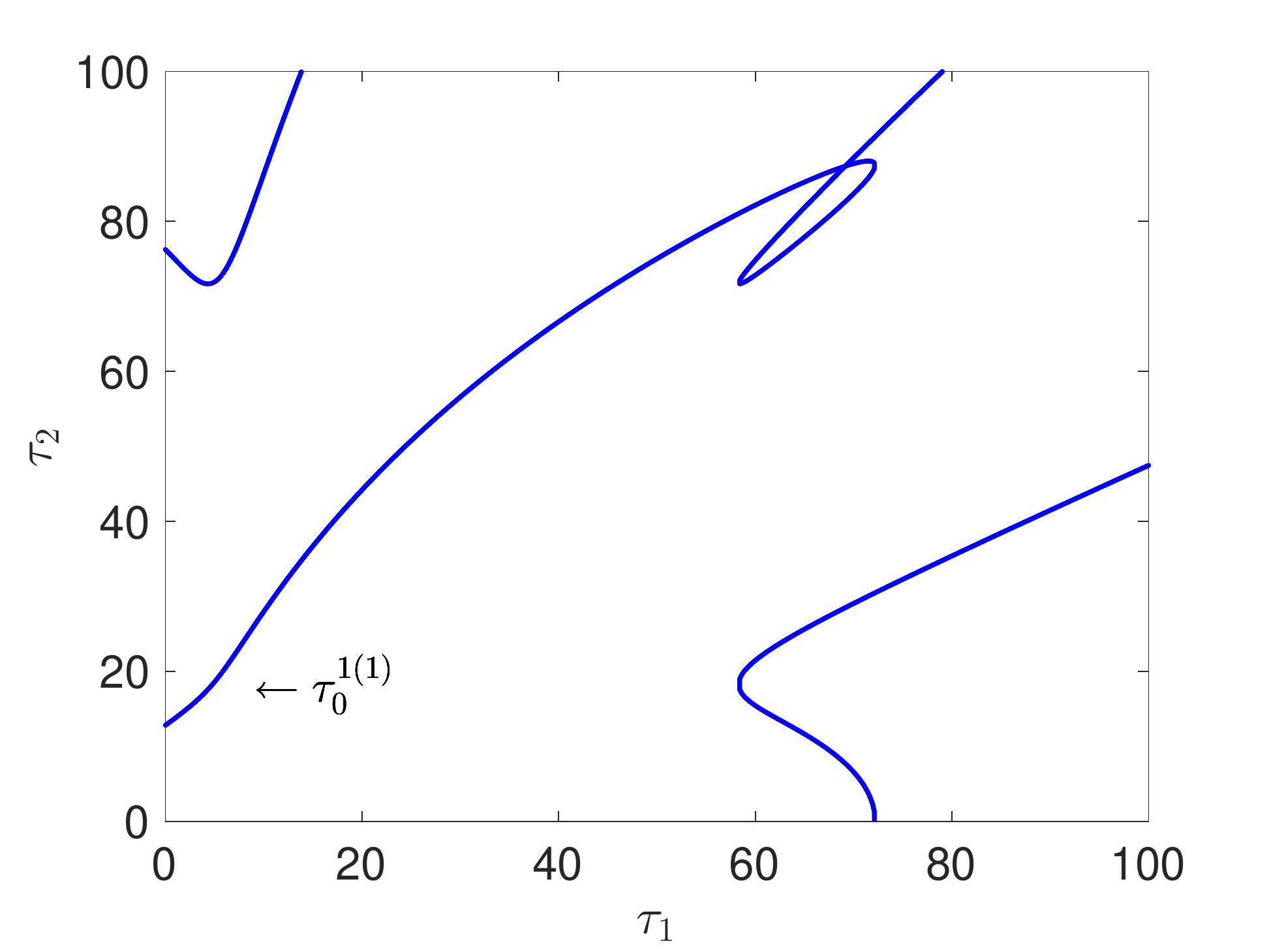}
  c) \includegraphics[width=0.3 \textwidth]{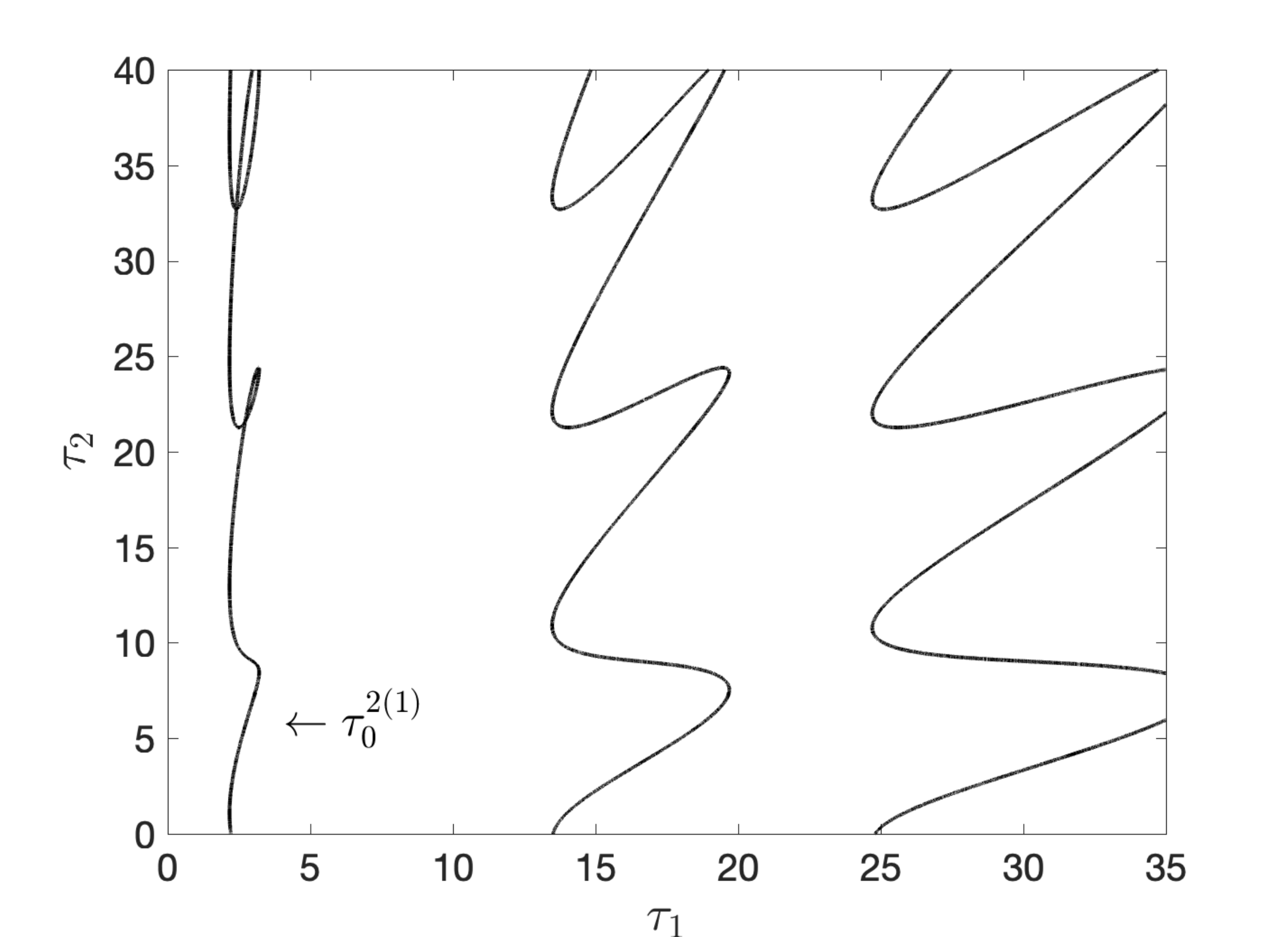}
 \caption{a) The crossing set   $\Omega_{1,0}\bigcup\Omega_{2,0}$. b) Stability switching curves $\mathcal{T}^1_0$. c) Stability switching curves $\mathcal{T}^2_0$.}
  \label{fig:F0T0}
 \end{figure}

\begin{figure}
     \centering
  a) \includegraphics[width=0.45\textwidth,height=0.40\textwidth]{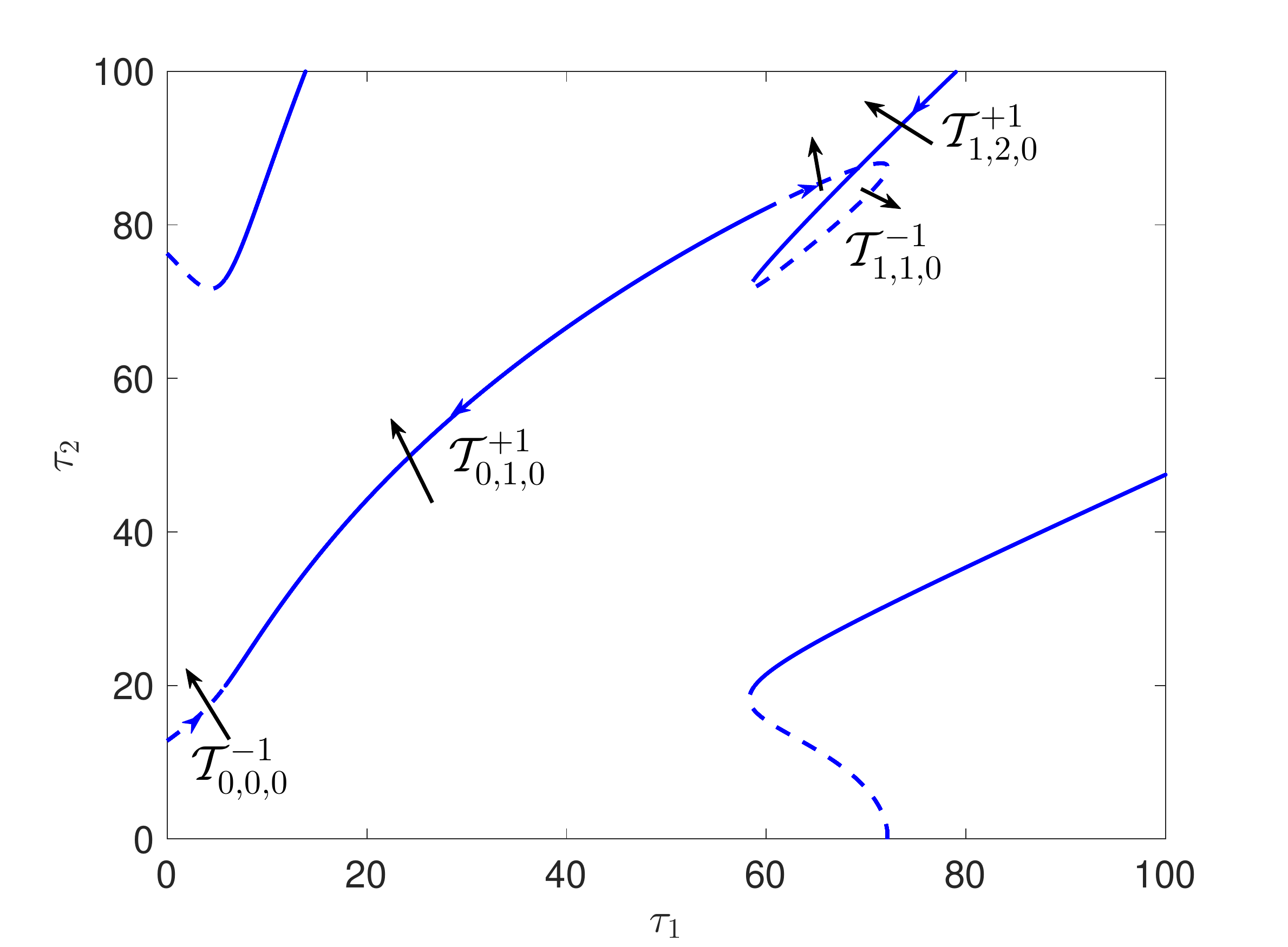}
  b) \includegraphics[width=0.45\textwidth,height=0.40\textwidth]{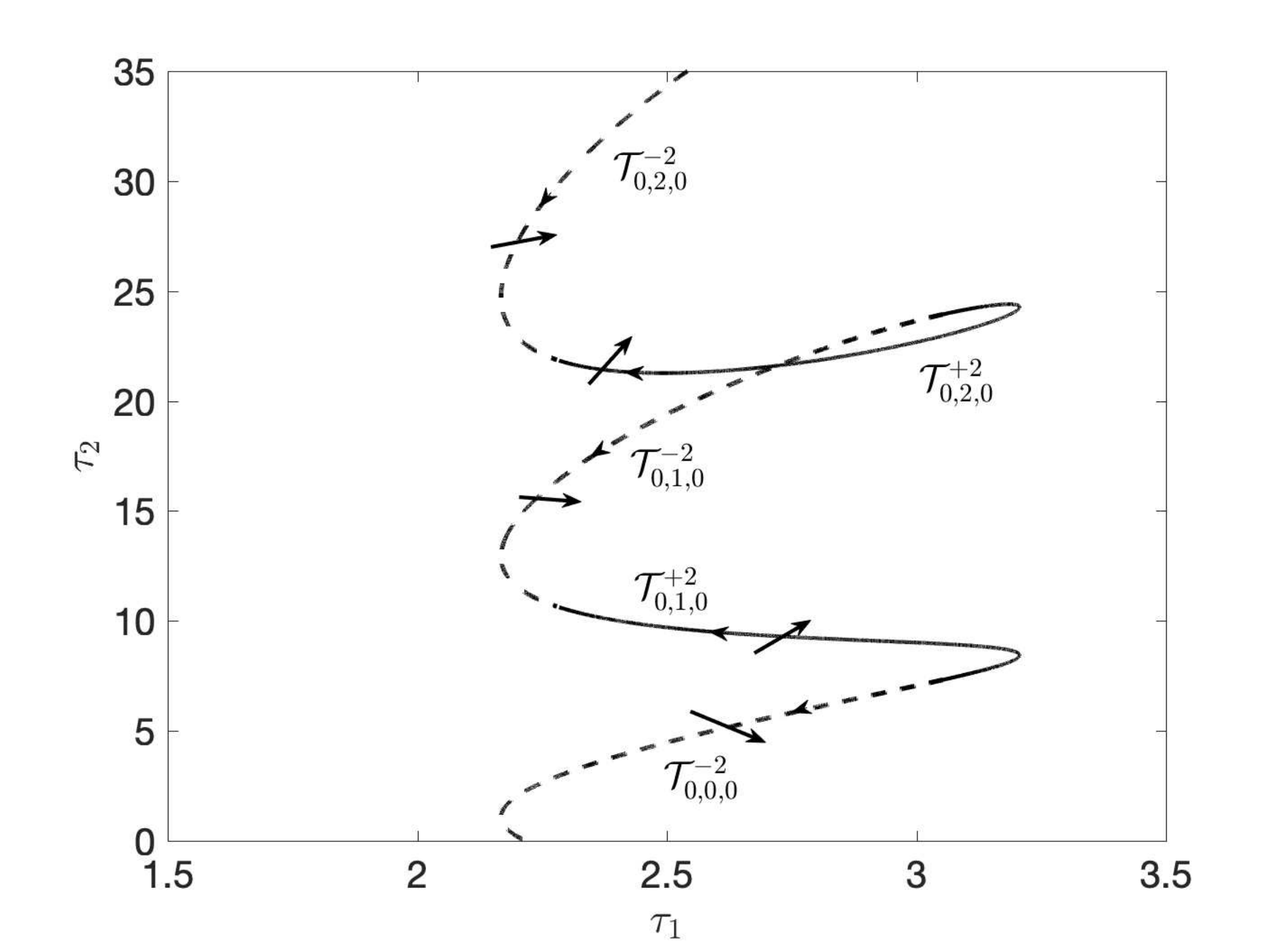}
 \caption{a) The detailed structure of  $\tau_0^{1(1)}$ in Fig. \ref{fig:F0T0} b).  b) The detailed structure of  $\tau_0^{2(1)}$ in Fig. \ref{fig:F0T0} c).  }
  \label{fig:T0direction}
 \end{figure}

  We can verify that  $F_0(0)>0$ and $F_0(\omega)=0$ has four roots (see  Fig. \ref{fig:F0T0} a)). Thus the crossing set $\Omega_0=$ $\Omega_{1,0}\bigcup\Omega_{2,0}=[a_{1,0},b_{1,0}]\bigcup[a_{2,0},b_{2,0}]$$=[0.0636,0.1184
]\bigcup[ 0.3781,0.5602]$.  From (\ref{Tzfk}) and (\ref{Tk}), we can get the stability switching curves $\mathcal{T}^1_0$ corresponding to $\Omega_{1,0}$ which is shown in Fig. \ref{fig:F0T0}  b).   
To show the structure of stability switching curves and the crossing direction clearly, we take a part of curves of $\mathcal{T}_0^1$ near the origin (i.e. $\tau_0^{1(1)}$ in  Fig. \ref{fig:F0T0}  b)) as an example, and draw the figure in Fig. \ref{fig:T0direction} a). From bottom to top, it starts with a part of $\mathcal{T}_{0,0,0}^{-1}$, which is connected to $\mathcal{T}_{0,1,0}^{+1}$ at $b_{1,0}$.  $\mathcal{T}_{0,1,0}^{+1}$ is linked to $\mathcal{T}_{1,1,0}^{-1}$ at $a_{1,0}$, which is again connected to $\mathcal{T}_{1,2,0}^{+1}$ at $b_{1,0}$ $\cdots$.    The numerical results support the analysis result in (\ref{connect}).
 Similarly,   the stability switching curves $\mathcal{T}^2_0$ corresponding to $\Omega_{2,0}$  are shown in Fig. \ref{fig:F0T0}  c). All the stability switching curves for $n=0$ are given by $\mathcal{T}_0=\mathcal{T}^1_0\cup\mathcal{T}^2_0$.
    We also draw the leftmost curve of $\mathcal{T}_0^2$ (marked $\tau_0^{2(1)}$ in Fig. \ref{fig:F0T0} c))  in  Fig. \ref{fig:T0direction} b).  In Fig. \ref{fig:T0direction}, the  arrows on the stability switching curves  represent their positive direction.  From Lemma \ref{direction}, we know that the regions on the right (left) of the solid (dashed) curves, which the black arrows point to, have two more characteristic roots with positive real parts.

\begin{figure}
 \centering
 a) \includegraphics[width=0.45\textwidth,height=0.30\textwidth]{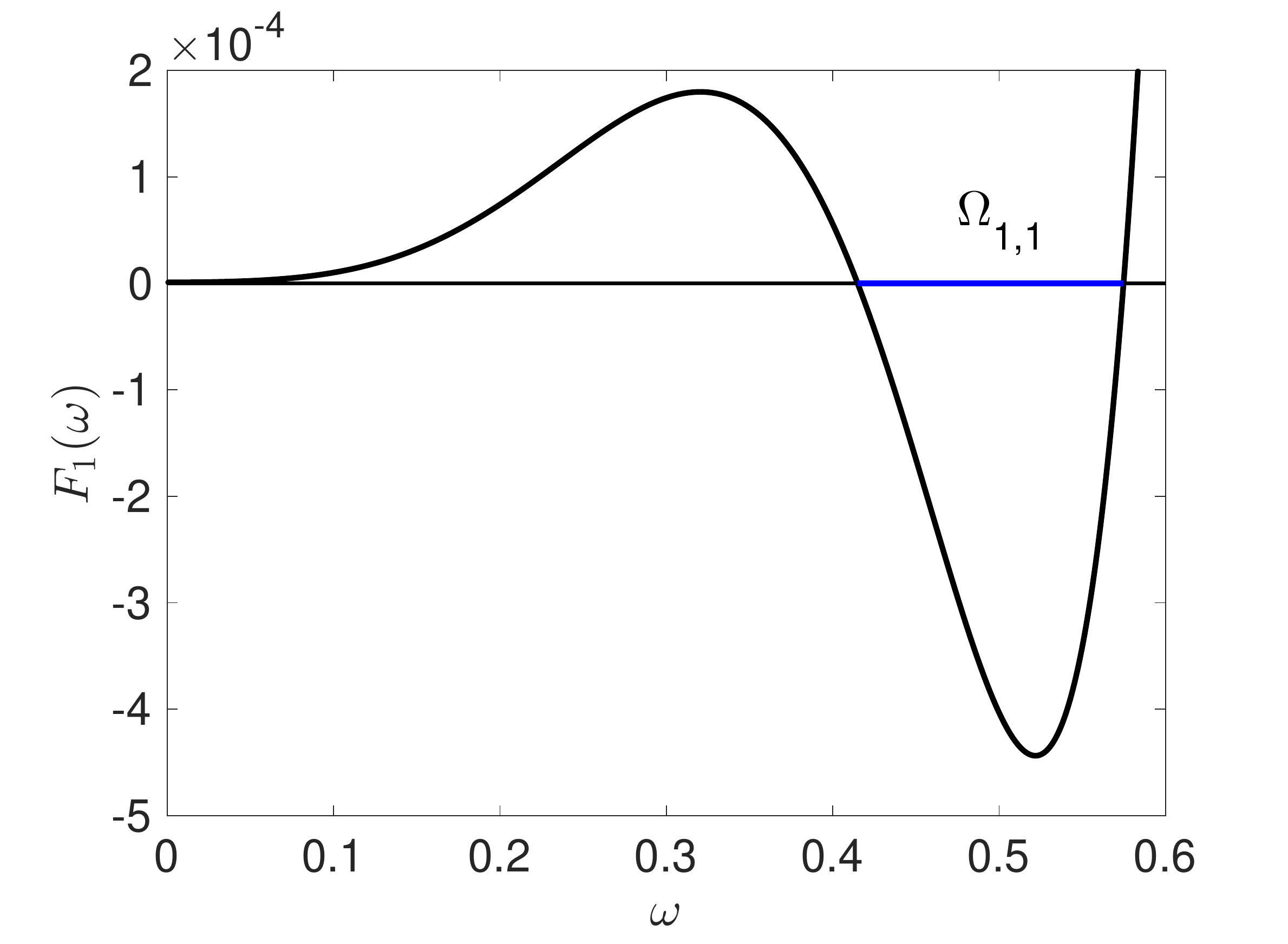}
   b)   \centering
 \includegraphics[width=0.45\textwidth,height=0.30\textwidth]{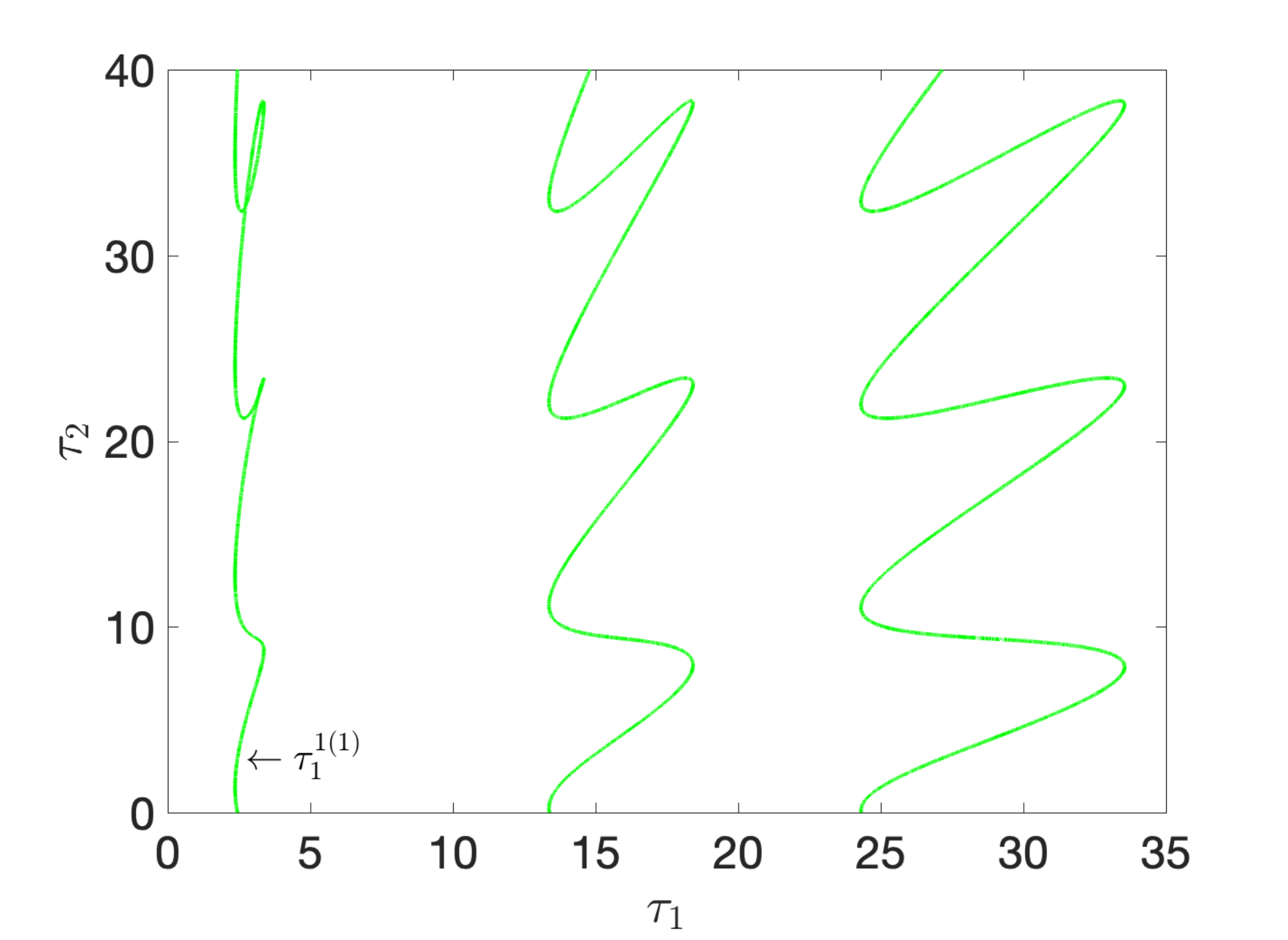}
 \caption{a) The crossing set   $\Omega_{1,1}$.  b) Stability switching curves $\mathcal{T}^1_1$.}
   \label{fig:F1T1}
       \end{figure}

When $n=1$, $F_1(0)>0$, and $F_1(\omega)=0$ has two roots.     The crossing set $\Omega_1=$  $\Omega_{1,1} =[0.4149,0.5747] $, which is shown in Fig. \ref{fig:F1T1} a).   We can get the stability switching curves $\mathcal{T}^1_1$,  which is shown in Fig. \ref{fig:F1T1} b).  Thus all the stability switching curves for $n=1$  are given by $\mathcal{T}^1_1$.

\begin{figure}
     \centering
 a) \includegraphics[width=0.45\textwidth,height=0.30\textwidth]{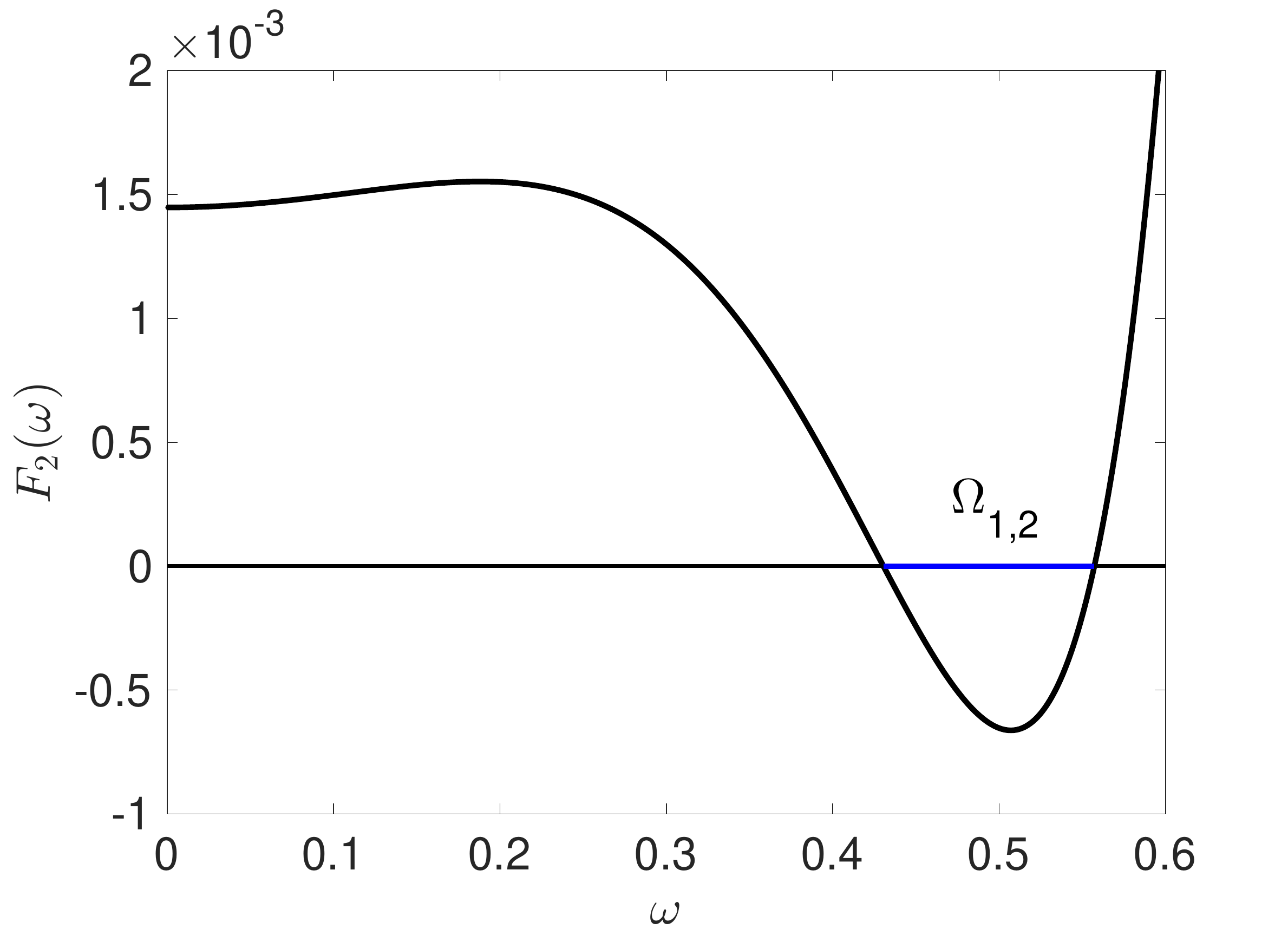}
  b)   \centering
  \includegraphics[width=0.45\textwidth,height=0.30\textwidth]{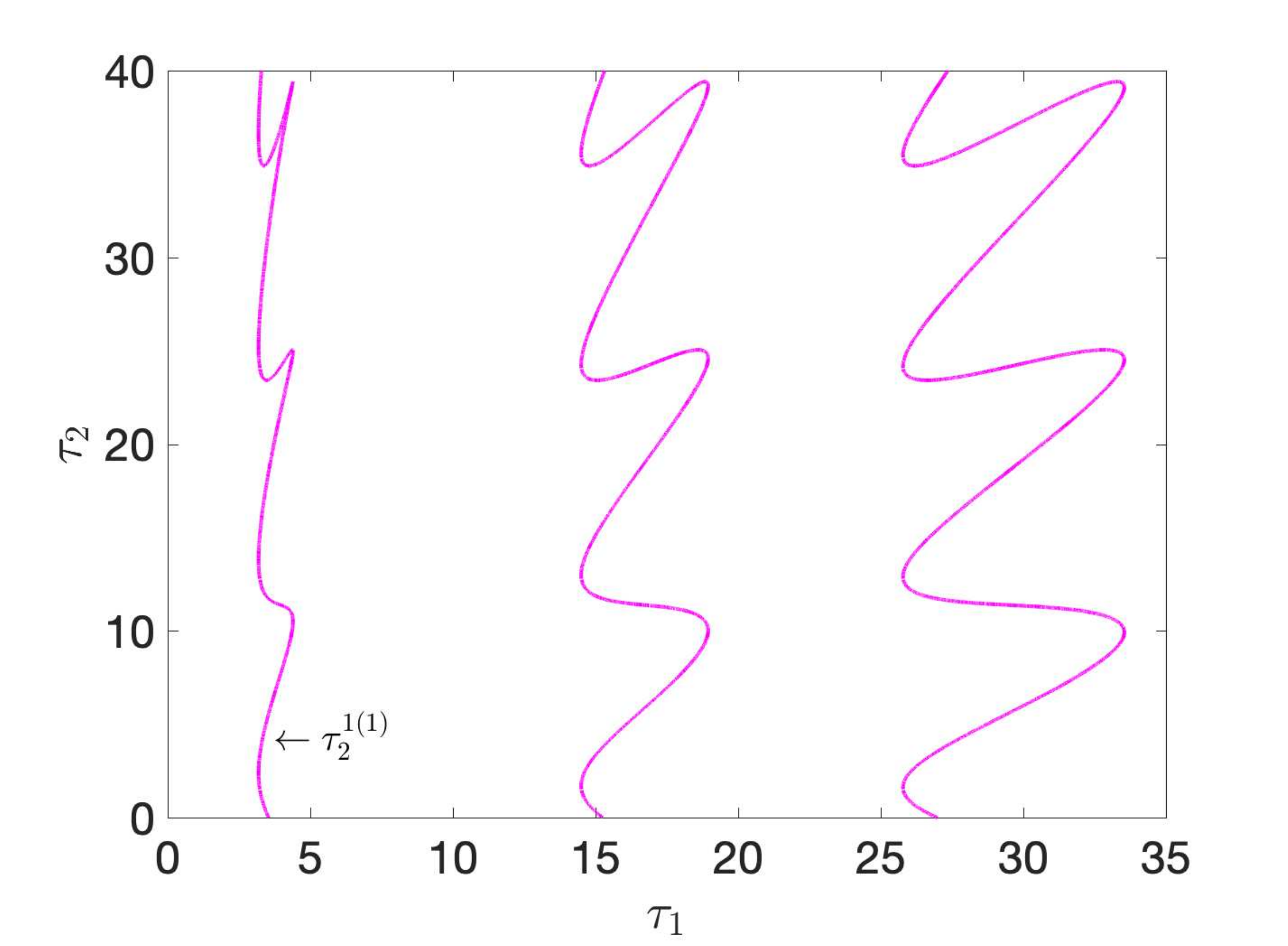}
  \caption{a) The crossing set   $\Omega_{1,2} $. b) Stability switching curves  $\mathcal{T}^1_2$. }
           \label{fig:F2T2}
        \end{figure}

When $n=2$, $F_2(0)>0$, and $F_2(\omega)=0$ has two roots,  and the crossing set is  $\Omega_{1,2}=[0.4304,0.5572]$, which is shown in Fig. \ref{fig:F2T2} a). The stability switching curves $\mathcal{T}^1_2$ corresponding to $\Omega_{1,2}$ is shown in Fig. \ref{fig:F2T2} b).  All the stability switching curves for $n=2$ are given by $\mathcal{T}_2=\mathcal{T}^1_2$.

\begin{figure}
         \centering
  a) \includegraphics[width=0.45\textwidth]{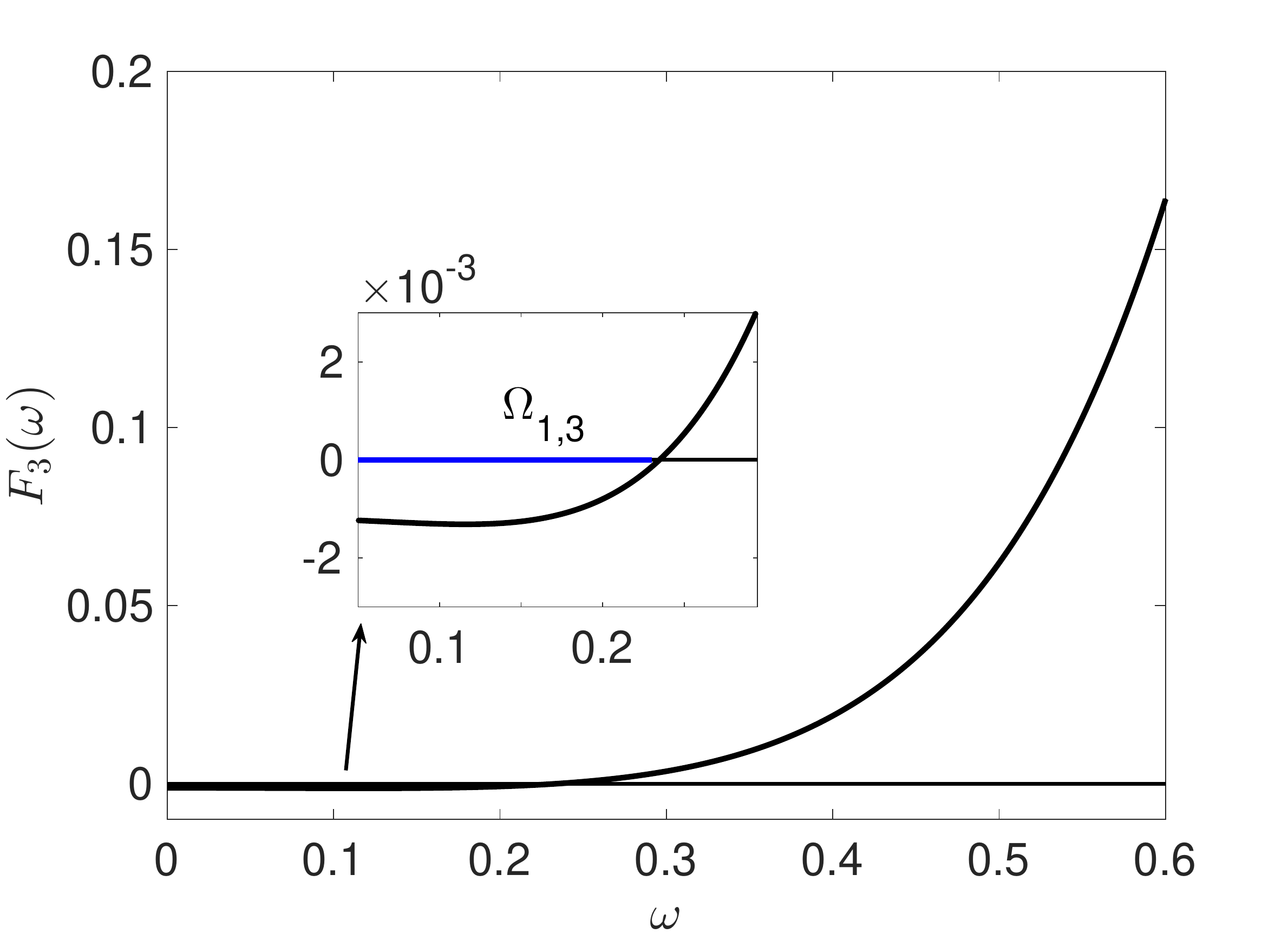}
      b)   \centering
  \includegraphics[width=0.45\textwidth]{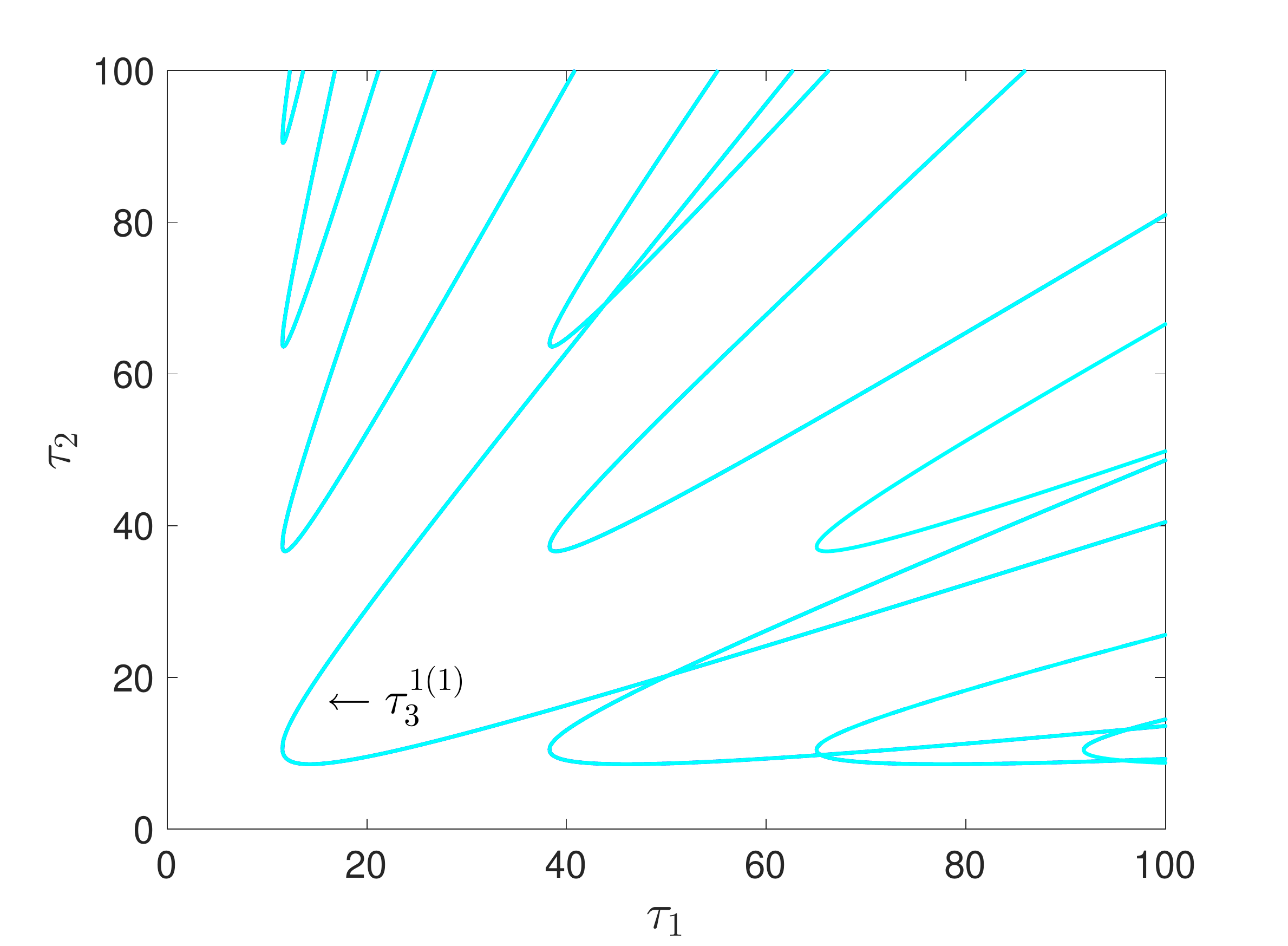}
 \caption{a) The crossing set   $\Omega_{1,3} $. b) Stability switching curves $\mathcal{T}^1_3$. }
        \label{fig:F3T3}
         \end{figure}

When $n=3$, $F_3(0)>0$,  and the crossing set such that  $F_3(\omega)<0$ is  $\Omega_{1,3}=[0,0.2348]$, which is shown in Fig. \ref{fig:F3T3} a). And the stability switching curves $\mathcal{T}^1_3$ corresponding to $\Omega_{1,3}$ is shown in Fig. \ref{fig:F3T3} b).  Thus all the stability switching curves for $n=3$ are given by $\mathcal{T}_3=\mathcal{T}^1_3$.

When $n\geq 4$, $F_n(\omega)>0$ for any $\omega$, thus there are no  stability switching curves on $(\tau_1,\tau_2)$ plane  for  $n\geq 4$.

\section*{Supplementary Material}
In the supplementary material, we give the detailed  calculation process of   normal forms  near double Hopf bifurcation induced by two delays.
\section*{Acknowledgments}
Y. Du is supported by  National Natural Science Foundation of China (grant No.11901369, No.61872227, and No.11971281) and Natural Science Basic Research Plan in Shaanxi Province of China (grant No. 2020JQ-699). B. Niu and J. Wei are supported by  National Natural Science Foundation of China (grant Nos.11701120 and 11771109) .



\end{document}